\documentclass[twocolumn]{autart}    
\let\theoremstyle\relax
\usepackage{amsthm}
\usepackage{amsmath,amsfonts}
\usepackage{array}
\usepackage{textcomp}
\usepackage{stfloats}
\usepackage{url}
\usepackage{verbatim}
\usepackage{cite}
\usepackage{booktabs} 
\usepackage{pifont} 
\usepackage{xcolor}

\usepackage{amsfonts}

\usepackage{color}
\newtheorem{Lemma}{Lemma}
\newtheorem{Remark}{Remark}
\newtheorem{Theorem}{Theorem}
\newtheorem{Definition}{Definition}

\usepackage{algorithm}
\usepackage{bm}
\usepackage{amssymb}
\usepackage{wasysym}
\usepackage{makecell}
\usepackage{cite}
\usepackage{xpatch}

\usepackage{etoolbox}

\usepackage{algpseudocode}
\usepackage{pifont}
\usepackage{multirow}  

\usepackage{graphicx}          

\begin{document}

\begin{frontmatter}


\title{Joint Cooperative and Non-Cooperative Localization in WSNs with Distributed SP-ADMM Algorithms\thanksref{footnoteinfo}} 
\thanks[footnoteinfo]{This paper was not presented at any IFAC 
meeting. Corresponding author Xiaojing Shen. e-mail: shenxj@scu.edu.cn}

\author[China]{Qiaojia Zhu}\ead{zhuqiaojiaa@163.com},   
\author[China]{Xiaojing Shen}\ead{shenxj@scu.edu.cn},                
\author[China]{Haiqi Liu}\ead{haiqiliu0330@163.com)},  
\author[USA]{Pramod K. Varshney}\ead{varshney@syr.edu}  
                                          
\address[China]{The Department of Mathematics, Sichuan University, Chengdu, Sichuan 610064, China}             
\address[USA]{The Department of Electrical Engineering and Computer Science, Syracuse University, Syracuse, NY 13244 USA}        


\begin{abstract}                          
The integration of cooperative and non-cooperative localization is fundamentally important, as these two modes frequently coexist in wireless sensor networks—especially when sensor positions are uncertain and targets are unable to communicate with the network. This paper presents a joint modeling approach that formulates cooperative and non-cooperative localization as a single optimization problem. By processing both tasks jointly, the proposed method eliminates the latency inherent in sequential approaches that perform cooperative localization first, followed by non-cooperative localization. However, this joint formulation introduces complex variable coupling, posing challenges in both modeling and optimization. To address this coupling, we introduce auxiliary variables that enable structural decoupling and facilitate distributed computation. Building on this formulation, we develop the Scaled Proximal Alternating Direction Method of Multipliers for Joint Cooperative and Non-Cooperative Localization (SP-ADMM-JCNL). Leveraging the structured design of the problem, we provide theoretical guarantees that the algorithm generates a sequence converging globally to a Karush–Kuhn–Tucker (KKT) point of the reformulated problem, and further to a critical point of the original non-convex objective function, with a sublinear convergence rate of $O(1/T)$. Experiments demonstrate that SP-ADMM-JCNL achieves accurate and reliable localization performance.
\end{abstract}

\end{frontmatter}

\section{Introduction}
\begin{table*}[t]
\caption{Comparison Between JCNL and SCNL Approaches}
\label{tab:jcvssep}
\centering
\footnotesize
{\renewcommand{\arraystretch}{1.45}
\begin{tabular}{lll}
\toprule
\textbf{Aspect}  & \textbf{SCNL (Separate)} & \textbf{JCNL (Joint)}\\
\midrule
Modeling & Two-stage: \textcolor{black}{sensors → targets} & \textcolor{black}{Joint: simultaneous estimation} \\
Optimization  & Sequential least-squares problems & Single least-squares problem\\
Information Flow  & Used only after sensor localization & \textcolor{black}{Targets} data may aid sensor localization\\
\textcolor{black}{Targets} Estimation Timeliness  & \ding{55} (\textcolor{black}{targets} wait for sensor positions) & \ding{51} (estimated simultaneously with sensors)\\
\bottomrule
\end{tabular}
}
\end{table*}

Wireless Sensor Networks (WSNs), equipped with sensing, processing and communication capabilities, have been widely deployed in various domains such as military operations, environmental monitoring, healthcare, aerospace and industrial systems \cite{george2013shooter, bonnet2000querying, milenkovic2006wireless, cevher2008distributed, petersen2007requirements}. Localization serves as a fundamental enabler of the functionality of WSNs and typically involves two primary tasks: cooperative localization and non-cooperative localization. The former focuses on estimating the positions of sensor sensors using a limited number of anchors with precisely known locations and noisy range measurements between neighboring sensors \cite{wymeersch2009cooperative}. The latter addresses the localization of external targets that are unable to communicate with the network, relying on passive observations collected by sensor sensors and their coordination with nearby sensors.

Both cooperative and non-cooperative localization problems have been extensively studied in the context of WSNs. These tasks typically rely on distance-related observations obtained through techniques such as time of arrival (TOA) \cite{yin2015cooperative, pun2021local}, time difference of arrival (TDOA) \cite{okello2011comparison}, and received signal strength (RSS) \cite{yin2017received, niu2018received}. Recently, optimization-based approaches have gained increasing attention due to their faster inference and theoretical guarantees on convergence. However, the corresponding maximum likelihood formulations usually result in non-convex and non-smooth least squares problems, which are difficult to solve efficiently and reliably \cite{simonetto2014distributed}. To mitigate this, many existing methods adopt convex relaxation techniques that approximate the original problem and enable convergence to a global optimum of the relaxed version \cite{biswas2006semidefinite, songchiyang2009efficient, songchiwang2011semidefinite}. Nonetheless, such solutions may not correspond to critical points of the original non-convex formulation. In addition, most relaxation-based methods are implemented in centralized frameworks, where all measurements must be aggregated at a central sensor. This centralized architecture may increase communication overhead and introduce vulnerability to the malfunction of any single sensor.

To overcome the limitations of centralized approaches, significant efforts have been devoted to developing distributed localization methods \cite{srirangarajan2008distributed, di2016next}. Among them, the Alternating Direction Method of Multipliers (ADMM) has been recognized as a suitable tool for distributed convex optimization and has also shown promise in certain non-convex scenarios. In particular, ADMM has been proven extremely suitable for distributed convex problems due to its decomposability and convergence guarantees \cite{wen2012alternating, magnusson2015distributed,  sun2013fully}, and has also demonstrated effectiveness in tackling structured non-convex formulations under specific conditions. In the context of cooperative localization, a common strategy reformulates the original problem as the problem of two-block non-convex optimization with linear equality constraints. Based on this reformulation, distributed ADMM-based methods have been proposed in \cite{andreani2008augmented, birgin2014practical} to solve the resulting problem. Building on these ideas, \cite{zhang2023distributed} introduced the Scaled Proximal ADMM (SP-ADMM) method—an efficient distributed algorithm with relatively low computational complexity—and established theoretical guarantees for its convergence and performance. More recently, graph neural network based approaches have also been explored for network localization \cite{yan2023attentional, yan2025attentional}. In non-cooperative localization, \cite{zhang2019sensor} addresses a sensor network-based event localization problem using ADMM, under the assumption that the sensor positions are known.

While cooperative and non-cooperative localization have each been extensively studied, they are typically handled separately.
In many practical scenarios, both cooperative and non-cooperative localization problems arise simultaneously. This often occurs in settings such as battlefield or environmental monitoring, where sensor positions are uncertain and \textcolor{black}{\textcolor{black}{targets}} cannot actively communicate with the network. A possible approach is to first perform cooperative localization and then estimate \textcolor{black}{target} positions based on the sensor estimates. For notational convenience, we refer to this sequential strategy as \textbf{SCNL} (Separate Cooperative and Non-Cooperative Localization). In contrast to this staged approach, joint approaches have shown clear advantages in various contexts. For instance, \cite{brambilla2022cooperative} investigates a joint solution to cooperative localization and multitarget tracking, while \cite{mortazavi2017robust, yuan2016cooperative, zheng2009joint} explore the integration of cooperative localization with clock synchronization.
Inspired by these, we propose a joint formulation— \textbf{JCNL} (Joint Cooperative and Non-Cooperative Localization)—that models and solves both localization tasks simultaneously. This integrated approach helps eliminate the \textcolor{black}{target} estimations delay inherent in sequential methods. Since the presence of targets often signals a potential risk or external threat, reducing such delays is particularly important in real-time applications.

In this article, to retain the advantages of the SP-ADMM method in terms of communication efficiency and computational cost, and to address the modeling challenges introduced by the non-convex coupling between sensor and \textcolor{black}{target} variables, we develop a structured formulation of the joint localization problem. By introducing variable decoupling techniques, the original non-convex problem is transformed into a tractable form amenable to distributed optimization. Moreover, we establish the theoretical equivalence between the reformulated and original problems, ensuring that the solution accuracy is not compromised. In addition, our joint formulation naturally reduces to the standard non-cooperative localization problem when all sensors are treated as anchors, implying that the theoretical guarantees presented in this work also apply to that special case. A detailed comparison between JCNL and SCNL is provided in Table~\ref{tab:jcvssep}.

The main contributions of this paper are as follows:
\begin{itemize}
    \item We present a distributed least squares formulation that jointly addresses cooperative (sensor self-localization) and non-cooperative (\textcolor{black}{targets} localization) tasks. Sensor and \textcolor{black}{passive target} positions are modeled together and estimated simultaneously under a unified optimization framework, instead of being decoupled into sequential stages. This joint formulation allows for consistent and concurrent estimation across heterogeneous sensor  types, enhancing the coherence and accuracy of localization in distributed settings where sensor positions are uncertain and targets cannot actively communicate with the network.
    
    \item Inspired by the work in \cite{zhang2023distributed}, we propose distributed SP-ADMM algorithms to address the more complex problem of joint localization involving both sensors and passive targets. The coupling between sensor and target variables introduces new structural challenges that prevent direct decomposition into simple update steps. To overcome this, we redesign several algorithmic components—including variable splitting and local update rules—so that the resulting distributed solver accommodates the jointly structured least squares model, maintains tractability for per-sensor computation, and preserves scalability with respect to network size.
    
    \item We establish that the sequence generated by the proposed algorithm globally converges to a KKT point of the reformulated problem and further to a critical point of the original non-convex objective function. This result is derived based on the problem formulation adopted for JCNL in this work. Compared with our earlier work \cite{zhang2023distributed}, the convergence analysis here deals with a significantly more general problem setting involving joint estimation and non-separable variables, requiring an extended proof framework to handle new coupling and regularity conditions.
   
\end{itemize}
The proposed SP-ADMM-JCNL algorithm was evaluated on synthetic and benchmark networks under various noise models and parameter settings. Results showed that it accurately estimated both sensor and \textcolor{black}{target} positions, with convergence behavior matching the theoretical $\mathcal{O}(1/T)$ rate. Compared to the separate two-stage method, SP-ADMM-JCNL eliminates \textcolor{black}{target} estimation delays and achieves lower overall computation time. By jointly updating sensor and \textcolor{black}{target} positions, it also delivers significantly improved localization performance during the convergence process.

\begin{table*}[!b]
\caption{Node  Types in WSNs for \textcolor{black}{Passive Targets} Localization}
\label{ta2}
\centering
\footnotesize
{\renewcommand{\arraystretch}{1.45}
\begin{tabular}{lcccc}
\toprule
\textbf{Nodes } & \textbf{Type} & \textbf{Position} & \textbf{Communication Enabled?} & \textbf{Number} \\
\midrule
\multirow{2}{*}{Sensor ($\mathcal{N}$)}
& Anchor ($\mathcal{A}$) & Known   & \ding{51} & $\textcolor{black}{A}$ \\
& Agent ($\mathcal{N}\setminus\mathcal{A}$)  & Unknown & \ding{51} & $\textcolor{black}{N-A}$ \\
\textcolor{black}{Target($\mathcal{M}$)} & -- & Unknown & \ding{55} & $\textcolor{black}{M}$ \\
\bottomrule
\end{tabular}
}
\end{table*}

\section*{Synopsis}
Section~\ref{sec:formulation} presents the mathematical formulation of the joint localization problem. Section~\ref{sec:algorithm} introduces the proposed SP-ADMM-JCNL algorithm. Section~\ref{sec:convergence} provides the convergence analysis. Section~\ref{sec:experiments} presents numerical evaluations, including comparisons with baseline methods. Finally, the conlusion is given in Section~\ref{con}.

\section{Problem Formulation}\label{sec:formulation}
\subsection{Notation}
Throughout this paper, column vectors are denoted by boldface lower-case letters (e.g.,~$\mathbf{a}$) and matrices by boldface upper-case letters (e.g.,~$\mathbf{A}$). For an integer $N_i$, the all-ones and all-zeros column vectors of length $N_i$ are denoted by $\mathbf{1}_{N_i}$ and $\mathbf{0}_{N_i}$, respectively. The identity and zero matrices of size $N_i \times N_i$ are written as $\mathbf{I}_{N_i}$ and $\mathbf{O}_{N_i}$. The Euclidean norm is denoted by $\|\cdot\|$ and the Kronecker product by $\otimes$. The operator $\operatorname{vec}(\mathbf{x}_i,\, i \in \mathcal{N})$ stacks the vectors $\mathbf{x}_i \in \mathbb{R}^n$ into one long column vector. The operator $\operatorname{Diag}(\mathbf{z})$ returns a diagonal matrix with vector $\mathbf{z}$ along its diagonal. The unit ball in $\mathbb{R}^n$ centered at the origin is defined as $\mathcal{B} := \left\{ \mathbf{x} \in \mathbb{R}^n \mid \| \mathbf{x} \| \leq 1 \right\}$, and the Cartesian product of $N_i$ such unit balls is denoted by $\mathcal{B}^{N_i}$.

The projection operator onto a convex set $\mathcal{B}^{N_i}$ is defined by
\[
\text{proj}_{\mathcal{B}^{N_i}}(\mathbf{w}_i^t) := \arg\min_{\mathbf{w}_i \in \mathcal{B}^{N_i}} \frac{1}{2} \| \mathbf{w}_i - \mathbf{w}_i^t \|^2.
\]
Let $f : \mathcal{C} \to (-\infty, +\infty]$ be a proper closed convex function and $\mathbf{W}$ a positive semidefinite matrix. The scaled proximal operator is
\[
\text{prox}_f^{\mathbf{W}}(\mathbf{z}) := \arg\min_{\mathbf{v} \in \mathcal{C}} f(\mathbf{v}) + \frac{1}{2} \| \mathbf{v} - \mathbf{z} \|_{\mathbf{W}}^2,
\]
where $\| \mathbf{z} \|_{\mathbf{W}}^2 := \langle \mathbf{z}, \mathbf{Wz} \rangle$ denotes the $\mathbf{W}$-induced norm for $\mathbf{z} $.
We also use $\delta_{\mathcal{C}}(\mathbf{x})$ to denote the indicator function of a closed convex set $\mathcal{C}$, defined as
\[
\delta_{\mathcal{C}}(\mathbf{x}) := 
\begin{cases}
0, & \text{if } \mathbf{x} \in \mathcal{C}, \\
+\infty, & \text{otherwise}.
\end{cases}
\]
The set $\mathcal{C}$ typically represents a feasible region or constraint set imposed on the variables in the optimization problem. Finally, we define with $blockdiag(\cdot)$ the block diagonal matrix whose diagonal contains the input blocks matrices.

\subsection{Problem Statement}
We consider a sensor network for \textcolor{black}{localizing multiple passive targets}, represented by an undirected connected graph \( G = (\mathcal{N}, \mathcal{E}) \), where the network topology is assumed to be known. The sensor set  \( \mathcal{N} = \{1, 2, \dots, N\} \) consists of \( N \) sensors, \textcolor{black}{some of
which are anchors with known true positions collected in the
set $\mathcal{A} = \{\mathbf{a}_{N-A+1}, \dots, \mathbf{a}_N\} \subset \mathcal{N} $.} The remaining sensors are unknown-position sensors that require localization. The true position of sensor  \( i \) is denoted as \( \mathbf{x}_i \in \mathbb{R}^n \), where \( n \) represents the spatial dimension. \textcolor{black}{Additionally, we assume that there are $M$ passive targets, indexed by the set \( \mathcal{M} = \{1, 2, \dots, M\} \). The true position of target $k\in\mathcal{M}$ is denoted by $\mathbf{y}_k\in \mathbb{R}^n$.} For clarity, Table~\ref{ta2} summarizes the types of sensors in WSNs considered for \textcolor{black}{targets} localization.
We define the concatenated position vector as 
\textcolor{black}{
\(\mathbf{x}=[\mathbf{x}_1^T,\mathbf{x}_2^T,\dots,\mathbf{x}_N^T]^T,
\mathbf{y}=[\mathbf{y}_1^T,\mathbf{y}_2^T,\dots,\mathbf{y}_M^T]^T.\)}
For each sensor $ i \in \mathcal{N} $, we define its neighborhood set as $\mathcal{N}_i = \{j \mid (i, j) \in \mathcal{E} \}\subseteq \mathcal{N}$ with cardinality $ N_i = |\mathcal{N}_i| $, where a communication link between neighboring sensors $ i $ and $ j $ enables information exchange; \textcolor{black}{similarly, we define the set of targets within line-of-sight for sensor $ i $ as $ \mathcal{M}_i \subseteq \mathcal{M} $ with $ M_i = |\mathcal{M}_i| $, allowing sensor $ i $ to obtain noisy range measurements to any target $ \textcolor{black}{k \in \mathcal{M}}_i $ by assuming that targets emit distinguishable signals for correct measurement association.}

\begin{figure}
\begin{center}
\includegraphics[height=6.6cm]{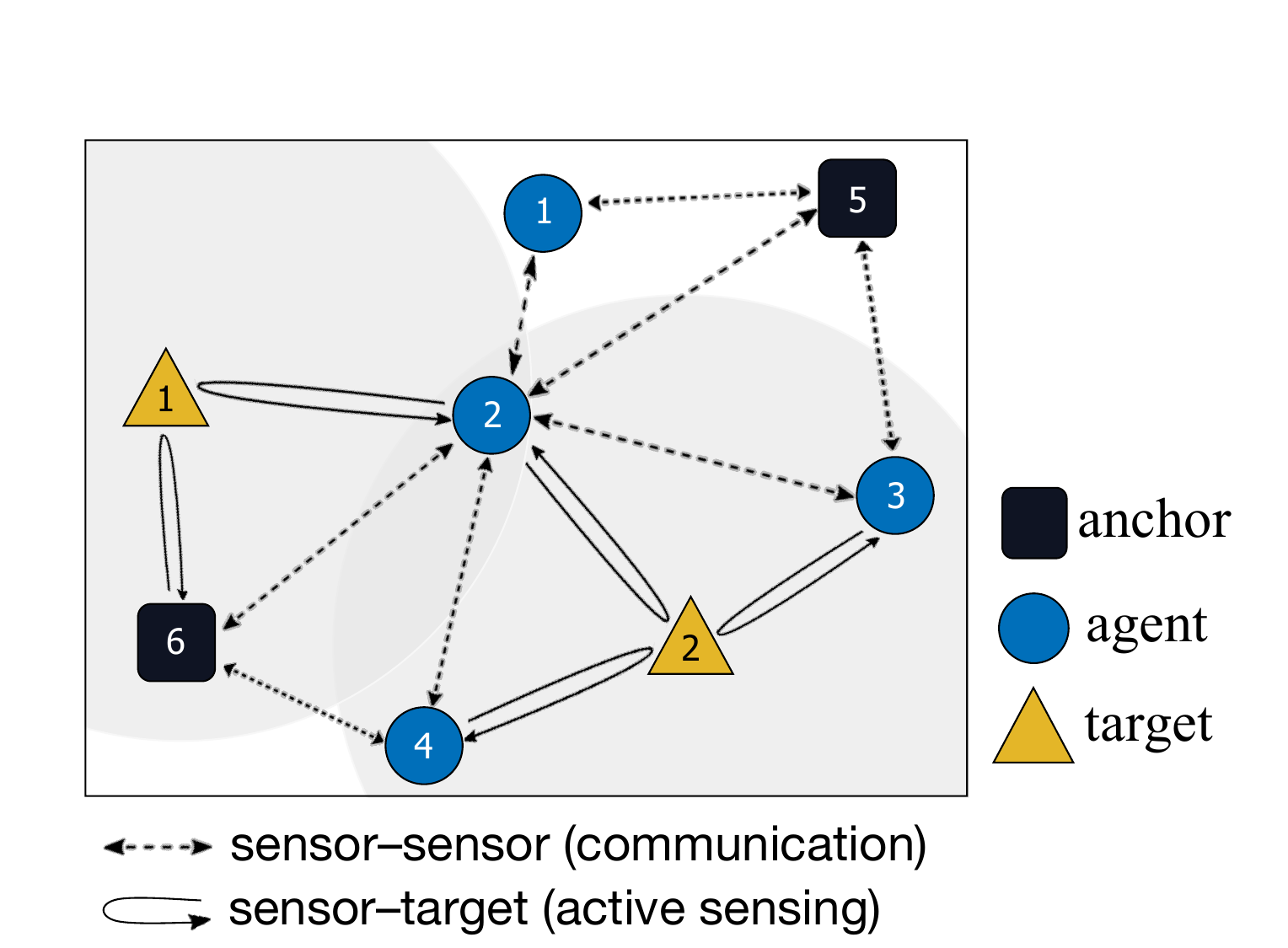}    
\caption{\textcolor{black}{Exemplary illustration of the considered scenario.}}  
\label{shili}                                 
\end{center}                                 
\end{figure}

\textcolor{black}{In the considered scenario, two types of range measurements are available: sensor–sensor measurements obtained via communication links between neighboring sensors, and sensor–target measurements acquired through active sensing (e.g., radar or sonar) for targets within line-of-sight. Targets are passive and cannot perform ranging; only sensors can measure distances to targets. An exemplary illustration of this scenario is provided in Fig.~\ref{shili}. The noisy distance measurements between sensor \( i \) and target \( k\in\mathcal{M}_i \), as well as between sensor \( i \) and its neighboring sensor \( j \in \mathcal{N}_i \), are defined as follows:
\begin{subequations}\label{1}
	\begin{align}
		r_{i, \textcolor{black}{k}} &= \| \mathbf{x}_i - \mathbf{y}_k \| + \tau_{i, \textcolor{black}{k}}, \quad i \in \mathcal{N},\; \textcolor{black}{k \in \mathcal{M}}_i, \\
		d_{i,j} &= \|\mathbf{x}_i - \mathbf{x}_j\| + \omega_{i,j}, \quad i \in \mathcal{N}, \, j \in \mathcal{N}_i.
	\end{align}
\end{subequations}
Here, \textcolor{black}{\( \tau_{i, \textcolor{black}{k}} \)} and \( \omega_{i,j} \) are zero-mean, independent Gaussian-distributed measurement noise terms. We assume that \( d_{i,j} = d_{j,i} \). }


The objective of this study is to estimate the positions of both the unknown sensors and the \textcolor{black}{multiple targets} based on the given measurement model. Using the notation defined above, the maximum likelihood estimator can be obtained by solving the following non-convex constrained optimization problem:
\begin{subequations}\label{yuanwenti}
	\begin{align}
		\min_{\mathbf{x}\in\mathbb{R}^{nN}, \mathbf{y}\in\mathbb{R}^{nM}}&\quad \frac{1}{2}\sum_{i\in\mathcal{N}}f_{i}(\mathbf{x},\mathbf{y}) \\
		\text{subject to }&\quad \mathbf{x}_l = \mathbf{a}_l, \quad \forall l \in \mathcal{A},
	\end{align}
\end{subequations}
where $\forall i \in \mathcal{N},$
\begin{align}\label{ss1}
f_{i}(\mathbf{x},\mathbf{y}) = &\sum_{j\in\mathcal{N}_i}(\|\mathbf{x}_i-\mathbf{x}_j\|-d_{i,j})^2 \nonumber \\
&+ \textcolor{black}{\sum_{k\in\mathcal{M}_i}(\|\mathbf{x}_i-\mathbf{y}_k\|-r_{i, \textcolor{black}{k}})^2}.
\end{align}
\textcolor{black}{We note that when $\mathcal{M} = \emptyset$, the problem reduces to the classic cooperative sensor network localization problem; and when $\mathcal{A} = \mathcal{N}$, problem \eqref{yuanwenti} becomes a pure non-cooperative localization task, where only the targets’ positions need to be estimated.}

\begin{Remark}
\textcolor{black}{ Unlike conventional cooperative localization where all nodes are cooperative, our framework distinguishes between active sensors (cooperative) and passive targets (non-cooperative). While sensors refine their own positions through inter-node communication and range updates, targets lack both the sensing and communication capabilities to participate in the network. Consequently, a target's position is not self-determined but is jointly estimated by its observing sensors. This necessitates a decoupled formulation where target localization relies on consensus constraints across multiple sensors (see Section \ref{pr} for detailed discussion). Beyond the standard spatial coupling via pairwise measurements $\|\mathbf{x}_i - \mathbf{x}_j\|$, any sensors $i$ and $j$ observing a common target $k \in (\mathcal{M}_i \cap \mathcal{M}_j)$ become implicitly coupled through the shared target variable $\mathbf{y}_k$. This multi-path architecture enables sensors to pool mutual geometric information, yet it intensifies the coupling relationships, posing challenges for distributed parallel processing.}
\end{Remark}


\subsection{Problem Reformulation}\label{pr}
Our initial step involves formulating an equivalent smooth and constrained reformulation of problem (\ref{yuanwenti}), which serves as the foundation for both the algorithm design and the convergence analysis in this study. Note that the objective function in (\ref{yuanwenti}) can be explicitly rewritten for all $i \in \mathcal{N}$ as follows:
\begin{align}\label{ns}
	f_{i}(\mathbf{x},\mathbf{y})=
	&\sum_{j\in\mathcal{N}_i} 
	\big( 
	\|\mathbf{x}_i-\mathbf{x}_{j}\|^2 
	- 2d_{i,j}\underbrace{\|\mathbf{x}_i-\mathbf{x}_{j}\|}_{nonsmooth} 
	+ d_{i,j}^2 \big)
	 \nonumber\\
	& + \textcolor{black}{\sum_{\textcolor{black}{k \in \mathcal{M}}_i} \big( \|\mathbf{x}_i-\mathbf{y}_k\|^2-2r_{i, \textcolor{black}{k}}\underbrace{\|\mathbf{x}_i-\mathbf{y}_k\|}_{nonsmooth} +r_{i, \textcolor{black}{k}}^2 \big)} .
\end{align}
It is evident that the above loss function contains non-smooth terms. Based on the existing literature, we apply the Cauchy-Schwarz inequality to obtain:
\begin{subequations}
	\begin{align}
		\|\mathbf{x}_i - \mathbf{x}_{j}\| &= \max_{{\mathbf{v}}_{i,j} \in \mathcal{B}} \mathbf{v}_{i,j}^T (\mathbf{x}_i - \mathbf{x}_{j}),\\
		\textcolor{black}{\|\mathbf{x}_i-\mathbf{y}_k\|} &\textcolor{black}{= \max_{{\mathbf{u}_{i, \textcolor{black}{k}} \in \mathcal{B}}} \mathbf{u}_{i, \textcolor{black}{k}}^T (\mathbf{x}_i - \mathbf{y}_k),}
	\end{align}
\end{subequations}
where \(\mathbf{v}_{i,j},\textcolor{black}{\mathbf{u}_{i, \textcolor{black}{k}}} \in \mathbb{R}^n\) are auxiliary variables. Then the objective function in problem (\ref{yuanwenti}a) can be conveniently reformulated as the minimization of a smooth function over a simple convex constraint set, as follows:
\begin{align}
	\min_{\mathbf{x}, \textcolor{black}{\mathbf{y}}, \mathbf{w}} \sum_{i \in \mathcal{N}}& \Bigg\{  \sum_{j \in \mathcal{N}_i} \left( \frac{1}{2} \|\mathbf{x}_i - \mathbf{x}_j\|^2 - d_{i,j} \mathbf{v}_{i,j}^T (\mathbf{x}_i - \mathbf{x}_j) \right) \nonumber \\
	 +& \textcolor{black}{\sum_{\textcolor{black}{k \in \mathcal{M}}_i} \left( \frac{1}{2} \|\mathbf{x}_i - \mathbf{y}_{k}\|^2 - r_{i, \textcolor{black}{k}} \mathbf{u}_{i, \textcolor{black}{k}}^T (\mathbf{x}_i - \mathbf{y}_{k}) \right)} \Bigg\} \\
	\text{subject to} \quad & \mathbf{w}_i \in \mathcal{B}^{N_i + \textcolor{black}{M}}, \quad \forall i \in \mathcal{N}, \nonumber
\end{align}
where \( \mathbf{w} := \operatorname{vec}(\mathbf{w}_i, i \in \mathcal{N}) \). For each \( i \in \mathcal{N} \), we define $\mathbf{w}_i := \operatorname{vec}(\mathbf{w}_{i,l},\; l \in \{1, \dots, N_i+M\}) \in \mathbb{R}^{(N_i+M)n}$, where \textcolor{black}{$\mathbf{w}_{i,l} = \mathbf{v}_{i,l}$} for \( l = 1, \dots, N_i \), and \textcolor{black}{$\mathbf{w}_{i,l} = \mathbf{u}_{i,{l-N_i}}$ for \( l = N_i+1, \dots, N_i+M \). To ensure a uniform computational structure, dummy variables $\mathbf{u}_{i, \textcolor{black}{k}}$ are included for $k \notin \mathcal{M}_i$, which will not affect the objective value.}

Let \( F_i : \mathbb{R}^{(N_i+1) n} \times \mathbb{R}^{M_in} \times \mathbb{R}^{(N_i+M)n} \to (-\infty, +\infty] \) be a properly defined lower semicontinuous function:
\begin{align}\label{pp}
	&F_i(\mathbf{x},\textcolor{black}{\mathbf{y}}, \mathbf{w}_i)\\
    &:=\sum_{ j \in \mathcal{N}_i}\left(\frac{1}{2}\left\|\mathbf{x}_i - \mathbf{x}_{j}\right\|^2 - d_{i,j}\mathbf{v}_{i,j}^T \left(\mathbf{x}_i - \mathbf{x}_j\right)\right)\nonumber\\
	&\quad+\textcolor{black}{\sum_{\textcolor{black}{k \in \mathcal{M}}_i} \left( \frac{1}{2} \left\|\mathbf{x}_i - \mathbf{y}_k\right\|^2 - r_{i, \textcolor{black}{k}} \mathbf{u}_{i, \textcolor{black}{k}}^T \left(\mathbf{x}_i -\mathbf{y}_k\right) \right)} \nonumber\\
    &\quad+ \delta_{\mathcal{B}^{\textcolor{black}{N_i+M}}}(\mathbf{w}_i),\quad \forall i\in\mathcal{N}
\end{align}
Using the notation introduced above, problem (\ref{yuanwenti}) is equivalently formulated as:
\begin{subequations}\label{ss}
	\begin{align}
		\min_{\mathbf{x}, \mathbf{y},\mathbf{w}} &\quad \sum_{i\in\mathcal{N}} F_i(\mathbf{x},\mathbf{y}, \mathbf{w}_i),\\
		\text{subject to }&\quad \mathbf{x}_l = \mathbf{a}_l, \quad \forall l \in \mathcal{A}.
	\end{align}
\end{subequations}

To facilitate a decentralized solution for problem (\ref{ss}), we let each sensor sensor $i \in \mathcal{N}$ maintain its own local estimate of the position of target \textcolor{black}{$\textcolor{black}{k \in \mathcal{M}}$, denoted by \(\mathbf{y}_{i, \textcolor{black}{k}} \in \mathbb{R}^n\).} To ensure that the local estimates across all sensors reach a common consensus on the true target positions without relying on a central controller, we enforce consistency between neighboring sensors in the communication graph:
\begin{align}\label{yueshu}
	\mathbf{y}_{i, \textcolor{black}{k}} = \mathbf{y}_{j,k}, \quad \forall (i,j) \in \mathcal{E}, \forall \textcolor{black}{k \in \mathcal{M}}.
\end{align}

\textcolor{black}{An important feature of this distributed formulation is its ability to handle limited sensing ranges. Specifically, even if target $k$ is not within the line-of-sight (LoS) of sensor $i$ (i.e., $k \notin \mathcal{M}_i$), sensor $i$ still maintains and updates the local variable $\mathbf{y}_{i, \textcolor{black}{k}}$. Through the iterative collaborative exchange prescribed by (\ref{yueshu}), target-related information propagates from sensors with direct observations to those without, eventually enabling every sensor  in the connected network to estimate the positions of all targets. }

\textcolor{black}{For a compact representation, we collect all local estimates at sensor $i$ into a concatenated vector \(\tilde{\mathbf{y}}_{i, \textcolor{black}{k}} := [\mathbf{y}_{i,1}^T, \dots, \mathbf{y}_{i,M}^T]^T \in \mathbb{R}^{Mn}\). Furthermore, let \(\tilde{\mathbf{y}} := [\tilde{\mathbf{y}}_1^T, \dots, \tilde{\mathbf{y}}_N^T]^T \in \mathbb{R}^{NMn}\) represent the collection of target estimates across the entire network. Based on these definitions, }problem (\ref{ss}) can be equivalently reformulated as the following constrained separable optimization problem:
\begin{subequations}\label{dis}
	\begin{align}
	\min_{\mathbf{x},\textcolor{black}{ \tilde{\mathbf{y}}}, \mathbf{w}} \quad&\sum_{i \in \mathcal{N}} \Big\{\sum_{ j \in \mathcal{N}_i}\left(\frac{1}{2}\left\|\mathbf{x}_i - \mathbf{x}_{j}\right\|^2 - d_{i,j}\mathbf{v}_{i,j}^T \left(\mathbf{x}_i - \mathbf{x}_j\right)\right)\nonumber\\
	&+\textcolor{black}{\sum_{\textcolor{black}{k \in \mathcal{M}}_i} \left( \frac{1}{2} \left\|\mathbf{x}_i - \mathbf{y}_{i, \textcolor{black}{k}}\right\|^2 - r_{i, \textcolor{black}{k}} \mathbf{u}_{i, \textcolor{black}{k}}^T \left(\mathbf{x}_i -\mathbf{y}_{i, \textcolor{black}{k}}\right) \right)} \nonumber\\
    &+ \delta_{\mathcal{B}^{\textcolor{black}{N_i+M}}}(\mathbf{w}_i)\Big\}, \\
		\text{subject}&\text{ to} \quad \mathbf{x}_l = \mathbf{a}_l, \quad \forall l \in \mathcal{A}, \\
		&\quad\quad\textcolor{black}{\mathbf{y}_{i, \textcolor{black}{k}}} = \textcolor{black}{\mathbf{y}_{j,k}}, \quad \forall (i,j) \in \mathcal{E}, \textcolor{black}{\forall \textcolor{black}{k \in \mathcal{M}}},
	\end{align}
\end{subequations}

To decouple the variables for the ADMM framework, we introduce auxiliary variables for both sensor positions and \textcolor{black}{multi-target estimates}. Let \( \mathbf{p}_{i,j}^+ \in \mathbb{R}^n \) and \( \mathbf{q}_{i, j, \textcolor{black}{k}}^+ \in \mathbb{R}^n \) be the copies of the position \( \mathbf{x}_j \) and the target estimate \( \mathbf{y}_{j,k} \) maintained by neighbor \( j \), respectively. Similarly, let \( \mathbf{p}_{i,j}^- \in \mathbb{R}^n \) and  \textcolor{black}{\( \mathbf{q}_{i, j, \textcolor{black}{k}}^- \in \mathbb{R}^n \) }be the copies of sensor $i$'s own position and target estimates assigned to neighbor $j$:
\begin{subequations}\label{fb}
    \begin{align}
        \mathbf{p}_{i,j}^+ :&= \mathbf{x}_j, \quad \textcolor{black}{\mathbf{q}_{i, j, \textcolor{black}{k}}^+ : = \mathbf{y}_{j,k}}, \quad j \in \mathcal{N}_i, \textcolor{black}{\textcolor{black}{k \in \mathcal{M}}}, \\
        \mathbf{p}_{i,j}^- :&= \mathbf{x}_i, \quad \textcolor{black}{\mathbf{q}_{i, j, \textcolor{black}{k}}^- : = \mathbf{y}_{i, \textcolor{black}{k}}}, \quad j \in \mathcal{N}_i,\textcolor{black}{ \textcolor{black}{k \in \mathcal{M}}}.
    \end{align}
\end{subequations}

By aggregating all variables associated with sensor $i$, we define the high-dimensional local variable \( \mathbf{z}_i \) as:
\textcolor{black}{\begin{equation}\label{fz}
    \mathbf{z}_i := \begin{bmatrix} \mathbf{z}^{\text{sen}}_i \\ \mathbf{z}^{\text{tar}}_i \end{bmatrix} \in \mathbb{R}^{(2N_i+1)(1+M)n}, \quad i \in \mathcal{N}.
\end{equation}}
 \textcolor{black}{Here, \( \mathbf{z}^{\text{sen}}_i := [ (\mathbf{x}_i)^T, (\tilde{\mathbf{p}}_i^-)^T, (\tilde{\mathbf{p}}_i^+)^T ]^T\in \mathbb{R}^{(2N_i+1)n}\) , \( \mathbf{z}^{\text{tar}}_i := [  (\mathbf{y}_{i,1})^T,(\tilde{\mathbf{q}}_{i,1}^-)^T, (\tilde{\mathbf{q}}_{i,1}^+)^T,...,(\mathbf{y}_{i,M})^T,(\tilde{\mathbf{q}}_{i,M}^-)^T, \)\\\((\tilde{\mathbf{q}}_{i,M}^+)^T ]^T\in \mathbb{R}^{(2N_i+1)Mn} \), where \( \mathbf{p}_i^- := \operatorname{vec} (\mathbf{p}_{i,j}^-, j \in \mathcal{N}_i) \in \mathbb{R}^{n N_i} \) and \(\forall k\in\mathcal{M}, \mathbf{q}_{i, \textcolor{black}{k}}^- = \operatorname{vec}(\mathbf{q}_{i, j, \textcolor{black}{k}}^-, j \in \mathcal{N}_i) \)  represent the collections of \( N_i \)values of the position \( \mathbf{x}_i \) and the \textcolor{black}{target} position estimate \( \mathbf{y}_{i, \textcolor{black}{k}} \), respectively. Similarly, \( \mathbf{p}_i^+ := \operatorname{vec} (\mathbf{p}_{i,j}^+, j \in \mathcal{N}_i) \in \mathbb{R}^{n N_i} \) and \( \mathbf{q}_{i, \textcolor{black}{k}}^+ := \operatorname{vec} (\mathbf{q}_{i, j, \textcolor{black}{k}}^+, j \in \mathcal{N}_i) \in \mathbb{R}^{n N_i} \) represent the collections of allvalues of \( \mathbf{x}_j \) and \( \mathbf{y}_{j,k} \) received from neighboring sensors \( j \in \mathcal{N}_i \), respectively. }

Using these notations, Equation (\ref{pp}) can be rewritten in a sensor-separable form as follows:
\begin{align}\label{2}
&\sum_{i\in\mathcal{N}}F_i(\mathbf{x}_i,\mathbf{y}_{i}, \mathbf{w}_i) \\=&\sum_{i\in\mathcal{N}}\bigg[\sum_{ j \in \mathcal{N}_i}\left(\frac{1}{2}\left\|\mathbf{x}_i - \mathbf{p}_{i,j}^+\right\|^2 - d_{i,j}\mathbf{v}_{i,j}^T \left(\mathbf{x}_i - \mathbf{p}_{i,j}^+\right)\right)\nonumber\\
    &\quad  +\textcolor{black}{\sum_{\textcolor{black}{k \in \mathcal{M}}_i} \left( \frac{1}{2} \left\|\mathbf{x}_i - \mathbf{y}_{i, \textcolor{black}{k}}\right\|^2 - r_{i, \textcolor{black}{k}} \mathbf{u}_{i, \textcolor{black}{k}}^T \left(\mathbf{x}_i -\mathbf{y}_{i, \textcolor{black}{k}}\right) \right)} \nonumber\\
    &\quad+ \delta_{\mathcal{B}^{\textcolor{black}{N_i+M}}}(\mathbf{w}_i)\bigg].
\end{align}

Problem (\ref{2}) can be rewritten in the following compact form:
\begin{equation}\label{7}
    \min_{\mathbf{z}, \mathbf{w}} \sum_{i \in \mathcal{N}} \left[ \frac{1}{2} \| \mathbf{H}_i \mathbf{z}_i \|^2 - \mathbf{w}_i^T \mathbf{D}_i \mathbf{H}_i \mathbf{z}_i + \delta_{\mathcal{B}^{\textcolor{black}{N_i+M}}} \left( \mathbf{w}_i \right) \right]
\end{equation}
where \( \mathbf{D}_i := \operatorname{Diag} \left(\operatorname{vec} (d_{i,j}, j \in \mathcal{N}_i), \operatorname{vec}( \textcolor{black}{r_{i, \textcolor{black}{k}}, \textcolor{black}{k \in \mathcal{M}}}) \right) \otimes \mathbf{I}_n \) is the measurement matrix. \textcolor{black}{For algebraic consistency and uniform computation, we define $r_{i, \textcolor{black}{k}} = 0$ if target $k \notin \mathcal{M}_i$. This convention ensures that the dimensions of \( \mathbf{D}_i \) and the corresponding mapping matrix \( \mathbf{H}_i \) remain invariant across all nodes. To this end, The mapping matrix \( \mathbf{H}_i \) is defined as:}
\textcolor{black}{\begin{equation}\label{qi}
    \mathbf{H}_i := \begin{bmatrix} \mathbf{H}_{i}^{\text{sen}} \\ \mathbf{H}_{i}^{\text{tar}} \end{bmatrix} \otimes \mathbf{I}_n \in \mathbb{R}^{n(N_i+M) \times n(2N_i+1)(1+M)}.
\end{equation}}
\textcolor{black}{where \( \mathbf{H}_{i}^{\text{sen}} \) captures the relative positions between sensors:
\(
\mathbf{H}_{i}^{\text{sen}} = \begin{bmatrix}
 \mathbf{1}_{N_i} & \mathbf{O}_{N_i} & -\mathbf{I}_{N_i} & \mathbf{0}_{N_i} & \dots & \mathbf{0}_{N_i} \end{bmatrix}\in \mathbb{R}^{N_i \times (2N_i+1)(1+M)},
\)
The target measurement matrix \( \mathbf{H}_{i}^{\text{tar}} = \operatorname{vec}(\mathbf{h}_{i, \textcolor{black}{k}}, \textcolor{black}{k \in \mathcal{M}}) \in \mathbb{R}^{M \times (2N_i+1)(1+M)} \) incorporates the row blocks \( \mathbf{h}_{i, \textcolor{black}{k}} \) for all \( \textcolor{black}{k \in \mathcal{M}} \). Specifically, we define the row block \( \mathbf{h}_{i, \textcolor{black}{k}} \) such that if target \( \textcolor{black}{k \in \mathcal{M}}_i \),
\begin{equation*}\label{hik}
    \mathbf{h}_{i, \textcolor{black}{k}} := \begin{bmatrix} 1 & 0 & \dots & 0 & \underset{\substack{\uparrow \\ \text{Column for } \mathbf{y}_{i, \textcolor{black}{k}}}}{-1} & 0 & \dots & 0 \end{bmatrix};
\end{equation*} otherwise, \( \mathbf{h}_{i, \textcolor{black}{k}} = \mathbf{0}_{(2N_i+1)(1+M)}^T \) for \( k \notin \mathcal{M}_i \). }


The linear constraints \eqref{yueshu} and \eqref{fb} are unified into the compact form
\begin{align}\label{ys1}
    \mathbf{A}_i \mathbf{z}_i = \mathbf{0},
\end{align}  
 \textcolor{black}{where
\begin{equation}\label{ai}
    \mathbf{A}_i := \operatorname{blockdiag} \left( \mathbf{A}_i^{\text{sen}}, \mathbf{A}_{i,1}^{\text{tar}}, \dots, \mathbf{A}_{i,M}^{\text{tar}} \right) \otimes \mathbf{I}_n,
\end{equation}}
 \textcolor{black}{while the first block handles sensor-related constraints:
\begin{equation*}
    \mathbf{A}_i^{\text{sen}} := \begin{bmatrix} \mathbf{1}_{N_i} & -\mathbf{I}_{N_i} & \mathbf{O}_{N_i} \end{bmatrix} \in \mathbb{R}^{N_i \times (1+2N_i)},
\end{equation*}
and the subsequent $M$ blocks enforce consistency for each target $\textcolor{black}{k \in \mathcal{M}}$:
\begin{equation*}
    \mathbf{A}_{i, \textcolor{black}{k}}^{\text{tar}} := \begin{bmatrix} 
        \mathbf{1}_{N_i} & -\mathbf{I}_{N_i} & \mathbf{O}_{N_i} \\ 
        \mathbf{1}_{N_i} & \mathbf{O}_{N_i} & -\mathbf{I}_{N_i} 
    \end{bmatrix} \in \mathbb{R}^{2N_i \times (1+2N_i)}.
\end{equation*}}

Then, based on (\ref{fz}), the anchor sensors in problem (\ref{yuanwenti}) can be represented in the following compact form:
\begin{align}\label{ys2}
	\mathbf{z}&:=\text{vec}(\mathbf{z}_i,i\in\mathcal{N})\in\mathcal{X},
\end{align}
where  $\mathcal{X} := \{ \mathbf{z} \mid \mathbf{E}_i \mathbf{z}_i =\mathbf{a}_i, \, \forall \, i \in \mathcal{A} \},\nonumber$ and $\mathbf{E}_i := \begin{bmatrix} 1, & \mathbf{0}_{ \textcolor{black}{(2N_i+1)M+2N_i}}^T \end{bmatrix} \otimes \mathbf{I}_n \in \mathbb{R}^{n \times  \textcolor{black}{(2N_i+1)(1+M)n}}.$

Moreover, for a pair of connected sensors \( i \) and \( j \), we can derive from equations (\ref{fb}b) that \( \mathbf{p}_{j,i}^- = \mathbf{x}_j \) and \( \mathbf{q}_{j, i, \textcolor{black}{k}}^- = \mathbf{y}_{j,k} \), respectively. Combining these two equations with (\ref{fb}a), we obtain the following constraints for all connected sensors:
\begin{equation}\label{ys3}
	\mathbf{z} \in \mathcal{Y}:= \mathcal{Y}_1\cap \mathcal{Y}_2,
\end{equation}
where
\begin{subequations}\label{12}
	\begin{align}
		\mathcal{Y}_1 &:= \{ \mathbf{z} \mid \mathbf{p}_{i,j}^+ = \mathbf{x}_j = \mathbf{p}_{j,i}^- , \, \forall \, i \in \mathcal{N}, \, j \in \mathcal{N}_i \},\nonumber\\
		\mathcal{Y}_2 &:= \{ \mathbf{z} \mid \mathbf{q}_{i, j, \textcolor{black}{k}}^+ = \mathbf{y}_{j,k} = \mathbf{q}_{j, i, \textcolor{black}{k}}^- , \, \forall k\in\mathcal{M}, i \in \mathcal{N}, \, j \in \mathcal{N}_i \}.\nonumber
	\end{align}
\end{subequations}
By combining equations (\ref{7}), (\ref{ys1}), (\ref{ys2}) and (\ref{ys3}), the non-convex optimization problem (\ref{yuanwenti}) can be reformulated into the following distributed and compact form:
\begin{subequations}\label{jin}
	\begin{align}
		\arg \min_{\mathbf{z},\mathbf{w}} &\sum_{i \in \mathcal{N}} \big\{\underbrace{\frac{1}{2} \| \mathbf{H}_i \mathbf{z}_i \|^2 - \mathbf{w}_i^T \mathbf{D}_i \mathbf{H}_i \mathbf{z}_i}_{G_i(\mathbf{z}_i,\mathbf{w}_i)} + \delta_{\mathcal{B}^{N_i+ \textcolor{black}{M}}} \left( \mathbf{w}_i \right) \big\} \\
		\text{subject to } & \mathbf{A}_i \mathbf{z}_i = \mathbf{0}, \quad \forall i \in \mathcal{N},\\
		&\mathbf{z} \in \mathcal{X}, \, \mathbf{z} \in \mathcal{Y}.
	\end{align}
\end{subequations}

\section{The Proposed SP-ADMM-JCNL Algorithm}\label{sec:algorithm}
The above reformulation introduces separability and structured constraints into the problem, transforming the original joint localization task into an optimization form suitable for the SP-ADMM method. Based on this reformulation, we develop an efficient distributed algorithm tailored to the joint localization setting. We first present the augmented Lagrangian function corresponding to problem (\ref{jin}) as follows:
\[
\mathcal{L}(\mathbf{z}, \mathbf{v}, \boldsymbol{\lambda}) = \sum_{i \in \mathcal{N}} \mathcal{L}_i \left( \mathbf{z}_i, \mathbf{w}_i, \boldsymbol{\lambda}_i \right),
\]
where
\begin{align}\label{lag}
	\mathcal{L}_i \left( \mathbf{z}_i, \mathbf{w}_i, \boldsymbol{\lambda}_i \right) :=& G_i \left( \mathbf{z}_i, \mathbf{w}_i \right) + \delta_{\mathcal{B}^{\textcolor{black}{N_i+M}}} \left( \mathbf{w}_i \right) + \left\langle \boldsymbol{\lambda}_i, \mathbf{A}_i \mathbf{z}_i \right\rangle \nonumber\\
    &+ \frac{c}{2} \| \mathbf{A}_i \mathbf{z}_i \|^2, \quad i \in \mathcal{N},
\end{align} 
In this formulation, \( c > 0 \) is a penalty parameter and \textcolor{black}{\(\boldsymbol{\lambda}_i :=[\boldsymbol{\lambda}_i^\text{sen}\quad\boldsymbol{\lambda}_i^\text{tar}]\in \mathbb{R}^{(1+2M)N_i n} \)}represents the corresponding Lagrange multipliers, Specifically,\textcolor{black}{ \(\boldsymbol{\lambda}_i^\text{sen} = \text{vec} \left( \lambda_{i,j}^0, j \in N_i \right) \otimes \mathbf{I}^n\)}corresponds to the multipliers for the sensing constraints, while \(\boldsymbol{\lambda}_i^\text{tar} = \text{vec} \left( \lambda_{i,j}^k, j \in N_i,k=1,2,...,2M \right) \otimes \mathbf{I}^n\) corresponds to those for the target constraints. The iterative update steps for the variables are obtained by applying the following SP-ADMM algorithm :
\begin{subequations}\label{gengxin1}
	\begin{align}
		&\mathbf{z}^{t+1} = \arg \min_{\mathbf{z} \in \mathcal{X}\cap \mathcal{Y}} \sum_{i \in \mathcal{N}} \mathcal{L}_i \left( \mathbf{z}_i, \mathbf{w}_i^t, \boldsymbol{\lambda}_i^t \right) + \frac{c}{2} \| \mathbf{z}_i - \mathbf{z}_i^{t} \|_{\mathbf{B}_i^T \mathbf{B}_i}^2,  \\
		&\mathbf{w}^{t+1} = \arg \min_{\mathbf{w}} \sum_{i \in \mathcal{N}} \mathcal{L}_i \left( \mathbf{z}_i^{t+1}, \mathbf{w}_i, \boldsymbol{\lambda}_i^t \right) + \frac{\rho}{2} \|\mathbf{w}_i - \mathbf{w}_i^t \|^2,  \\
		&\boldsymbol{\lambda}_i^{t+1} = \boldsymbol{\lambda}_i^t + c \mathbf{A}_i \mathbf{z}_i^{t+1}, \quad \forall i \in \mathcal{N}. 
	\end{align}
\end{subequations}
where \( \rho > 0 \) is a penalty parameter and \( t \) denotes the iteration index.

To facilitate distributed implementation, we design the augmented Lagrangian function \( \mathcal{L} \) such that it is separable with respect to each local variable \( \mathbf{z}_i \), and remains convex when other variables are fixed. These structural properties enable efficient decentralized updates across networked sensors and support scalability in implementation.In each subproblem, scaled proximal terms are incorporated to improve numerical stability and support the convergence analysis. This treatment follows the methodology adopted in our earlier work and ensures that each subproblem remains strongly convex, despite the overall non-convexity of the joint formulation. Specifically, they ensure that each subproblem remains strongly convex and thus admits a unique solution:
\begin{itemize}
	\item \(\mathbf{z}_i \to \mathcal{L}_i \left( \mathbf{z}_i, \mathbf{w}_i^t, \boldsymbol{\lambda}_i^t \right) + \frac{c}{2} \|\mathbf{z}_i - \mathbf{z}_i^{t}\|_{\mathbf{B}_i^T \mathbf{B}_i}^2\) is strongly convex,
	\item \(\mathbf{w}_i \to \mathcal{L}_i \left( \mathbf{z}_i^{t+1}, \mathbf{w}_i, \boldsymbol{\lambda}_i^t \right) + \frac{\rho}{2} \|\mathbf{w}_i - \mathbf{w}_i^t\|^2\) is strongly convex.
\end{itemize}

To see the first point, we reformulate the objective function of problem (\ref{gengxin1}a). Combining (\ref{jin}a) and (\ref{lag}) and rearranging the quadratic term, we have
\begin{align}\label{18}
	&\mathcal{L}_i \left( \mathbf{z}_i, \mathbf{w}_i^t, \boldsymbol{\lambda}_i^t \right) + \frac{c}{2} \|\mathbf{z}_i - \mathbf{z}_i^t\|_{\mathbf{B}_i^T \mathbf{B}_i}^2 \nonumber\\
	=&\frac{1}{2} \|\mathbf{z}_i\|_{\mathbf{U}_i}^2 - \left< \mathbf{H}_i^T \mathbf{D}_i \mathbf{w}_i^t - \mathbf{A}_i^T \boldsymbol{\lambda}_i^t + c \mathbf{B}_i^T \mathbf{B}_i \mathbf{z}_i, \mathbf{z}_i\right> \nonumber\\
	&+ \delta_{\mathcal{B}^{\textcolor{black}{N_i+M}}} \left( \mathbf{w}_i^t \right) + \frac{c}{2} \|\mathbf{z}_i\|_{\mathbf{B}_i^T \mathbf{B}_i}^2,
\end{align}
where
\begin{align}\label{19}
	\mathbf{U}_i = \mathbf{H}_i^T \mathbf{H}_i + c \mathbf{A}_i^T \mathbf{A}_i + c \mathbf{B}_i^T \mathbf{B}_i.
\end{align}
Assuming that the selection matrix \( \mathbf{B}_i^T \mathbf{B}_i \) is chosen such that \( \mathbf{U}_i \) is a positive definite matrix, it follows from (\ref{18}) that the objective function of subproblem (\ref{gengxin1}a) is smooth and strongly convex, thereby ensuring the uniqueness of its solution. By completing the square, we obtain:
\begin{align}\label{20}
	&\mathcal{L}_i \left( \mathbf{z}_i, \mathbf{w}_i^t, \boldsymbol{\lambda}_i^t \right) + \frac{c}{2} \|\mathbf{z}_i - \mathbf{z}_i^{t}\|_{\mathbf{B}_i^T \mathbf{B}_i}^2\nonumber\\
    =& \frac{1}{2} \|\mathbf{z}_i - \tilde{\mathbf{z}}_i^{t+1}\|_{\mathbf{U}_i}^2 - \frac{1}{2} \|\tilde{\mathbf{z}}_i^{t+1}\|_{\mathbf{U}_i}^2 \nonumber\\
    &+ \delta_{\mathcal{B}^{\textcolor{black}{N_i+M}}}(\mathbf{w}_i^{t+1})  + \frac{c}{2} \|\mathbf{z}_i^t\|_{\mathbf{B}_i^T \mathbf{B}_i}^2,
\end{align}
where
\begin{equation}\label{zwan}
	\tilde{\mathbf{z}}_i^{t+1} = \mathbf{U}_i^{-1} \left( \mathbf{H}_i^T \mathbf{D}_i \mathbf{w}_i^{t+1} - \mathbf{A}_i^T \boldsymbol{\lambda}_i^{t+1} + c \mathbf{B}_i^T \mathbf{B}_i \mathbf{z}_i^t \right).
\end{equation}

Substituting (\ref{20}) into (\ref{gengxin1}a) and omitting constant terms, we obtain:
\begin{equation}\label{22}
	\mathbf{z}^{t+1} = \arg \min_{\mathbf{z} \in \mathcal{X}\cap\mathcal{Y}} \sum_{i \in \mathcal{N}} \frac{1}{2}\|\mathbf{z}_i - \tilde{\mathbf{z}}_i^{t+1}\|_{\mathbf{U}_i}^2.
\end{equation}
However, computing $\mathbf{U}_i^{-1}$ can be computationally intensive. To address this, we follow a structured design for $\mathbf{B}_i^T \mathbf{B}_i$ such that $\mathbf{U}_i$ becomes diagonal\cite{zhang2023distributed} :
\begin{equation}\label{23}
	c \mathbf{B}_i^T \mathbf{B}_i = c|\mathbf{A}_i^T \mathbf{A}_i| +| \mathbf{H}_i^T \mathbf{H}_i |,
\end{equation}
where \( |\cdot| \) denotes element-wise absolute value. Based on the definitions of \( \mathbf{A}_i \) and \( \mathbf{H}_i \) in (\ref{ai}) and (\ref{qi}), we can easily obtain the specific form of matrix $c \mathbf{B}_i^T \mathbf{B}_i$. By substituting it into Equation (\ref{19}), we ensure that \( \mathbf{U}_i \) is a positive definite diagonal matrix of the form:
\begin{align}\label{27}
    &\mathbf{U}_i=\nonumber\\
    & 2 \cdot \text{Diag} \Bigg( \Big[\textcolor{black}{ (c+1)N_i+M_i}, c \cdot \mathbf{1}_{N_i}, \mathbf{1}_{N_i}, \textcolor{black}{2cN_i+\mathbb{I}_{\mathcal{M}_i}(1)},\nonumber \\
     &c \cdot \mathbf{1}_{N_i}, c \cdot \mathbf{1}_{N_i},...,\textcolor{black}{2cN_i+\mathbb{I}_{\mathcal{M}_i}(M)}, c \cdot \mathbf{1}_{N_i}, c \cdot \mathbf{1}_{N_i} \Big] \Bigg) \otimes  \mathbf{I}_n,
\end{align}
\textcolor{black}{where $$\forall k\in\mathcal{M},\quad\mathbb{I}_{\mathcal{M}_i}(k) = 
\begin{cases} 
1, & \textcolor{black}{k \in \mathcal{M}}_i, \\ 
0, & k \notin \mathcal{M}_i.
\end{cases}$$}
This diagonal structure significantly reduces the computational burden of updating \( \tilde{\mathbf{z}}_i^{t+1} \).

Furthermore, to provide an explicit expression for \( \mathbf{z}_i^{t+1} \), similar to the definition of \( \mathbf{z}_i \), we partition \textcolor{black}{ \( \tilde{\mathbf{z}}_i := [ (\mathbf{x}_i)^T, (\tilde{\mathbf{p}}_i^-)^T, (\tilde{\mathbf{p}}_i^+)^T, (\mathbf{y}_{i,1})^T,(\tilde{\mathbf{q}}_{i,1}^-)^T, (\tilde{\mathbf{q}}_{i,1}^+)^T,...,(\mathbf{y}_{i,M})^T\)\\\(, (\tilde{\mathbf{q}}_{i,M}^-)^T, (\tilde{\mathbf{q}}_{i,M}^+)^T ]^T\) into $1+M$ components}, where: \(
\tilde{\mathbf{p}}_i^- := \text{vec} \left( \tilde{\mathbf{p}}_{i,j}^-, j \in \mathcal{N}_i \right), \quad \textcolor{black}{\tilde{\mathbf{q}}_{i, k}^- := \text{vec} \left( \tilde{\mathbf{q}}_{i, j, k}^-, j \in \mathcal{N}_i \right)},\)\\\(\tilde{\mathbf{p}}_i^+ := \text{vec} \left( \tilde{\mathbf{p}}_{i,j}^+, j \in \mathcal{N}_i \right),
 \quad \textcolor{black}{\tilde{\mathbf{q}}_{i, k}^+ := \text{vec} \left( \tilde{\mathbf{q}}_{i, j, k}^+, j \in \mathcal{N}_i \right).}\)
\begin{Remark}\label{remark1}
    Given the matrix \( \mathbf{B}_i, i\in\mathcal{N} \) satisfying Equation (\ref{23}), the optimal solution \( \mathbf{z}_i^{t+1} \) for the optimization problem (\ref{22}) has the following closed-form expression:
\begin{subequations}\label{28}
\textit{Sensor-related:}
\begin{align}
\mathbf{x}_i^{t+1} &= \begin{cases} \tilde{\mathbf{x}}_i^{t+1}, & i \notin \mathcal{A} \\ \mathbf{a}_i, & i \in \mathcal{A}, \end{cases} \\
(\mathbf{p}_{i}^-)^{t+1} &= \operatorname{vec} \bigl( c_1 (\tilde{\mathbf{p}}_{i,j}^-)^{t+1} + c_2 (\tilde{\mathbf{p}}_{j,i}^+)^{t+1},\ j\in\mathcal{N}_i \bigr), \\ 
(\mathbf{p}_{i}^+)^{t+1} &= \operatorname{vec} \bigl( c_1 (\tilde{\mathbf{p}}_{i,j}^+)^{t+1} + c_2 (\tilde{\mathbf{p}}_{j,i}^-)^{t+1},\ j\in\mathcal{N}_i \bigr),
\end{align}
where $c_1=\frac{c}{c + 1}, c_2=\frac{1}{c + 1}.$\\
\textit{Target-related:} $\forall \textcolor{black}{k}\in\mathcal{M},$
\begin{align}
\mathbf{y}_{i,\textcolor{black}{k}}^{t+1} &= \tilde{\mathbf{y}}_{i,\textcolor{black}{k}}^{t+1}, \\
(\mathbf{q}_{i,\textcolor{black}{k}}^-)^{t+1} &= \operatorname{vec} \left( \frac{1}{2} \bigl( (\tilde{\mathbf{q}}_{i,j,\textcolor{black}{k}}^-)^{t+1} + (\tilde{\mathbf{q}}_{j,i,\textcolor{black}{k}}^+)^{t+1} \bigr),\ j\in\mathcal{N}_i \right), \\ 
(\mathbf{q}_{i,\textcolor{black}{k}}^+)^{t+1} &= \operatorname{vec} \left( \frac{1}{2} \bigl( (\tilde{\mathbf{q}}_{i,j,\textcolor{black}{k}}^+)^{t+1} + (\tilde{\mathbf{q}}_{j,i,\textcolor{black}{k}}^-)^{t+1} \bigr),\ j\in\mathcal{N}_i \right).
\end{align}
\end{subequations}
\end{Remark} 
By examining equations (\ref{28}a)–(\ref{28}f), it can be observed that the update of \( \mathbf{z}_i^{t+1} \) relies solely on the estimated values of its neighboring sensors \( \tilde{\mathbf{z}}^{t+1}, j \in \mathcal{N}_i \). Consequently, the variable updates can be performed in a distributed manner. The detailed proof is provided in See Appendix \ref{sa1}.
For the optimization problem (\ref{gengxin1}b), omitting the terms that are independent of \( \mathbf{w}_i \), it can be simplified as
\begin{equation}
\begin{aligned}
\arg \min_{\mathbf{w}} \sum_{i \in \mathcal{N}} \bigg( 
    & -\langle \mathbf{H}_i^T \mathbf{D}_i \mathbf{w}_i, \mathbf{z}_i^{t+1} \rangle 
    + \delta_{\mathcal{B}^{\textcolor{black}{N_i+M}}} (\mathbf{w}_i) \\
    & + \frac{\rho}{2} \|\mathbf{w}_i - \mathbf{w}_i^{t} \|^2 
\bigg).
\end{aligned}
\end{equation}
The above objective function is separable with respect to sensor  \( i \) and requires projection onto the unit ball. Therefore, 
\begin{align}\label{30}
    \mathbf{w}_i^{t+1} = \text{proj}_{\mathcal{B}^{\textcolor{black}{N_i+M}}} \left( \mathbf{w}_i^{t} + \frac{1}{\rho} \mathbf{D}_i \mathbf{H}_i \mathbf{z}_i^{t+1} \right).
\end{align}
And its analytical expression is given by
\begin{align}\label{31}
    \mathbf{w}_{i,l}^{t+1} = \frac{\mathbf{w}_{i,l}^{t}}{\max \{ 1, \|\mathbf{w}_{i,l}^{t}\| \}}, \quad l \in\{1, ..., N_i+M\} ,
\end{align}
where $\tilde{\mathbf{w}}_i^{t+1} = \text{vec} \left( \tilde{\mathbf{w}}_{i,l}^{t+1}, l \in\{1, ..., N_i+M\} \right) =\mathbf{w}_i^t + \frac{1}{\rho} \mathbf{D}_i \mathbf{H}_i \mathbf{z}_i^{t+1}.$

The complete distributed joint localization method is detailed in Algorithm~\ref{al1}:
\begin{algorithm}[H]
\caption{{SP-ADMM-JCNL}}\label{al1}
\begin{algorithmic} 
\State \textbf{Input:} Parameters $c$, $\rho$, initial values $\mathbf{z}_i^0$, $\mathbf{w}_i^0$, $\boldsymbol{\lambda}_i^0 = \mathbf{0},$\\\hspace{3.5em}$ \forall i = 1, \dots, N$. Let $\mathbf{U}_i$ be as in \eqref{27}.
\For{$t = 0$ to $T$}
    \ForAll{sensors $i =1, 2, \dots, N$ \textbf{in parallel}}
        \State Update $\tilde{\mathbf{z}}_i^{t+1}$ via (\ref{zwan})
        \If{$i \in \mathcal{A}$}
            \State $\tilde{x}_i^{t+1} = \mathbf{a}_i$.
        \EndIf
    \EndFor
    \State \textbf{Communication:}
    \State \hspace{0.2em} Broadcast $(\tilde{\mathbf{p}}_{i,j}^-)^{t+1}, ( \tilde{\mathbf{q}}_{i, j, \textcolor{black}{k}}^- )^{t+1}, k\in\mathcal{M}$ 
    to $j \in \mathcal{N}_i$;
    \State \hspace{0.2em} Receive $( \tilde{\mathbf{p}}_{j,i}^- )^{t+1}, ( \tilde{\mathbf{q}}_{j,i,k}^+ )^{t+1}, k\in\mathcal{M}$ 
    from $j \in \mathcal{N}_i$.
    \ForAll{sensors  $i = 1, 2, \dots, N$ \textbf{in parallel}}
        \State Update $\mathbf{z}_i^{t+1}$ via (\ref{28}a)--(\ref{28}f);
        \State Update $\mathbf{w}_i^{t+1}$ via (\ref{31});
        \State $\boldsymbol{\lambda}_i^{t+1} = \boldsymbol{\lambda}_i^t + c \mathbf{A}_i \mathbf{z}_i^{t+1}$.
    \EndFor
\EndFor
\end{algorithmic}
\end{algorithm}

\section{Convergence Analysis}\label{sec:convergence}

We analyze the convergence behavior of the proposed SP-ADMM-JCNL algorithm under the more general joint localization setting. Compared to SP-ADMM, which was designed for cooperative localization, the current formulation involves both sensor and target variables, leading to additional coupling and more general constraints. Under standard assumptions, we show that the sequence generated by SP-ADMM-JCNL globally converges to the KKT point of the reformulated problem, which further corresponds to the critical point of the original non-convex objective function. This result is derived based on the problem formulation adopted for JCNL in this work. The convergence analysis extends and refines our earlier work to accommodate the broader setting considered here.

To facilitate the convergence analysis, we first introduce some notation and auxiliary definitions that will be used throughout this section. Let \( \Delta \mathbf{z}_i^{t} := \mathbf{z}_i^{t} - \mathbf{z}_i^{t-1} \), \( \Delta \mathbf{w}_i^{t} := \mathbf{w}_i^{t} - \mathbf{w}_i^{t-1} \), \( \Delta \tilde{\mathbf{z}}_i^{t} := \tilde{\mathbf{z}}_i^{t} - \tilde{\mathbf{z}}_i^{t-1} \) and \( \Delta {\mathbf{\lambda}}_i^{t} := {\mathbf{\lambda}}_i^{t} - {\mathbf{\lambda}}_i^{t-1} ,\quad\forall t\geq1\) denote the updates of the corresponding variables at node \( i \) from iteration \( t-1 \) to \( t \).

With these definitions in place, we proceed to analyze the descent behavior of the augmented Lagrangian. However, it turns out that the function value \( \mathcal{L}(\mathbf{z}, \mathbf{w}, \boldsymbol{\lambda}) \) does not necessarily decrease over iterations proceed, regardless of the penalty parameters \( c \) and \( \rho \). This means \( \mathcal{L} \) alone cannot serve as a suitable descent measure for convergence analysis. To address this limitation, and inspired by the approach in \cite{zhang2023distributed}, we incorporate constraint violation and proximity terms into the augmented Lagrangian and design a new potential function as follows:
\begin{align}\label{43}
	\mathcal{P}^t &= \sum_{i \in \mathcal{N}} \frac{c}{2} \Big[ 
	\kappa_1 \| \mathbf{A}_i \mathbf{z}_i^t \|^2 
	+ \kappa_2 \| \mathbf{A}_i \mathbf{z}_i^{t-1} \|^2 
	+ \frac{\rho}{2c} \|\Delta \mathbf{w}_i^{t}  \|^2 \nonumber\\
	&\quad + (\kappa_1 + \kappa_2) \| \Delta \mathbf{z}_i^{t} \|_{\mathbf{B}_i^T \mathbf{B}_i}^2 \Big] 
	+ \mathcal{L} (\mathbf{z}^t, \mathbf{w}^t, \boldsymbol{\lambda}^t), 
\end{align}
where \( \kappa_1, \kappa_2, c, \rho > 0 \) are positive constants and their specific ranges can be determined by the following lemma.
\begin{Lemma}\label{lem12}
  Suppose that the sequence \( \{ (\mathbf{z}_i^t, \mathbf{w}_i^t, \boldsymbol{\lambda}_i^t) \}_{t \geq 1} \) is generated by Algorithm~\ref{al1} and that \( c \mathbf{B}_i^T \mathbf{B}_i \) follows the form given in Equation(\ref{23}). Then, the following results hold:
\begin{align}\label{44}
	&\mathcal{P}^{t+1} - \mathcal{P}^t  \nonumber\\
    \leq& \sum_{i \in \mathcal{N}} \Bigg[ -\frac{ (\kappa_1 - 1)}{2} (c\|\Delta \mathbf{z}_i^{t+1} \|_{\mathbf{A}_i^T \mathbf{A}_i}^2 
	+  \| \Delta \mathbf{z}_i^{t+1} \| _{\mathbf{H}_i^T \mathbf{H}_i}^2) \nonumber \\
	-&\frac{1}{2} \| \Delta \mathbf{z}_i^{t+1} \| _{\mathbf{U}_i}^2- \frac{c (\kappa_1 - 6 L_2)}{2} \| \Delta \tilde{\mathbf{z}}_i^{t+1} - \Delta \mathbf{z}_i^{t+1}\|_{\mathbf{A}_i^T \mathbf{A}_i}^2 \nonumber \\
	-& \frac{\rho}{4} \|\Delta \mathbf{w}_i^{t+1} \|^2-\frac{c \kappa_1 - 6 L_1}{2c} \|\mathbf{H}_i \Delta \tilde{\mathbf{z}}_i^{t+1} - \mathbf{D}_i \Delta \mathbf{w}_i^{t}\|^2 \nonumber \\
	-& \left( \frac{c \kappa_2 \tilde{\tau}_{\min\quad}}{2(2M+1)N_{\text{sum}}n (c + 1)^2} - \frac{L_2c \kappa_1}{2} \right) \|\tilde{\mathbf{z}}_i^{t+1} - \mathbf{z}_i^t\|^2 \nonumber \\
	-& \frac{c \kappa_1 - 6 (1+c)L_1}{2} \|\Delta \tilde{\mathbf{z}}_i^{t+1} - \Delta \mathbf{z}_i^{t}\|_{\mathbf{B}_i^T \mathbf{B}_i}^2 \nonumber \\
	-& \left( \frac{\rho}{4} - d_{\max}^2 (\kappa_1 + \kappa_2) \right) \|\Delta \mathbf{w}_i^{t}\|^2
	\Bigg],
\end{align}
\textcolor{black}{where $L_1 = N_{\max} + M_{\max} + 1, L_2 = 2N_{\max} + 1$, with $N_{\max} := \max \{ N_i : i \in \mathcal{N} \}$, $M_{\max} := \max \{ M_i : i \in \mathcal{M} \}$.} While \( N_{\text{sum}} := \sum_{i \in \mathcal{N}} N_i \) represents the total number of neighboring nodes in the sensor network, \( d_{\max} := \max\left(\{d_{i,j} \mid i \in \mathcal{N}, j \in \mathcal{N}_i \} \cup \{r_{i, \textcolor{black}{k}} \mid i \in \mathcal{N},\textcolor{black}{k\in\mathcal{M}_i\}}\right) \) denotes the maximum distance measurement, and \(\tilde{\tau}_{\min} := \min \left\{\textcolor{black}{ (c+1)^2[(6M+1)N_i^2+N_i]}, \; i \in \mathcal{N} \right\}.\)
\end{Lemma}
\vspace{-1.1em}
\begin{proof}
See Appendix \ref{sa2}.
   \end{proof}

From Lemma \ref{lem12}, it is evident that as long as the parameters \( \kappa_1, \kappa_2 \) and \( \rho \) are sufficiently large, the right-hand side of the inequality in Equation (\ref{44}) becomes negative. Consequently, the value of the potential function \( \mathcal{P}^t \) decreases with each iteration of the algorithm. Next, we determine the specific ranges of \( \kappa_1, \kappa_2 \) and \( \rho \) that ensure the monotonic decrease of \( \mathcal{P}^t \). First, based on the first five terms on the right-hand side of the inequality in Equation (\ref{44}), a sufficient condition for ensuring the monotonic decrease of \( \mathcal{P}^t \) with respect to \( \kappa_1 \) (for any given \( c > 0 \)) is:
\begin{equation}\label{45}
	\kappa_1 \geq 6 (2N_{\max}+\textcolor{black}{M_{\max}} + 1) \left( 1 + \frac{1}{c} \right). 
\end{equation}
Second, for given values of \( c \) and \( \kappa_1 \), from the sixth term on the right-hand side of the inequality in Equation (\ref{44}), it follows that:
\begin{equation}\label{46}
	\kappa_2 \geq \frac{\textcolor{black}{(2M+1)}N_{\text{sum}} n(c + 1)^2 (2N_{\max} + 1) \kappa_1}{ \tilde{\tau}_{\min\quad}}. 
\end{equation}
Finally, for given values of \( c, \kappa_1, \kappa_2 \), based on the last term on the right-hand side of the inequality in Equation (\ref{44}), it is required that:
\begin{equation}\label{47}
	\rho \geq 4 d_{\max}^2 (\kappa_1 + \kappa_2). 
\end{equation}
\begin{Remark}
Thus, if the parameters satisfy conditions (\ref{45})–(\ref{47}), then the potential function \( \mathcal{P}^t \) decreases monotonically. Unless otherwise specified, all subsequent theorems are derived under the assumption that the parameters are appropriately chosen to satisfy conditions (\ref{45})–(\ref{47}). Moreover, \( \mathcal{P}^t \) is lower-bounded (i.e. there exists a constant \( \underline{\mathcal{P}} > -\infty \) such that \( \mathcal{P}^t \geq \underline{\mathcal{P}} \) for all \( t > 0 \). )This fact will be used in establishing the sublinear convergence rate.
\end{Remark}
Furthermore, we can establish the following result:%

\begin{Theorem}\label{the3}\textbf{(Eventual Consensus)} \\
Consider the sequence $\{ (\mathbf{z}_i^t, \mathbf{w}_i^t, \boldsymbol{\lambda}_i^t) \}_{t \geq 1}, \forall i \in \mathcal{N}$ produced by Algorithm~\ref{al1}. Then:
	\begin{align*}
	  	&\lim_{t \to \infty}  \Delta \mathbf{z}_i^{t}  \to 0, \quad \lim_{t \to \infty} \Delta \boldsymbol{\lambda}_i^{t}  \to 0, \\ 
        &\lim_{t \to \infty} \Delta \mathbf{w}_i^{t} \to 0, \quad \lim_{t \to \infty} (\mathbf{A}_i \mathbf{z}_i^t) \to 0.  
	\end{align*}
\end{Theorem}
\begin{proof}
    See Appendix \ref{sa5}.
\end{proof}

To establish the convergence result, we define the optimality gap function \( \mathcal{G}(\mathbf{z}^t, \mathbf{w}^t, \boldsymbol{\lambda}^t) \) for problem (\ref{jin}) as follows:
\begin{align*}\label{48}
	\mathcal{G}(\mathbf{z}^t, \mathbf{w}^t, \boldsymbol{\lambda}^t) \nonumber:= \sum_{i \in \mathcal{N}} \Bigg[&
	\| \mathbf{z}_i^t - \text{proj}_{\mathcal{X}, \mathcal{Y}} (\mathbf{z}_i^t - (\nabla_{\mathbf{z}_i} G_i (\mathbf{z}_i^t, \mathbf{w}_i^t) \notag \\+ &\mathbf{A}_i^T \boldsymbol{\lambda}_i^t)) \|^2 
	+ \| \mathbf{A}_i \mathbf{z}_i^t \|^2 + \|\Delta \mathbf{w}_i^{t}  \|^2
	\Bigg].
\end{align*}
The following lemma establishes the relationship between the function \( \mathcal{G}(\mathbf{z}^t, \mathbf{w}^t, \boldsymbol{\lambda}^t) \), the Karush-Kuhn-Tucker (KKT) points of problem~(\ref{jin}), and the critical points of the original problem~(\ref{yuanwenti}). This result follows from the formulation of the JCNL problem adopted in this work.

\begin{Lemma}\label{lem15}
 When \( \mathcal{G}(\mathbf{z}^t, \mathbf{w}^t, \boldsymbol{\lambda}^t) = 0 \), the tuple \( (\mathbf{z}^t, \mathbf{w}^t, \boldsymbol{\lambda}^t) \) is a KKT point of problem (\ref{jin}), i.e., it satisfies the following conditions:
\begin{subequations}\label{49}
	\begin{align}
		 \mathbf{z}^t &\in \arg \min_{\mathbf{z} \in \mathcal{Y}, \mathbf{z} \in \mathcal{X}}\sum_{i \in \mathcal{N}} \Big[ G_i (\mathbf{z}_i, \mathbf{w}_i^t) + \langle \boldsymbol{\lambda}_i^t, \mathbf{A}_i \mathbf{z}_i \rangle \Big], \\
		 \mathbf{0} &\in \nabla_{\mathbf{w}_i} G_i (\mathbf{z}_i^t, \mathbf{w}_i^t) + \partial \delta_{\mathcal{B}^{\textcolor{black}{N_i+M}}} (\mathbf{w}_i^t), \quad\forall i\in\mathcal{N} \\
		 \mathbf{0} &= \mathbf{A}_i \mathbf{z}_i^t, \quad\forall i\in\mathcal{N}.
	\end{align}
\end{subequations}
Moreover, $(\mathbf{x}^t, \mathbf{y}_1^t,\mathbf{w}^t)$ is a critical point of problem (\ref{ss}), and $(\mathbf{x}^t, \mathbf{y}_1^t)$ is a critical point of the original problem (\ref{yuanwenti}), where $\mathbf{y}_1^t=\text{vec}(\mathbf{y}_{1,k},k\in\mathcal{M})$.
\end{Lemma}
\begin{proof}
    See Appendix \ref{sa4}.
\end{proof}
\begin{Remark} Since \( \mathbf{z} \in \mathcal{Y} \), if \( (\mathbf{z}^t, \mathbf{w}^t, \boldsymbol{\lambda}^t) \) satisfies condition (\ref{49}c), then it follows that the constraint \( \mathbf{y}_{i, \textcolor{black}{k}} = \mathbf{y}_{j,k}, \forall (i,j) \in \mathcal{E},\textcolor{black}{k\in\mathcal{M}} \) holds. Given that the sensor network is connected, this implies the consistency of the entire sensor network's estimation of the target position, i.e., for any \( i \in \mathcal{N},k\in\mathcal{M} \), we have \( \mathbf{y}_{i, \textcolor{black}{k}}^t = \mathbf{y}_{1,k}^t\). Here, \( \mathbf{y}_1^t \) can be replaced with any \( \mathbf{y}_{i}^t, i \in \mathcal{N} \).
\end{Remark}

Lemma \ref{lem15} indicates that proving the sublinear convergence of Algorithm~\ref{al1} to a critical point of the original optimization problem  (\ref{ss})can be achieved by demonstrating that the sequence \( \{ (\mathbf{z}^t, \mathbf{w}^t, \boldsymbol{\lambda}^t) \}_{t \geq 1} \) sublinearly converges to the critical point of problem (\ref{jin}). 

Combining the following definition and Lemma~\ref{lem17}, we establish the sublinear convergence rate of Algorithm~\ref{al1}.

\begin{Definition}[$\epsilon$-Solution]\label{def:eps_solution}
A tuple \( (\mathbf{z}, \mathbf{w}, \boldsymbol{\lambda}) \) is called an \( \epsilon \)-solution of problem~(\ref{jin}) if it satisfies \( \| \mathbf{A}_i \mathbf{z}_i \| < \epsilon \) and there exist vectors $\boldsymbol{\eta}_1 \in \nabla_{\mathbf{z}_i} G_i(\mathbf{z}_i, \mathbf{w}_i) + \mathbf{A}_i \boldsymbol{\lambda}_i + \partial \delta_{\mathcal{X} \times \mathcal{Y}} (\mathbf{z}_i)$, $\boldsymbol{\eta}_2 \in \nabla_{\mathbf{z}_i} G_i(\mathbf{z}_i, \mathbf{w}_i) + \partial \delta_{\mathcal{B}^{\textcolor{black}{N_i+M}}} (\mathbf{w}_i)$ such that \( \| \boldsymbol{\eta}_1  \| \leq \epsilon \) and \( \| \boldsymbol{\eta}_2 \| \leq \epsilon \).
\end{Definition}

\begin{Lemma}\label{lem17}
If the parameters satisfy conditions (\ref{45})–(\ref{47}), then for any \( t \geq 1 \), there exists a constant \( M > 0 \) and an index \( s \in \{1, 2, \ldots, t-1\} \) such that \( (\mathbf{z}_i^{s+1}, \mathbf{w}_i^{s+1}, \boldsymbol{\lambda}_i^{s+1}) \) is an \( M/\sqrt{t-1} \)-solution. This implies that an \( \epsilon \)-solution can be obtained within at most \( M^2/\epsilon^2 \) iterations.
\end{Lemma}
\begin{proof}
    See Appendix \ref{sa6}.
\end{proof}

Based on the above lemmas, we now summarize the convergence properties of Algorithm~\ref{al1} in the following theorems. The results provide theoretical support for the algorithm’s validity and underpin its empirical performance observed in the experiments. 

\begin{Theorem}\label{the1}
For each node \( i \in \mathcal{N} \), suppose the sequence \( \{ (\mathbf{z}_i^t, \mathbf{w}_i^t, \boldsymbol{\lambda}_i^t) \}_{t \geq 1} \) is generated by Algorithm~\ref{al1} and the parameters are appropriately set to satisfy conditions (\ref{45})–(\ref{47}). Then, the following conclusions hold:
\begin{itemize}
    \item \textbf{Convergence to Stationary Points:}  
    \( \mathcal{G}(\mathbf{z}^t, \mathbf{w}^t, \boldsymbol{\lambda}^t) \to 0 \) and the selected sequence 
    \( \{ (\mathbf{z}_i^t, \mathbf{w}_i^t, \boldsymbol{\lambda}_i^t) \} \) converges to a KKT point of problem~(\ref{jin}).
    
    \item \textbf{Sublinear Convergence Rate:}  
    For any given \( \epsilon_1 > 0 \), suppose that at iteration \( T \), the function value 
    \( \mathcal{G}(\mathbf{z}^t, \mathbf{w}^t, \boldsymbol{\lambda}^t) \) first falls below \( \epsilon_1 \), i.e.
    \[
    T := \arg \min_t \mathcal{G}(\mathbf{z}^t, \mathbf{w}^t, \boldsymbol{\lambda}^t) \leq \epsilon_1.
    \]
    Then there exists a constant \( \epsilon_2 > 0 \) such that 
    \( \epsilon_1 \leq \frac{\epsilon_2}{T-1} \),
    implying that the convergence rate of 
    \( \mathcal{G}(\mathbf{z}^t, \mathbf{w}^t, \boldsymbol{\lambda}^t) \)
    is \( \mathcal{O}(1/T) \).
\end{itemize}
\end{Theorem}
\begin{proof}
    See Appendix \ref{sa7}.
\end{proof}
\begin{Theorem}\label{the2}\textbf{ (Global Convergence).} \\
Let the sequence \( \{ \mathbf{s}^t = (\mathbf{z}^t, \mathbf{w}^t, \boldsymbol{\lambda}^t) \}_{t \geq 1} \) be generated by Algorithm~\ref{al1} and suppose that the parameters satisfy conditions (\ref{45})–(\ref{47}). Then, the following results hold:
\begin{itemize}
	\item The sequence \( \{ \mathbf{s}^t \}_{t \geq 1} \) is bounded and has finite length, meaning that:
	\[
	\sum_{t=1}^\infty \| \mathbf{s}^{t+1} - \mathbf{s}^t \| < +\infty.
	\]
	\item The sequence \( \{ \mathbf{s}^t \}_{t \geq 1} \) converges to a KKT point of problem (\ref{jin}).
\end{itemize}
\end{Theorem}
\begin{proof}
    See Appendix \ref{sa8}.
\end{proof}

\section{Numerical Experiments}\label{sec:experiments} 
We conduct experiments on a two-dimensional localization problem to validate the proposed SP-ADMM-JCNL algorithm. To examine the algorithm’s performance and applicability under different settings, the experimental evaluation is organized into two parts: 
(i) performance assessment on a synthetic sensor network with controllable topology and noise levels, which enables flexible investigation of the algorithm's behavior under different parameter settings and network configurations; and comparison with a two-stage separate localization strategy (described in Algorithm~\ref{alg:scnl}), where sensor and \textcolor{black}{target} positions are estimated in sequential steps.
\textcolor{black}{(ii) analysis of the impact of different parameter configurations on the performance and convergence of JCNL.}

To quantitatively assess the results, we separately report the root mean squared error (RMSE) for sensor self-localization and \textcolor{black}{targets} localization. This distinction allows for a clearer understanding of the respective performance in both cooperative and non-cooperative components of the joint localization task.


\begin{itemize}
    \item The \textbf{sensor localization RMSE} at iteration \( t \) is defined as:
    \begin{equation}
    \text{RMSE}_{\text{sensor}}(t)  = \sqrt{\frac{\sum\limits_{i \in (\mathcal{N}/\mathcal{A})} \| \mathbf{x}_i^t - \mathbf{x}_i \|^2}{N - \textcolor{black}{A}}},
     \end{equation}
    where \( \mathbf{x}_i^t \) denotes the estimated position of sensor \( i \) and \( \mathbf{x}_i \) is its ground-truth position.

    \item The \textbf{\textcolor{black}{targets} localization RMSE} at iteration \( t \) is given by:
    \begin{equation}
    \text{RMSE}_{\text{\textcolor{black}{target}}}(t) =  \sqrt{\frac{\sum\limits_{\textcolor{black}{k \in \mathcal{M}}}\sum\limits_{i \in \mathcal{N}} \| \mathbf{y}_{i, \textcolor{black}{k}}^t - \mathbf{y}_{k} \|^2}{MN}},
    \end{equation}
    where \( \mathbf{y}_{i, \textcolor{black}{k}}^t \) represents the estimated \textcolor{black}{target} position by sensor \( i \) and \( \mathbf{y} \) is the true \textcolor{black}{target} location.
    \end{itemize}
These two error metrics are evaluated separately to characterize the localization performance for both sensors and the \textcolor{black}{target} event.

The following metrics are defined to evaluate convergence:
\begin{align}
S(t) &= \sum\limits_{i \in \mathcal{N}} \| \nabla_{\mathbf{z}_i} F_i (\mathbf{z}_i^t, \mathbf{w}_i^t) + \mathbf{A}_i^T \boldsymbol{\lambda}_i^t \|^2, \\
W(t) &= \sum\limits_{i \in \mathcal{N}} \| \mathbf{w}_i^t -  \mathbf{w}_i^{t-1} \|^2, \\
P(t) &= \sum\limits_{i \in \mathcal{N}} \| \mathbf{A}_i \mathbf{z}_i^t \|^2.
\end{align}
Here, \( S(t) \) captures the stationarity gap, \( W(t) \) reflects the change in iterates (update gap) and \( P(t) \) measures constraint feasibility. These metrics jointly reflect both optimization convergence and constraint satisfaction performance.

\subsection{Synthetic Network Evaluation}
\begin{figure}[htbp]
\centering
\includegraphics[height=5.6cm]{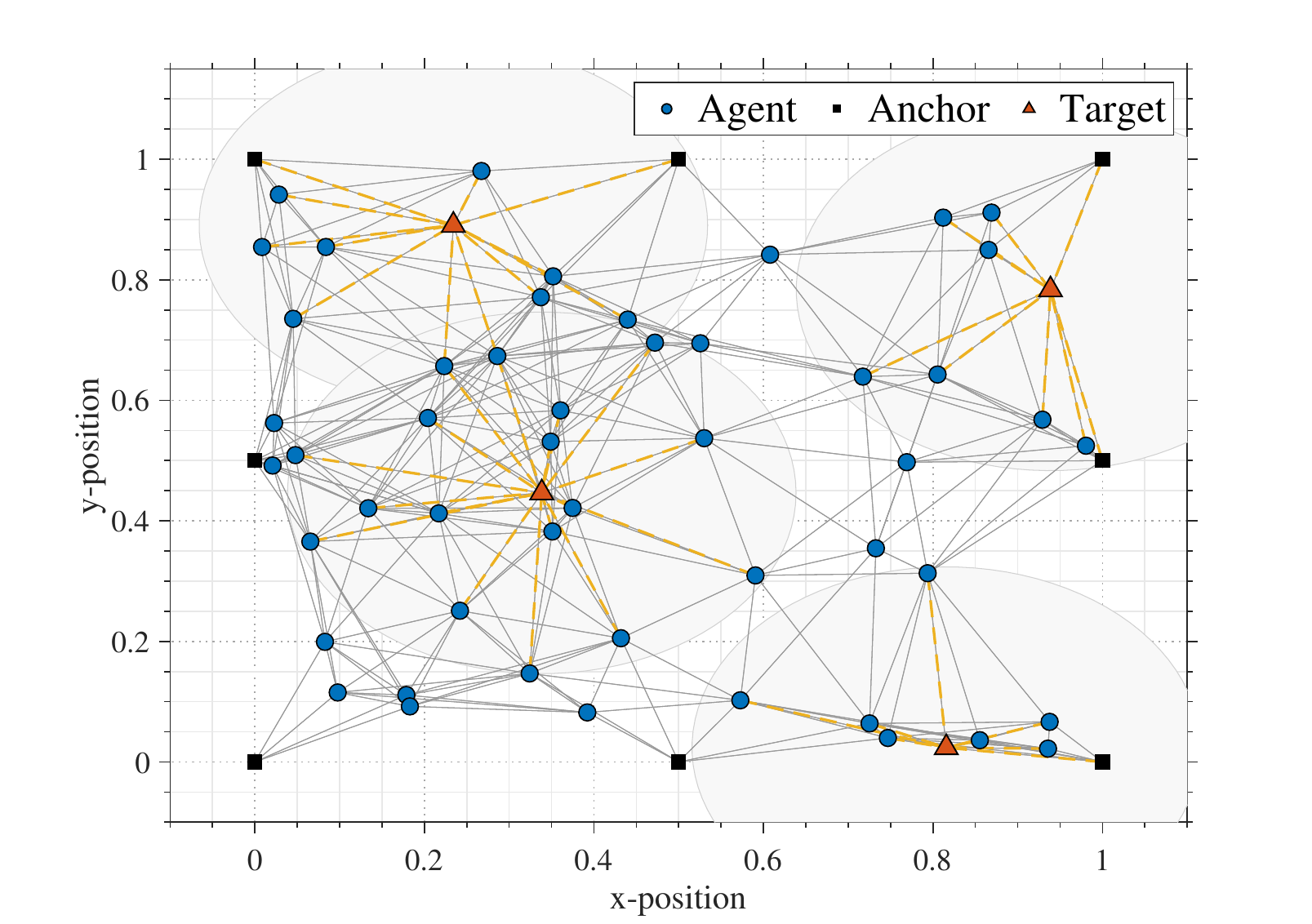}
\caption{Synthetic network layout.}
\label{fig:SN}
\end{figure}

We begin by constructing a synthetic sensor network to evaluate the convergence behavior and localization accuracy of the proposed SP-ADMM-JCNL algorithm in a controlled environment. Specifically, a total of  \textcolor{black}{50 agents}, 8 anchors, and \textcolor{black}{four targets} are randomly deployed within a square region $[0,1] \times [0,1] $. The communication range is set to $C_{\text{range}} = 0.3$, which results in a moderately connected network topology with sparse local neighborhoods. This level of connectivity reflects realistic conditions often encountered in low-power wireless sensor deployments. The corresponding network topology used in our experiments is shown in Fig.~\ref{fig:SN}.

Following \cite{jin2021exploiting}, we consider two types of measurement noise to simulate different sensing conditions. The first type is additive white Gaussian noise (AWGN), where the measurement noise is modeled as a zero-mean Gaussian variable with fixed variance, i.e., $\sigma_{i,j}^2 = \sigma_{\text{add}}^2$. This noise model captures the baseline level of measurement uncertainty commonly assumed in localization literature. The second type is range-dependent Gaussian noise, where the variance increases proportionally with the square of the Euclidean distance between sensor  pairs. Specifically, the noise is modeled as zero-mean Gaussian with $\sigma_{i,j}^2 = \sigma_{\text{add}}^2 \| \mathbf{p}_i - \mathbf{p}_j \|^2$, where $\| \mathbf{p}_i - \mathbf{p}_j \|$ denotes the true distance between sensors $i$ and $j$. This setting reflects more realistic scenarios where measurement errors grow with distance. In this part, we set $\sigma_{\text{add}}^2 = 0.02$. These two noise types are tested separately in order to evaluate the robustness and adaptability of the proposed algorithm under varying noise characteristics.

For the synthetic network experiments, the penalty parameters in the augmented Lagrangian are set to $\rho = c = 0.26$ and the initial dual variable is set as $\mathbf{u}^0 = \mathbf{0}$. Initial estimates are drawn from a uniform distribution: \( \mathbf{z}_0 \sim \text{Unif}(-1, 1)^{(M+1)(4|\mathcal{E} |+ 2N)} \), where $\mathcal{E}$ denotes the set of communication links and $N$ is the number of sensors.

\begin{figure}[htbp]
    \centering

    \begin{minipage}{0.236\textwidth}
        \centering
        \includegraphics[width=\linewidth]{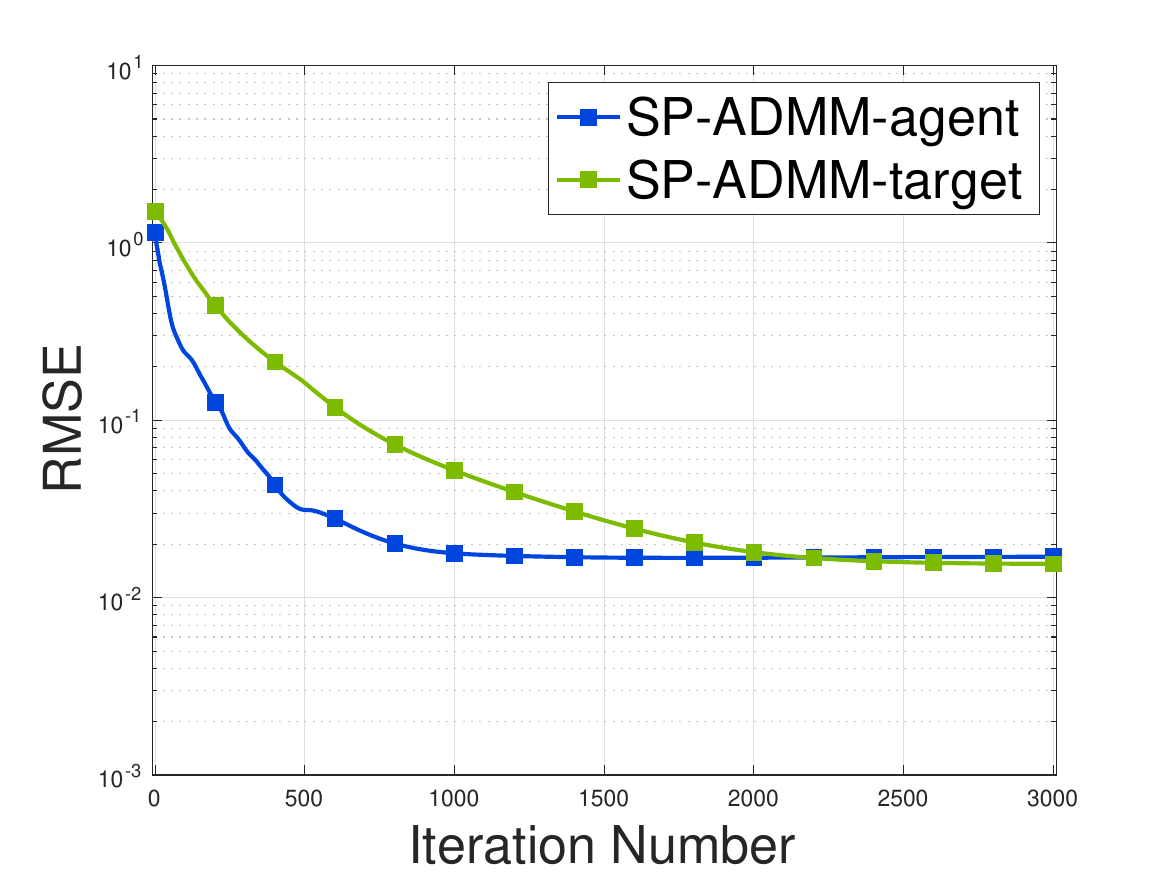}
        \vspace{0.2em}
        \scriptsize (a)AWGN
    \end{minipage}
    \hfill
    \begin{minipage}{0.236\textwidth}
        \centering
        \includegraphics[width=\linewidth]{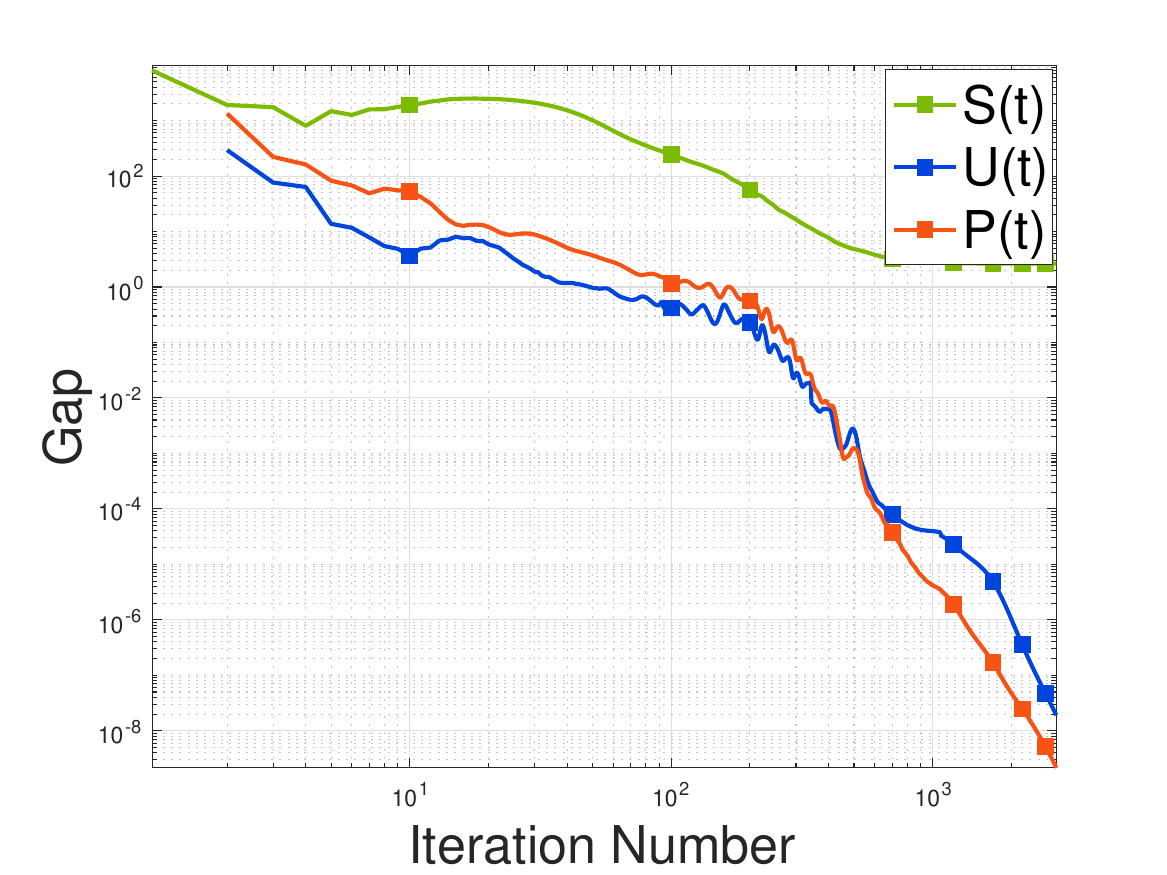}
        \vspace{0.2em}
        \scriptsize (b)AWGN
    \end{minipage}

    \vspace{0.5em}

    \begin{minipage}{0.236\textwidth}
        \centering
        \includegraphics[width=\linewidth]{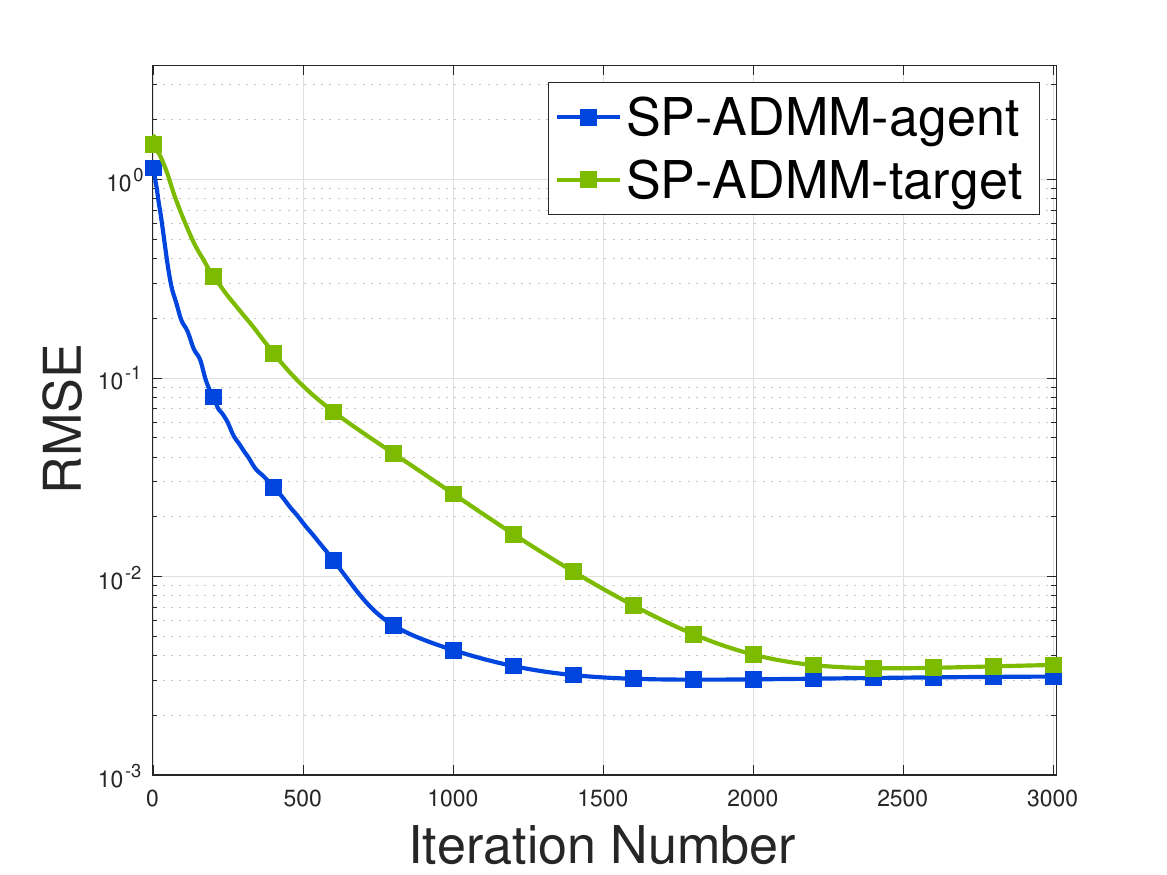}
        \vspace{0.2em}
        \scriptsize (c)Range-dependent noise
    \end{minipage}
    \hfill
    \begin{minipage}{0.236\textwidth}
        \centering
        \includegraphics[width=\linewidth]{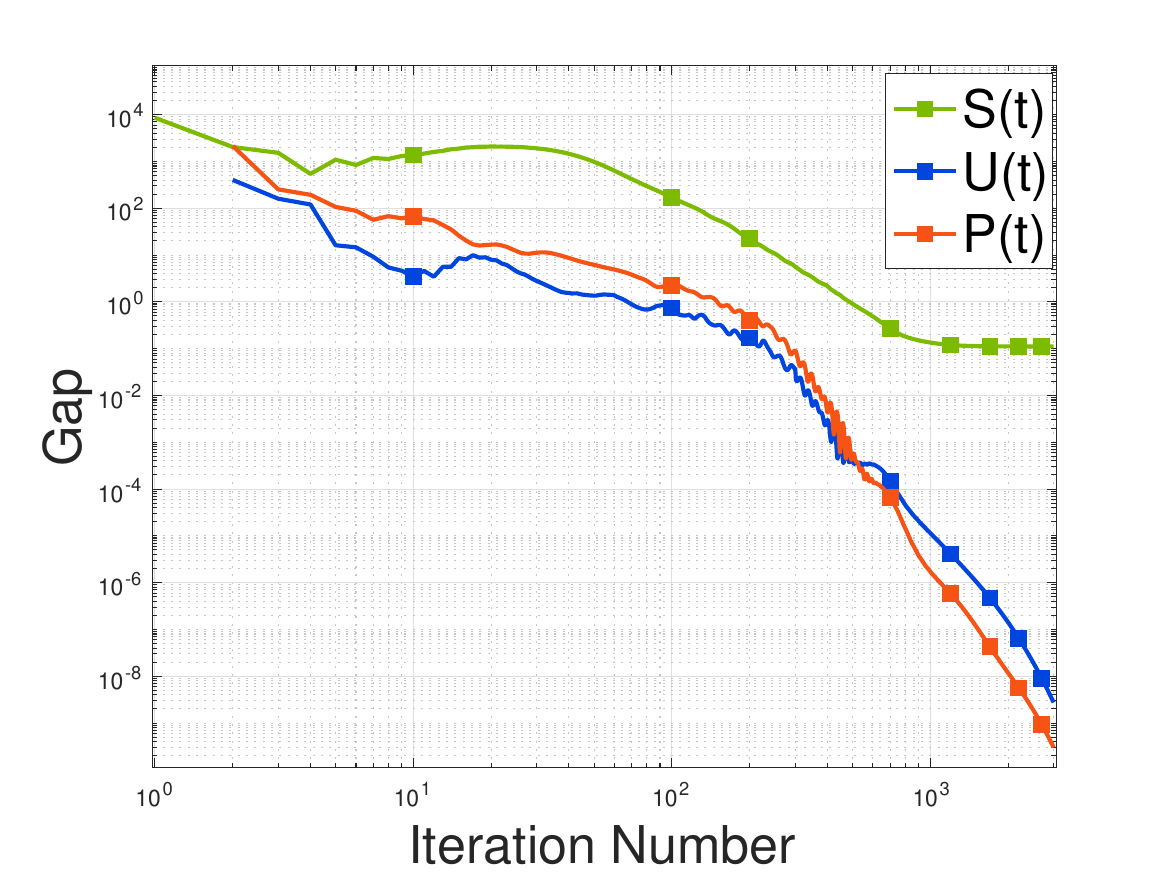}
        \vspace{0.2em}
        \scriptsize (d)Range-dependent noise
    \end{minipage}

    \caption{Performance with synthetic network. Measurement noise: AWGN (first row) and range-dependent Gaussian noise $\sigma_{i,j}^2 = \sigma_{\text{add}}^2 \| \mathbf{p}_i - \mathbf{p}_j \|^2$ (second row).RMSE value(left column); feasibility gap and stationarity gap(right column).}
    \label{fig:syrmse}
\end{figure}

\begin{figure}[htbp]
    \centering
        \begin{minipage}{0.236\textwidth}         \centering         \includegraphics[width=\linewidth]{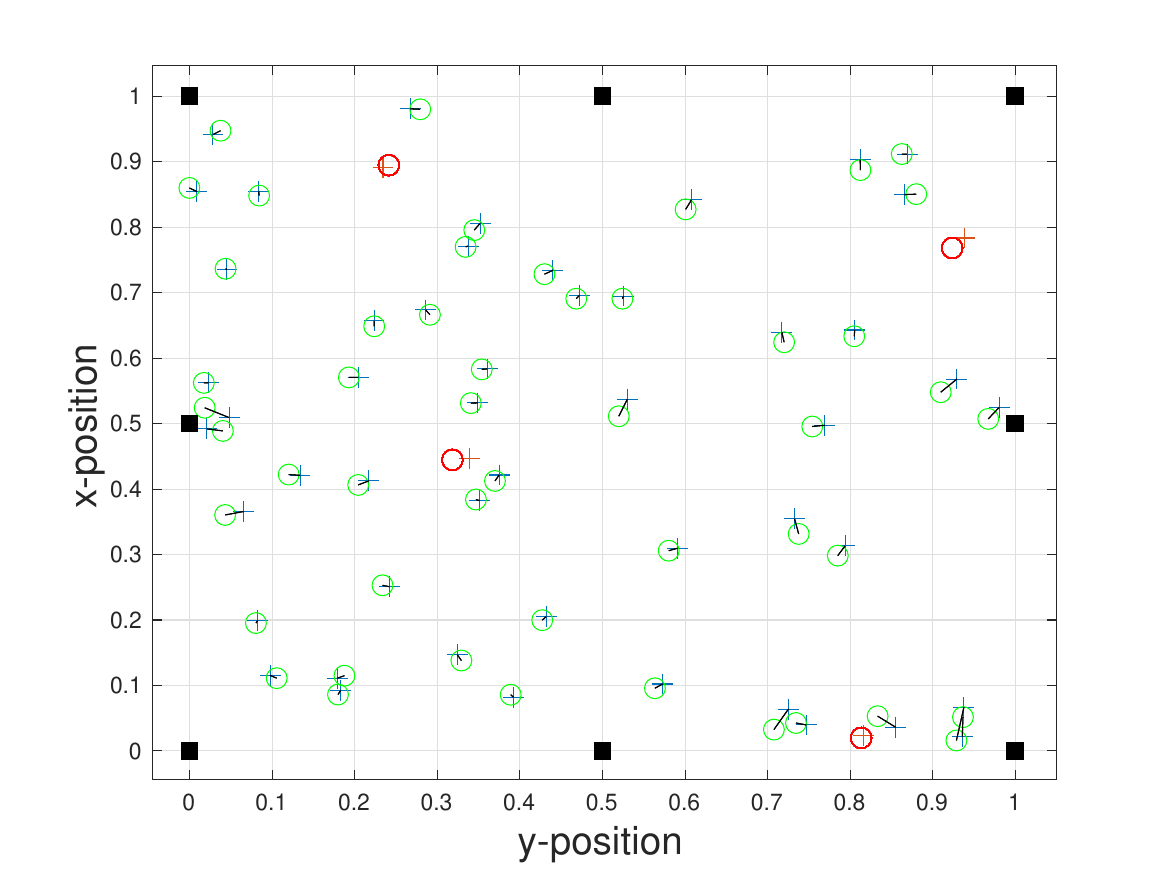}         
        \vspace{0.2em}         
        \scriptsize (a)AWGN     
        \end{minipage}      
        \begin{minipage}{0.236\textwidth}         \centering         \includegraphics[width=\linewidth]{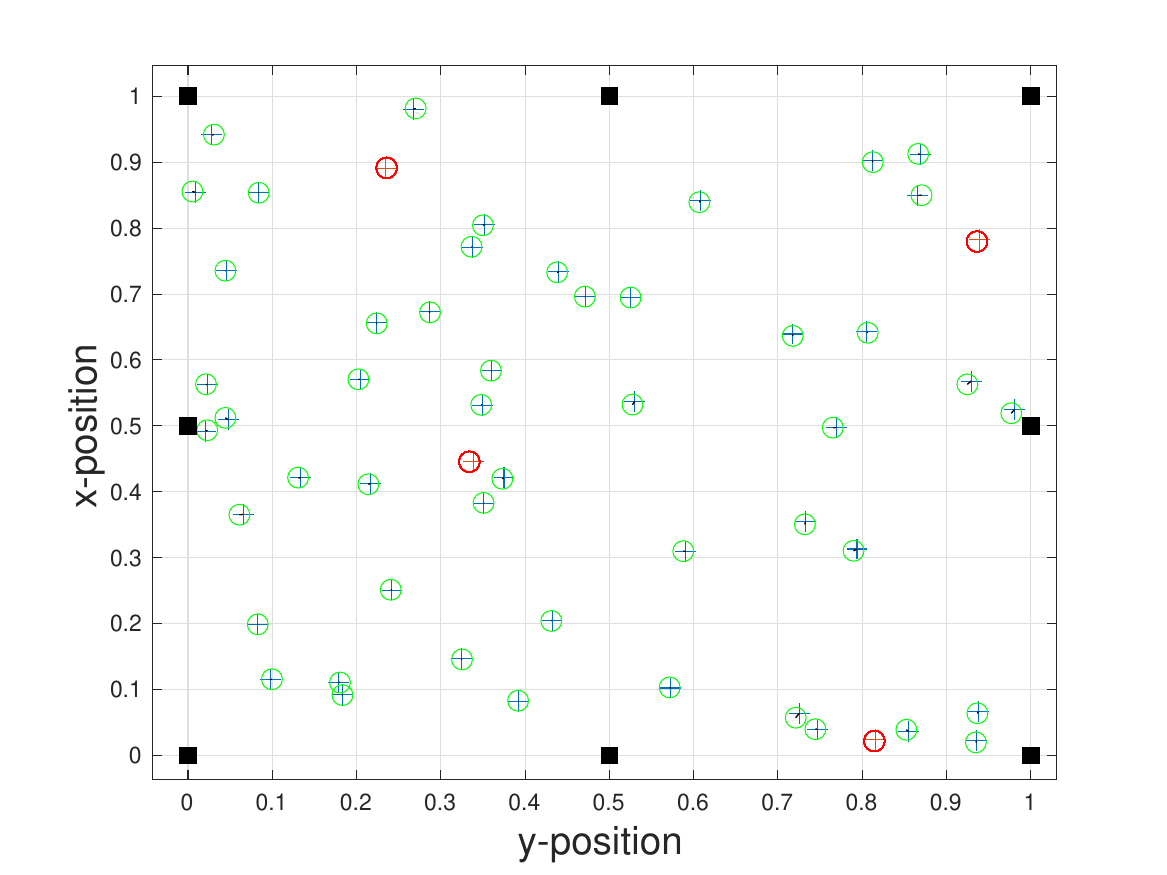}         \vspace{0.2em}     
        \scriptsize (b)Range-dependent noise 
        \end{minipage}

\caption{Localization results with synthetic network. Measurement noise: AWGN and range-dependent Gaussian noise. Anchors: {\color{black}\rule{0.8ex}{0.8ex}}; true agent positions: {\color{blue}\textbf{+}}; true \textcolor{black}{target} positions: {\color{red}\textbf{+}}; estimated agent positions: {\color{green}\textbf{$\circ$}}; estimated \textcolor{black}{target} positions: {\color{orange}\textbf{$\circ$}}.}
\label{fig:4subfigs1}
\end{figure}

We first examine the convergence behavior of the proposed algorithm. Figures~\ref{fig:syrmse}(a) and~\ref{fig:syrmse}(c) respectively show the RMSE of the sensors and the \textcolor{black}{targets} as a function of the number of iterations when additive white Gaussian noise (AWGN) and distance-dependent Gaussian noise are used. To further examine the convergence of the SP-ADMM algorithm, we plot the optimality gap of problem (\ref{jin}) versus the number of iterations in Figures~\ref{fig:syrmse}(b) and~\ref{fig:syrmse}(d) . It can be clearly observed that the performance gap of the proposed SP-ADMM algorithm decreases as the number of iterations increases. This result validates Theorem~1, which states that the iterative sequence $\{(\mathbf{z}^t, \mathbf{u}^t, \boldsymbol{\lambda}^t)\}$ generated by Algorithm~1 converges to a KKT critical point of problem (\ref{jin}), with a convergence rate of $\mathcal{O}(1/T)$. Figure~\ref{fig:4subfigs1} shows the joint localization results of sensors and \textcolor{black}{targets} obtained by the SP-ADMM-JCNL algorithm in the scenario depicted in Figure~\ref{fig:SN}. The first row presents results under additive white Gaussian noise (AWGN), while the second row corresponds to range-dependent Gaussian noise model. In both settings, the proposed method is able to accurately recover the sensor and \textcolor{black}{target} positions. The estimated locations closely match the ground truth, demonstrating the robustness of the algorithm across different noise types. These visual results align well with the quantitative RMSE evaluations discussed earlier.

Moreover, we evaluate the performance of the proposed SP-ADMM-JCNL algorithm against a baseline that separately addresses cooperative and non-cooperative localization tasks. For clarity, we refer to this baseline as SP-ADMM-SCNL. The corresponding two-stage procedure is summarized in Algorithm~\ref{alg:scnl}, where sensor self-localization is performed first using anchor-sensor range measurements, followed by \textcolor{black}{targets} localization based on the estimated sensor positions.

We consider a two-dimensional localization problem with a randomly generated network comprising $100$ sensors, $12$ anchors, and \textcolor{black}{$10$ targets} uniformly distributed in $[0, 1]\times[0, 1]$. Range measurements are corrupted by additive white Gaussian noise with variance $\sigma_{\text{add}}^2=0.01$, and the communication range is set to $C_{range} = 0.25$. For both JCNL and SCNL, identical noisy measurements are used as input. All experiments are repeated over $100$ Monte Carlo trials, and the average RMSE for sensor and target localization is reported. For both JCNL and SCNL, the same noisy measurements are used as input to ensure fair comparison. All experiments are repeated over multiple Monte Carlo trials, and the average RMSE for sensors and \textcolor{black}{targets} localization is reported to assess accuracy and robustness.

\begin{algorithm}[H]
\caption{SP-ADMM-SCNL}\label{alg:scnl}
\begin{algorithmic}
\State \textbf{Input:} Anchor positions \( \{ \mathbf{a}_l \}_{l=1}^{\textcolor{black}{A}} \), sensor -to-sensor  distances, parameters $c$, $\rho$, initial values $\mathbf{z}_i^0$, $\mathbf{w}_i^0$, $\boldsymbol{\lambda}_i^0 = \mathbf{0}$.

\State \textbf{Stage 1: Sensor Self-Localization}
\For{$t = 0$ to $T_1$}
\State Apply SP-ADMM to estimate sensor positions \( \{ \hat{\mathbf{x}}_i \} \) based on sensor-to-sensor distances, ignoring \textcolor{black}{target}-related information.
\EndFor

\State \textbf{Stage 2: \textcolor{black}{targets} localization}
\For{$t = T_1$ to $T_2$}
\State Treat the estimated \( \{ \hat{\mathbf{x}}_i \} \) as virtual anchors (i.e., set \( \textcolor{black}{A} = N \)) then apply SP-ADMM-JCNL to estimate the \textcolor{black}{targets} position \(\{ \hat{\mathbf{y}}_k\} \).
\EndFor
\end{algorithmic}
\end{algorithm}
 Initial estimates are drawn from a uniform distribution: \( \mathbf{z}_0 \sim \text{Unif}(-1, 1)^{(M+1)(4|\mathcal{E} |+ 2N)} \), where $\mathcal{E}$ denotes the set of communication links and $N$ is the number of sensors. The penalty parameters in the augmented Lagrangian are set to \( c = \rho = 0.18 \). For the proposed joint localization method (SP-ADMM-JCNL), we set the iteration cap to \( T = 5000 \). For the two-stage separate localization strategy (SP-ADMM-SCNL), we allocate \( T_1 = 2000 \) iterations for the sensor self-localization stage and \( T_2 = 3000 \) iterations for the subsequent \textcolor{black}{targets} localization stage, keeping the total iteration budget consistent. To ensure statistical reliability and account for randomness in initialization and noise, we conduct 100 independent Monte Carlo trials. The final results are presented in terms of average RMSE values over these runs.

\begin{figure}[htbp]
\centering
\includegraphics[height=5.8cm]{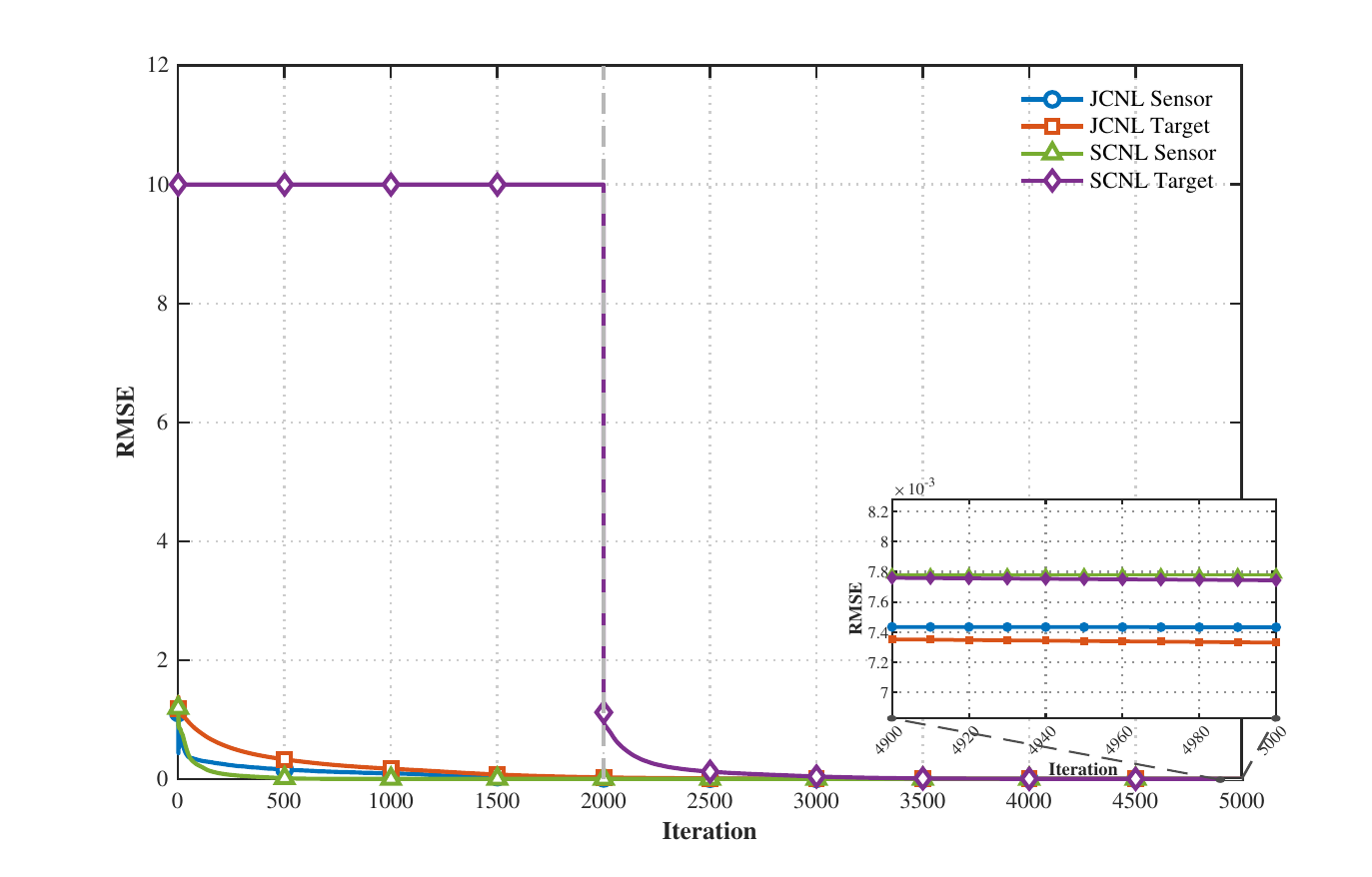}
\caption{Comparison of RMSE between JCNL and SCNL. }
\label{fig:duibi}
\end{figure}
Figure~\ref{fig:duibi} illustrates the localization results for both sensor self-localization and \textcolor{black}{targets} localization under the two different processing strategies: the separate approach (SP-ADMM-SCNL) and the proposed joint approach (SP-ADMM-JCNL). For the separate method, the average computation time was approximately 0.0258 seconds for the first stage (cooperative localization of sensors) and 0.0863 seconds for the second stage (non-cooperative localization of the \textcolor{black}{targets}), resulting in a total runtime of around 0.1643 seconds. In contrast, the joint method required an average computation time of about 0.1768 seconds to complete both tasks simultaneously, which is slightly lower than the cumulative cost of the two-stage approach. The
joint method exhibits faster convergence in the initial iterations, particularly in estimating the \textcolor{black}{target} positions.
This improvement is crucial in real-time applications. Moreover, unlike the separate approach, the joint
strategy avoids \textcolor{black}{target} estimation delays by updating
sensor and \textcolor{black}{target} locations simultaneously.

\subsection{\textcolor{black}{Parameter Selection and Sensitivity Analysis}}

\textcolor{black}{This section evaluates the localization performance of the proposed SP-ADMM-JCNL algorithm under different choices of penalty parameters $c$, $\rho$ and initial values $\mathbf{u}^0$. The goal is to investigate the sensitivity of JCNL to these algorithmic parameters and to provide empirical guidance for selecting suitable values in practice. These results can help users reduce the effort required for parameter tuning when applying JCNL to new localization tasks.}

\subsubsection{\textcolor{black}{Parameter Selection for \(\rho\) and \(c\)}}
\textcolor{black}{According to Theorem \ref{the1}, if the parameter \(\rho\) is selected to satisfy conditions (\ref{45})–(\ref{47}) for a given \(c\), the SP-ADMM algorithm is guaranteed to converge. We test the algorithm on the network configuration depicted in Fig.~\ref{fig:SN}, which consists of $N=50$ sensors, $m=8$ anchors, and $M=4$ targets where nodes are uniformly distributed over a \( [0,1] \times [0,1] \) square region. }
\begin{figure}[htbp]
        \centering
    \includegraphics[width=\linewidth]{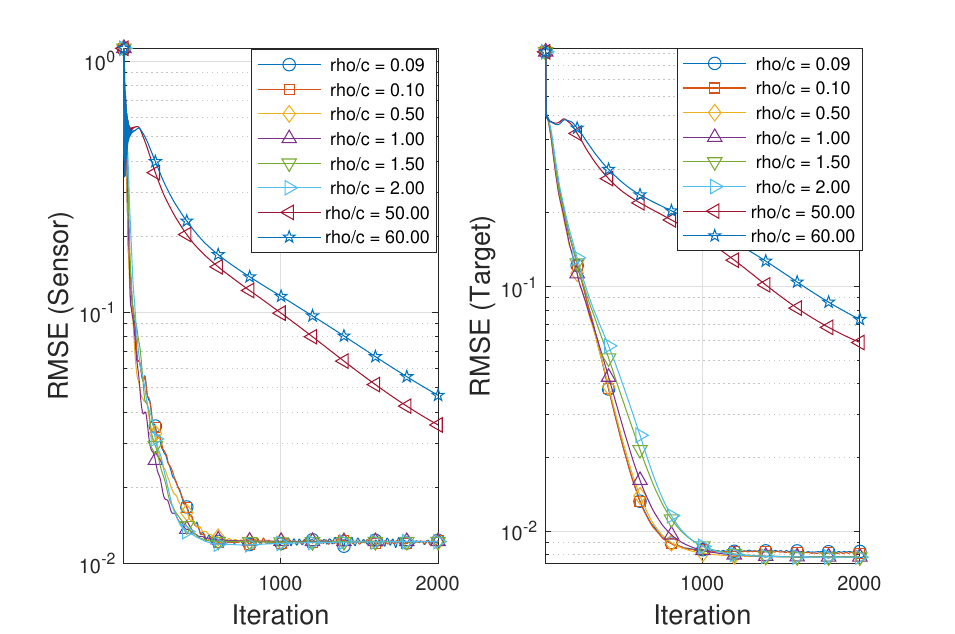}
    \caption{\textcolor{black}{RMSE curves achieved by our proposed methods with different penalty parameters $\rho$ under the benchmark network (N = 50, m = 8, M=4) with Range-dependent noise. $RMSE_{sensor}$(left column); $RMSE_{target}$(right column).}}
    \label{fig:diff_rho}
\end{figure}

\begin{figure}[htbp]
    \centering        \includegraphics[width=\linewidth]{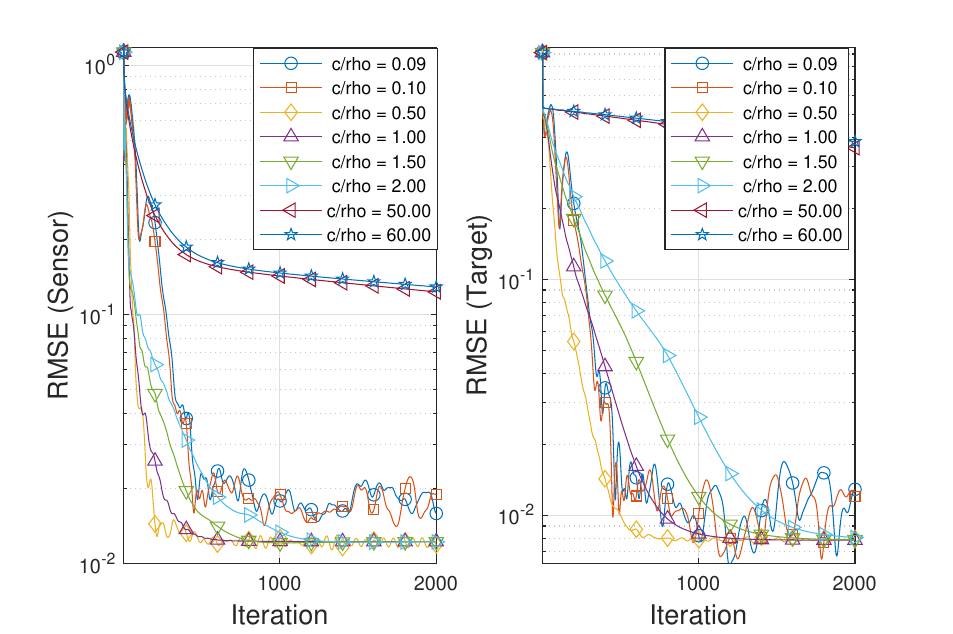}
    \caption{\textcolor{black}{RMSE curves achieved by our proposed methods with different penalty parameters $c$ under the benchmark network (N = 50, m = 8, M=4) with Range-dependent noise. $RMSE_{sensor}$(left column); $RMSE_{target}$(right column).}}
    \label{fig:diff_c}
\end{figure}

\textcolor{black}{Figure~\ref{fig:diff_rho} respectively shows the RMSE performance of sensor self-localization and target localization achieved by JCNL under different values of $\rho$ with a fixed $c = 0.18$. The results indicate that the choice of $\rho$ and $c$ affects both the convergence speed and final localization accuracy of the proposed algorithm. Specifically, an excessively large $\rho$ leads to slower convergence and reduced localization accuracy, as evidenced by the red and light blue curves in Figure~\ref{fig:diff_rho}. It can be observed that when $\rho$ lies within the range $[0.09c, 2.0c]$, the performance of JCNL remains relatively stable. And figure~\ref{fig:diff_c} respectively presents the RMSE performance for sensor self-localization and target localization obtained under different values of $c$ with a fixed $\rho = 0.18$. It is observed that both excessively large and excessively small values of $c$ result in slower convergence and degraded localization accuracy, as shown by the red, light blue, dark blue and orange curves in Figure~\ref{fig:diff_c}. In particular, when $c$ is selected within the range $[0.5\rho, 2.0\rho]$, the algorithm demonstrates stable and favorable performance.}

\subsubsection{\textcolor{black}{Localization Performance With $u_0$ Selection:}}

\begin{figure}[htbp]
    \centering
    \includegraphics[width=\linewidth]{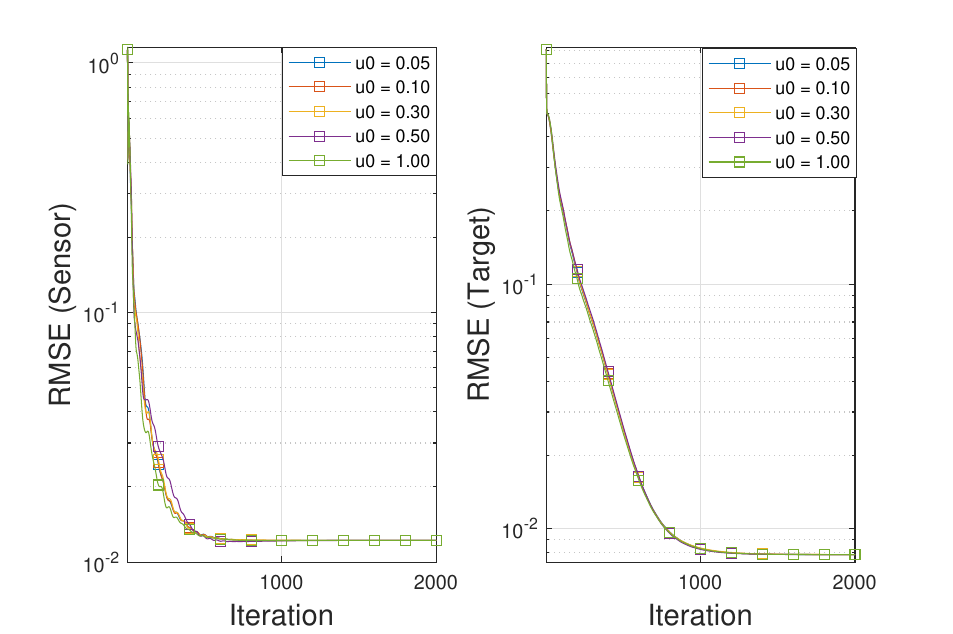}
    \caption{\textcolor{black}{RMSE curves achieved by our proposed methods with different penalty parameters $\rho$ under the benchmark network (N = 50, m = 8, M=4) with Range-dependent noise. $RMSE_{sensor}$(left column); $RMSE_{target}$(right column).}}
    \label{fig:diff_u0}
\end{figure}

\textcolor{black}{Figure~\ref{fig:diff_u0} illustrates the influence of the initialization parameter $u_0$ on sensor self-localization and target localization, respectively. It is observed that the performance exhibits only minor fluctuations with respect to the choice of $u_0$. These results suggest that the proposed SP-ADMM-JCNL algorithm is relatively insensitive to the initialization of the parameter $u_0$.}

\section{Conclusion}\label{con}
In this work, we proposed SP-ADMM-JCNL, a distributed algorithm that jointly addresses cooperative and non-cooperative localization tasks in wireless sensor networks. By leveraging a joint modeling framework and carefully designed variable decoupling strategies, the original non-convex problem is reformulated into a structure amenable to distributed optimization. We provided theoretical guarantees of global convergence to a KKT point of the reformulated problem and to a critical point of the original objective, with a sublinear $\mathcal{O}(1/T)$ rate. Experimental results demonstrate that SP-ADMM-JCNL achieves accurate localization for both sensors and targets, eliminates the estimation delay inherent in two-stage methods and improves computational efficiency. As a natural next step, future work will explore extending the framework to handle multiple targets and accommodate dynamic or multi-modal sensing scenarios commonly encountered in real-world deployments.



\begin{ack}                               
Partially supported by the Roman Senate.  
\end{ack}

\bibliographystyle{unsrt}       
\bibliography{autosam-2/main}           



\clearpage
\onecolumn

\renewcommand\thesection{A.\arabic{section}}  
\setcounter{section}{0}                      

\section{Solution to Remark \ref{remark1}}\label{sa1}

\begin{proof}

By substituting the explicit form of equation (\ref{27}) into problem (\ref{22}), the optimization problem (\ref{22}) can be reformulated as:
\begin{align*}\label{29}
\min_{\mathbf{z} \in \mathcal{X}} \sum_{i \in \mathcal{N}} \Big(& [(c + 1) N_i+M_i] \|\mathbf{x}_i - \tilde{\mathbf{x}}_i^{t+1}\|^2 + \sum_{j \in \mathcal{N}_i} c\|\mathbf{p}_{i,j}^- - (\tilde{\mathbf{p}}_{i,j}^-)^{t+1} \|^2 \nonumber\\
&+ \sum_{j \in \mathcal{N}_i}  \|\mathbf{p}_{i,j}^+ - (\tilde{\mathbf{p}}_{i,j}^+)^{t+1} \|^2  +\sum_{k\in\mathcal{M}}((2cN_i+\mathbb{I}_{\mathcal{M}_i}(k)) \|\mathbf{y}_{i, \textcolor{black}{k}} - \tilde{\mathbf{y}}_{i, \textcolor{black}{k}}^{t+1}\|^2 \nonumber\\
&+ \sum_{j \in \mathcal{N}_i} c\|\mathbf{q}_{i, j, \textcolor{black}{k}}^- - (\tilde{\mathbf{q}}_{i, j, \textcolor{black}{k}}^-)^{t+1} \|^2 + \sum_{j \in \mathcal{N}_i} c \|\mathbf{q}_{i, j, \textcolor{black}{k}}^+ - (\tilde{\mathbf{q}}_{i, j, \textcolor{black}{k}}^+)^{t+1} \|^2 )\Big).
\end{align*}

Combined with the definition of the constraint set $\mathcal{X}$ in equation (\ref{ys2}), the optimal solution for $\mathbf{x}_i^{t+1}$ is:
\[
\mathbf{x}_i^{t+1} = \begin{cases} \tilde{\mathbf{x}}_i^{t+1}, & \text{for } i \notin \mathcal{A}, \\ \mathbf{a}_i, & \text{for } i \in \mathcal{A}. \end{cases}
\]

The optimal solution for $\mathbf{y}_{i, \textcolor{black}{k}}^{t+1}$ is:
\[
\mathbf{y}_{i, \textcolor{black}{k}}^{t+1} = \tilde{\mathbf{y}}_{i, \textcolor{black}{k}}^{t+1}, \text{for } i \in \mathcal{N}
\]

The terms in the optimization problem related to $\mathbf{p}_{i,j}^-$ are:
\[
c\|\mathbf{p}_{i,j}^- - (\tilde{\mathbf{p}}_{i,j}^-)^{t+1} \|^2+ \|\mathbf{p}_{i,j}^- - (\tilde{\mathbf{p}}_{j,i}^+)^{t+1}\|^2 , \quad \forall i \in \mathcal{N}, j \in \mathcal{N}_i.
\]
The terms related to $\mathbf{p}_{i,j}^+$ are:
\[
c\|\mathbf{p}_{i,j}^+ - (\tilde{\mathbf{p}}_{i,j}^-)^{t+1} \|^2+ \|\mathbf{p}_{i,j}^+ -(\tilde{\mathbf{p}}_{j,i}^+)^{t+1}\|^2 , \quad \forall i \in \mathcal{N}, j \in \mathcal{N}_i.
\]
The terms related to $\mathbf{q}_{i, j, \textcolor{black}{k}}^-$ are:
\[
c\|\mathbf{q}_{i, j, \textcolor{black}{k}}^- - (\tilde{\mathbf{q}}_{i, j, \textcolor{black}{k}}^-)^{t+1} \|^2+ c\|\mathbf{q}_{i, j, \textcolor{black}{k}}^- - (\tilde{\mathbf{q}}_{j, i, \textcolor{black}{k}}^+)^{t+1}\|^2 , \quad \forall i \in \mathcal{N}, j \in \mathcal{N}_i, k\in\mathcal{M}.
\]
The terms related to $\mathbf{q}_{i, j, \textcolor{black}{k}}^+$ are:
\[
c\|\mathbf{q}_{i, j, \textcolor{black}{k}}^+ - (\tilde{\mathbf{q}}_{i, j, \textcolor{black}{k}}^-)^{t+1} \|^2+ c\|\mathbf{q}_{i, j, \textcolor{black}{k}}^+ - (\tilde{\mathbf{q}}_{j, i, \textcolor{black}{k}}^+)^{t+1}\|^2 , \quad \forall i \in \mathcal{N}, j \in \mathcal{N}_i, k\in\mathcal{M}.
\]

According to the first-order optimality condition, we obtain:
\[
c \left( (\mathbf{p}_{i,j}^-)^{t+1} - (\mathbf{p}_{i,j}^+)^{t+1} \right) + (\mathbf{p}_{i,j}^-)^{t+1} - (\mathbf{p}_{i,j}^+)^{t+1} = 0.
\]

Rearranging the above equation yields:
\[
(\mathbf{p}_{i,j}^-)^{t+1} = \frac{1}{c + 1} \left( c \cdot (\tilde{\mathbf{p}}_{i,j}^-)^{t+1} + (\tilde{\mathbf{p}}_{j,i}^+)^{t+1} \right).
\]
Similarly, for $i \in \mathcal{N}, j \in \mathcal{N}$, the closed-form solution for $\mathbf{p}_{i,j}^+$ is:
\[
(\mathbf{p}_{i,j}^+)^{t+1} = \frac{1}{c + 1} \left( (\tilde{\mathbf{p}}_{i,j}^+)^{t+1} +  c \cdot(\tilde{\mathbf{p}}_{j,i}^-)^{t+1} \right), \quad \forall i \in \mathcal{N}, j \in \mathcal{N}_i.
\]
For $i \in \mathcal{N}, j \in \mathcal{N}, k\in\mathcal{M} $, the closed-form solution for $\mathbf{q}_{i, j, \textcolor{black}{k}}^+$ is:
\[
(\mathbf{q}_{i, j, \textcolor{black}{k}}^-)^{t+1} = \frac{1}{2} \left((\tilde{\mathbf{q}}_{i, j, \textcolor{black}{k}}^-)^{t+1} + (\tilde{\mathbf{q}}_{j, i, \textcolor{black}{k}}^+)^{t+1} \right), \quad \forall i \in \mathcal{N}, j \in \mathcal{N}_i.
\]
And the closed-form solution for $\mathbf{q}_{i, j, \textcolor{black}{k}}^-$ is:
\[
(\mathbf{q}_{i, j, \textcolor{black}{k}}^+)^{t+1} = \frac{1}{2} \left((\tilde{\mathbf{q}}_{i, j, \textcolor{black}{k}}^+)^{t+1} + (\tilde{\mathbf{q}}_{j, i, \textcolor{black}{k}}^-)^{t+1} \right), \quad \forall i \in \mathcal{N}, j \in \mathcal{N}_i.
\]

\end{proof}

\section{(Proof of Lemma\ref{lem12})}\label{sa2}
\begin{proof}
The following lemmas can be derived from the results presented in \cite{zhang2023distributed}.
\begin{Lemma}\label{lem1}
	  Assume that the matrix \( c \mathbf{B}_i^T \mathbf{B}_i \) follows the form given in Equation (\ref{23}). Let \( \{ (\mathbf{z}_i^t, \mathbf{w}_i^t, \boldsymbol{\lambda}_i^t) \}_{t \geq 1} \) be the sequence generated by Algorithm~\ref{al1}. Then, for all \( t \geq 1 \), the following inequality holds:
\begin{align*}\label{42}
	&\mathcal{L} \left( \mathbf{z}^{t+1}, \boldsymbol{\lambda}^{t+1}, \mathbf{w}^{t+1} \right) - \mathcal{L} \left( \mathbf{z}^t, \boldsymbol{\lambda}^t, \mathbf{w}^t \right) \notag \\
	&\leq \sum_{i \in \mathcal{N}} \Bigg[ 
	-\frac{c}{2} \| \Delta \mathbf{z}_i^{t+1} \| _{\mathbf{B}_i^T \mathbf{B}_i}^2 
	- \frac{\rho}{2} \| \Delta \mathbf{w}_i^{t+1} \| ^2 \notag \\
	&\quad + \frac{3 \left( N_{\max} + M_{\max}+1 \right)}{c} \| \mathbf{H}_i  \Delta \tilde{\mathbf{z}}_i^{t+1} - \mathbf{D}_i \Delta \mathbf{w}_i^{t} \|^2 \notag \\
	&\quad + 3c^2 \left( 2N_{\max} + 1 \right) \|  \Delta \tilde{\mathbf{z}}_i^{t+1} - \Delta \mathbf{z}_i^{t+1} \|_{\mathbf{A}_i^T \mathbf{A}_i}^2 \notag \\
	&\quad + 3 c\left( 1 + c \right) \left( 1+M_{\max} + 2N_{\max} \right) \| \Delta \tilde{\mathbf{z}}_i^{t+1}  -\Delta \mathbf{z}_i^{t} \|_{\mathbf{B}_i^T \mathbf{B}_i}^2 
	\Bigg],\tag{B.1}
\end{align*}
where $N_{\max} := \max \{ N_i, i \in \mathcal{N} \}, M_{\max} := \max \{ M_i, i \in \mathcal{M} \}$.
\end{Lemma}
\begin{proof}
According to the definitions of $\mathbf{A}_i$ and $\mathbf{H}_i$ in equations (\ref{ai}) and (\ref{qi}), we have
\begin{align*}\label{24}
    \mathbf{A}_i^T \mathbf{A}_i &=  \operatorname{blockdiag} \left( {\mathbf{A}_i^{\text{sen}}}^T\mathbf{A}_i^{\text{sen}}, {\mathbf{A}_{i,1}^{\text{tar}}}^T\mathbf{A}_{i,1}^{\text{tar}}, \dots, {\mathbf{A}_{i,M}^{\text{tar}}}^T\mathbf{A}_{i,M}^{\text{tar}} \right) \otimes \mathbf{I}_n\nonumber\\
    &\in \mathbb{R}^{n \cdot (2+4N_i) \times n \cdot (2+4N_i)},\tag{B.2}
\end{align*}
where, \begin{align*}
    {\mathbf{A}_i^{\text{sen}}}^T\mathbf{A}_i^{\text{sen}} = \begin{bmatrix}
        N_i & -\mathbf{1}_{N_i}^T & \mathbf{0}_{N_i}^T \\
        -\mathbf{1}_{N_i} & \mathbf{I}_{N_i} & \mathbf{O}_{N_i}  \\
        \mathbf{0}_{N_i} & \mathbf{O}_{N_i} & \mathbf{O}_{N_i}  
    \end{bmatrix} \otimes \mathbf{I}_n \in \mathbb{R}^{n \cdot (2+4N_i) \times n \cdot (2+4N_i)},
    \end{align*}
and $\forall k\in\mathcal{M}$,
    \begin{align*}{\mathbf{A}_{i, \textcolor{black}{k}}^{\text{tar}}}^T\mathbf{A}_{i, \textcolor{black}{k}}^{\text{tar}} = \begin{bmatrix}
         2N_i & -\mathbf{1}_{N_i}^T & -\mathbf{1}_{N_i}^T \\
         -\mathbf{1}_{N_i} & \mathbf{I}_{N_i} & \mathbf{O}_{N_i} \\
         -\mathbf{1}_{N_i} & \mathbf{O}_{N_i} & \mathbf{I}_{N_i}
    \end{bmatrix} \otimes \mathbf{I}_n \in \mathbb{R}^{n \cdot (2+4N_i) \times n \cdot (2+4N_i)};
\end{align*}
\begin{align*}\label{25}
    \mathbf{H}_i^T \mathbf{H}_i &= \begin{bmatrix}
        N_i+M_i & \mathbf{0}_{N_i}^T & -\mathbf{1}_{N_i}^T & -\textcolor{black}{\mathbb{I}_{\mathcal{M}_i}(1)} &  \mathbf{0}_{N_i}^T &  \mathbf{0}_{N_i}^T &...&\textcolor{black}{\mathbb{I}_{\mathcal{M}_i}(M)} &  \mathbf{0}_{N_i}^T &  \mathbf{0}_{N_i}^T \\
        \mathbf{0}_{N_i} & \mathbf{O}_{N_i} & \mathbf{O}_{N_i} & \mathbf{0}_{N_i} & \mathbf{O}_{N_i} & \mathbf{O}_{N_i}&...&\mathbf{0}_{N_i} & \mathbf{O}_{N_i} & \mathbf{O}_{N_i} \\
        -\mathbf{1}_{N_i} & \mathbf{O}_{N_i} & \mathbf{I}_{N_i} & \mathbf{0}_{N_i} & \mathbf{O}_{N_i} & \mathbf{O}_{N_i} &...&\mathbf{0}_{N_i} & \mathbf{O}_{N_i} & \mathbf{O}_{N_i}\\
        -\textcolor{black}{\mathbb{I}_{\mathcal{M}_i}(1)} & \mathbf{0}_{N_i}^T & \mathbf{0}_{N_i}^T & \textcolor{black}{\mathbb{I}_{\mathcal{M}_i}(1)} & \mathbf{0}_{N_i}^T &  \mathbf{0}_{N_i}^T&...&0 & \mathbf{O}_{N_i} & \mathbf{O}_{N_i} \\ 
        \mathbf{0}_{N_i} & \mathbf{O}_{N_i} & \mathbf{O}_{N_i} & \mathbf{0}_{N_i} & \mathbf{O}_{N_i} & \mathbf{O}_{N_i} &...&\mathbf{0}_{N_i} & \mathbf{O}_{N_i} & \mathbf{O}_{N_i} \\
        \mathbf{0}_{N_i} & \mathbf{O}_{N_i} & \mathbf{O}_{N_i} & \mathbf{0}_{N_i} & \mathbf{O}_{N_i} & \mathbf{O}_{N_i}&...&\mathbf{0}_{N_i} & \mathbf{O}_{N_i} & \mathbf{O}_{N_i} \\
        \vdots & \vdots & \vdots & \vdots & \vdots & \vdots &\ddots&\vdots & \vdots & \vdots \\
        -\textcolor{black}{\mathbb{I}_{\mathcal{M}_i}(M)} & \mathbf{0}_{N_i}^T & \mathbf{0}_{N_i}^T & 0 & \mathbf{0}_{N_i}^T &  \mathbf{0}_{N_i}^T &...&\textcolor{black}{\mathbb{I}_{\mathcal{M}_i}(M)}& \mathbf{O}_{N_i} & \mathbf{O}_{N_i} \\
        \mathbf{0}_{N_i} & \mathbf{O}_{N_i} & \mathbf{O}_{N_i} & \mathbf{0}_{N_i} & \mathbf{O}_{N_i} & \mathbf{O}_{N_i} &...&\mathbf{0}_{N_i} & \mathbf{O}_{N_i} & \mathbf{O}_{N_i} \\
        \mathbf{0}_{N_i} & \mathbf{O}_{N_i} & \mathbf{O}_{N_i} & \mathbf{0}_{N_i} & \mathbf{O}_{N_i} & \mathbf{O}_{N_i}&...&\mathbf{0}_{N_i} & \mathbf{O}_{N_i} & \mathbf{O}_{N_i} 
    \end{bmatrix} 
    \otimes \mathbf{I}_n \\&\in \mathbb{R}^{n  \cdot (2+4N_i) \times n \cdot (2+4N_i)}.\tag{B.3}
\end{align*}

Subtracting $\mathbf{z}_i = \mathbf{z}_i^t$ from $\mathbf{z}_i = \mathbf{z}_i^{t+1}$ in equation (\ref{20}), we obtain
\begin{align}\label{a21}
    &\mathcal{L}_i \left( \mathbf{z}_i^{t+1}, \mathbf{w}_i^t, \boldsymbol{\lambda}_i^t \right) + \frac{c}{2} \| \mathbf{z}_i^{t+1}-\mathbf{z}_i^{t} \|_{\mathbf{B}_i^T \mathbf{B}_i}^2
- \mathcal{L}_i \left( \mathbf{z}_i^t, \mathbf{w}_i^t, \boldsymbol{\lambda}_i^t \right)\nonumber\\
=&
\frac{1}{2} \|  \mathbf{z}_i^{t+1} -\tilde{\mathbf{z}}_i^{t+1} \|_{\mathbf{U}_i}^2 - \frac{1}{2} \| \mathbf{z}_i^t - \tilde{\mathbf{z}}_i^{t+1}\|_{\mathbf{U}_i}^2. \tag{B.4}
\end{align}

According to Remark~\ref{remark1}, $\mathbf{z}_i^{t+1}$ is the optimal solution to problem (\ref{22}), and both $\mathbf{z}_i^t \in \mathcal{Y}_i \cap \mathcal{X}$ and $\mathbf{z}_i^{t+1} \in \mathcal{Y}_i \cap \mathcal{X}$, thus
\begin{align*}\label{a22}
    \frac{1}{2} \| \mathbf{z}_i^{t+1} -\tilde{\mathbf{z}}_i^{t+1}\|_{\mathbf{U}_i}^2 \leq \frac{1}{2} \| \mathbf{z}_i^t - \tilde{\mathbf{z}}_i^{t+1} \|_{\mathbf{U}_i}^2. \tag{B.5}
\end{align*}
Substituting equation (\ref{a22}) into equation (\ref{a21}), we derive
\begin{align*}\label{a23}
    \mathcal{L}_i \left( \mathbf{z}_i^{t+1}, \mathbf{w}_i^t, \boldsymbol{\lambda}_i^t \right) - \mathcal{L}_i \left( \mathbf{z}_i^t, \mathbf{w}_i^t, \boldsymbol{\lambda}_i^t \right) 
\leq -\frac{c}{2} \| \mathbf{z}_i^{t+1}-\mathbf{z}_i^{t} \|_{\mathbf{B}_i^T \mathbf{B}_i}^2. \tag{B.6}
\end{align*}
Applying the same technique to the update step of variable $\mathbf{w}$ in equation (\ref{gengxin1}b), we obtain
\begin{align*}\label{a24}
    \mathcal{L}_i \left( \mathbf{z}_i^{t+1}, \mathbf{w}_i^{t+1}, \boldsymbol{\lambda}_i^t \right) - \mathcal{L}_i \left( \mathbf{z}_i^{t+1}, \mathbf{w}_i^t, \boldsymbol{\lambda}_i^t \right) 
\leq -\frac{\rho}{2} \| \mathbf{w}_i^{t+1} - \mathbf{w}_i^t \|^2. \tag{B.7}
\end{align*}
Similarly, from equation (\ref{gengxin1}c), i.e., $\boldsymbol{\lambda}_i^{t+1} = \boldsymbol{\lambda}_i^t + c \mathbf{A}_i \mathbf{z}_i^{t+1}$, it follows that
\begin{align*}\label{a25}
    \mathcal{L}_i \left( \mathbf{z}_i^{t+1}, \mathbf{w}_i^{t+1}, \boldsymbol{\lambda}_i^{t+1} \right) - \mathcal{L}_i \left( \mathbf{z}_i^{t+1}, \mathbf{w}_i^{t+1}, \boldsymbol{\lambda}_i^t \right) 
= \frac{1}{c} \| \boldsymbol{\lambda}_i^{t+1} - \boldsymbol{\lambda}_i^t \|^2. \tag{B.8}
\end{align*}
Left-multiplying equation (\ref{zwan}) by matrix $\mathbf{U}_i$ and rearranging terms, we get
\begin{align*}\label{a26}
\mathbf{A}_i^T \boldsymbol{\lambda}_i^{t} = \mathbf{H}_i^T \mathbf{D}_i \mathbf{w}_i^{t} - \mathbf{U}_i \tilde{\mathbf{z}}_i^{t+1} + c \mathbf{B}_i^T \mathbf{B}_i \mathbf{z}_i^t. \tag{B.9}
\end{align*}
Adding $c \mathbf{A}_i^T \mathbf{A}_i \mathbf{z}_i^{t+1}$ to both sides of the above equation and combining with equation (\ref{a26}), we obtain
\begin{align*}\label{a27}
    \mathbf{A}_i^T \boldsymbol{\lambda}_i^{t+1} 
&= \mathbf{H}_i^T \mathbf{D}_i \mathbf{w}_i^{t} - \mathbf{U}_i\tilde{\mathbf{z}}_i^{t+1} + c \mathbf{B}_i^T \mathbf{B}_i \mathbf{z}_i^t + c \mathbf{A}_i^T \mathbf{A}_i \mathbf{z}_i^{t+1}  \\
&= \mathbf{H}_i^T \mathbf{D}_i \mathbf{w}_i^{t} - \mathbf{H}_i^T \mathbf{H}_i \tilde{\mathbf{z}}_i^{t+1} - c \mathbf{A}_i^T \mathbf{A}_i (\tilde{\mathbf{z}}_i^{t+1}-\mathbf{z}_i^{t+1}) - c \mathbf{B}_i^T \mathbf{B}_i (\tilde{\mathbf{z}}_i^{t+1}-\mathbf{z}_i^{t}). \tag{B.10}
\end{align*}
The above equation holds due to the definition of $\mathbf{U}_i$ in equation (\ref{19}). Let $\sigma_{\min}$ denote the smallest non-zero eigenvalue of $\mathbf{A}_i^T \mathbf{A}_i$, then
\[
\sigma_{\min} \| \boldsymbol{\lambda}_i^{t+1} - \boldsymbol{\lambda}_i^t \|^2 \leq \| \mathbf{A}_i^T \left( \boldsymbol{\lambda}_i^{t+1} - \boldsymbol{\lambda}_i^t \right) \|^2.
\]

According to the definition of $\mathbf{A}_i^T \mathbf{A}_i$ in equation (\ref{24}), we have $\sigma_{\min} = 1$. Combining this inequality with equation (\ref{a27}), we obtain
\begin{align*}\label{a28}
    \| \boldsymbol{\lambda}_i^{t+1} - \boldsymbol{\lambda}_i^t \|^2 \leq &
\| \mathbf{H}_i^T \mathbf{D}_i \left( \mathbf{w}_i^{t} - \mathbf{w}_i^{t-1} \right) - \mathbf{H}_i^T \mathbf{H}_i \left( \tilde{\mathbf{z}}_i^{t+1} - \tilde{\mathbf{z}}_i^t \right) 
- c \mathbf{A}_i^T \mathbf{A}_i \left[\tilde{\mathbf{z}}_i^{t+1} - \mathbf{z}_i^{t+1} - (\tilde{\mathbf{z}}_i^t - \mathbf{z}_i^t) \right] \\
&- c \mathbf{B}_i^T \mathbf{B}_i \left( \tilde{\mathbf{z}}_i^{t+1} - \mathbf{z}_i^{t+1} - \left( \tilde{\mathbf{z}}_i^t - \mathbf{z}_i^{t-1} \right) \right) \|^2\\
\leq &
3(N_i+M_i+1)\| \mathbf{D}_i \left( \mathbf{w}_i^{t} - \mathbf{w}_i^{t-1} \right) - \mathbf{H}_i^T \left( \tilde{\mathbf{z}}_i^{t+1}  - \tilde{\mathbf{z}}_i^t \right) \|^2\\
&+3 c^2(2N_i+1)\|\tilde{\mathbf{z}}_i^{t+1} - \mathbf{z}_i^{t+1} - (\tilde{\mathbf{z}}_i^t - \mathbf{z}_i^t) \| _{\mathbf{A}_i^T \mathbf{A}_i }\\
&+3c^2\| \mathbf{B}_i^T \mathbf{B}_i \|\|\tilde{\mathbf{z}}_i^{t+1} - \mathbf{z}_i^{t} - (\tilde{\mathbf{z}}_i^t - \mathbf{z}_i^{t-1})\|_{\mathbf{B}_i^T \mathbf{B}_i}
\tag{B.11}
\end{align*}

The last inequality above holds based on equations (\ref{24}),(\ref{25}), the triangle inequality of norms, and the facts
$\| \mathbf{H}_i^T \mathbf{H}_i \| =N_i+M_i+1$, $ \| \mathbf{A}_i^T \mathbf{A}_i \| = 2N_i + 1$. Here, $\| \mathbf{A}_i^T \mathbf{A}_i \|$ denotes the spectral norm of the matrix $\mathbf{A}_i^T \mathbf{A}_i$. Combining with the form of $c \mathbf{B}_i^T \mathbf{B}_i$ , we have
\begin{align*}\label{a29}
\| \mathbf{B}_i^T \mathbf{B}_i \| \leq \frac{\left( 1 + c \right) \left( 1+M_{\max} + 2N_{\max} \right)}{c}, \tag{B.12}
\end{align*}
where $N_{\max} = \max \left\{ N_i, i \in \mathcal{N} \right\}, M_{\max} = \max \left\{ M_i, i \in \mathcal{M} \right\}$. Substituting equation (\ref{a29}) into equation (\ref{a28}) and applying the result to equation (\ref{a25}), we get
\begin{align*}\label{a210}
    &\mathcal{L}_i \left( \mathbf{z}_i^{t+1}, \mathbf{w}_i^{t+1}, \boldsymbol{\lambda}_i^{t+1} \right) - \mathcal{L}_i \left( \mathbf{z}_i^t, \mathbf{w}_i^t, \boldsymbol{\lambda}_i^t \right)\\
\leq &
\frac{3(N_{max}+M_{max}+2)}{c}\| \mathbf{D}_i \left( \mathbf{w}_i^{t} - \mathbf{w}_i^{t-1} \right) - \mathbf{H}_i^T \left( \tilde{\mathbf{z}}_i^{t+1}  - \tilde{\mathbf{z}}_i^t \right) \|^2\\
&+3 c^2(2N_{max}+1)\|\tilde{\mathbf{z}}_i^{t+1} - \mathbf{z}_i^{t+1} - (\tilde{\mathbf{z}}_i^t - \mathbf{z}_i^t) \| _{\mathbf{A}_i^T \mathbf{A}_i }\\
&+3(1+c)\left( 1+M_{\max} + 2N_{\max} \right)\|\tilde{\mathbf{z}}_i^{t+1} - \mathbf{z}_i^{t} - (\tilde{\mathbf{z}}_i^t - \mathbf{z}_i^{t-1})\|_{\mathbf{B}_i^T \mathbf{B}_i} \tag{B.13}
\end{align*}

Combining equations (\ref{a23})–(\ref{a24}) and equation (\ref{a210}), and summing over nodes $i \in \mathcal{N}$, we derive inequality (\ref{42}). 
\end{proof}

\begin{Lemma}\label{lemc1}
Suppose $c\mathbf{B}_i^T \mathbf{B}_i$ is defined as in equation (\ref{23}). Then for Algorithm \ref{al1}, the following holds:
\begin{align*}
&\sum_{i \in \mathcal{N}} \Bigg[ \frac{c}{2} \| \mathbf{A}_i \tilde{\mathbf{z}}_i^{t+1} \|^2 + \frac{c}{2} \| \mathbf{z}_i^{t+1} - \mathbf{z}_i^t \|_{\mathbf{B}_i^T \mathbf{B}_i}^2 \Bigg] \\
\leq & \sum_{i \in \mathcal{N}} \Bigg[ \frac{c}{2} \| \mathbf{A}_i \tilde{\mathbf{z}}_i^t \|^2 + \frac{c}{2} \| \mathbf{z}_i^t - \mathbf{z}_i^{t-1} \|_{\mathbf{B}_i^T \mathbf{B}_i}^2 - \frac{c}{2} \| \mathbf{z}_i^{t+1} - \mathbf{z}_i^t \|_{\mathbf{A}_i^T \mathbf{A}_i}^2 - \frac{1}{2} \| \mathbf{z}_i^{t+1} - \mathbf{z}_i^t \|_{\mathbf{H}_i^T \mathbf{H}_i}^2 \\
& + \frac{d_{\max}^2}{2} \| \mathbf{w}_i^t - \mathbf{w}_i^{t-1} \|^2 - \frac{1}{2} \| \mathbf{H}_i (\tilde{\mathbf{z}}_i^{t+1} - \tilde{\mathbf{z}}_i^t) - \mathbf{D}_i (\mathbf{w}_i^t - \mathbf{w}_i^{t-1}) \|^2 \\
& - \frac{c}{2} \| \tilde{\mathbf{z}}_i^{t+1} - \mathbf{z}_i^{t+1} - (\tilde{\mathbf{z}}_i^t - \mathbf{z}_i^t) \|_{\mathbf{A}_i^T \mathbf{A}_i}^2 + \frac{c}{2} (2N_i + 1) \| \mathbf{z}_i^t - \tilde{\mathbf{z}}_i^{t+1} \|^2 \\
& - \frac{c}{2} \| \tilde{\mathbf{z}}_i^{t+1} - \mathbf{z}_i^t - (\tilde{\mathbf{z}}_i^t - \mathbf{z}_i^{t-1}) \|_{\mathbf{B}_i^T \mathbf{B}_i}^2 \Bigg],
\end{align*}
where $d_{\max} = \max\left(\{d_{i,j} \mid i \in \mathcal{N}, j \in \mathcal{N}_i \} \cup \{r_{i, \textcolor{black}{k}} \mid i \in \mathcal{N},k\in\mathcal{M}_i\}\right)$.
\end{Lemma}
\begin{proof}
From equation (\ref{a27}), we have:
\begin{align*}\label{a31}
\langle \mathbf{A}_i^T (\boldsymbol{\lambda}_i^{t+1} - \boldsymbol{\lambda}_i^t), \tilde{\mathbf{z}}_i^{t+1} - \tilde{\mathbf{z}}_i^t \rangle
= &\langle -\mathbf{U}_i (\mathbf{z}_i^{t+1} - \tilde{\mathbf{z}}_i^t) + \mathbf{H}_i^T \mathbf{D}_i (\mathbf{w}_i^t - \mathbf{w}_i^{t-1}), \mathbf{z}_i^{t+1} - \tilde{\mathbf{z}}_i^t \rangle\\
&+ \langle c \mathbf{A}_i^T \mathbf{A}_i (\mathbf{z}_i^{t+1} - \mathbf{z}_i^t) + c \mathbf{B}_i^T \mathbf{B}_i (\mathbf{z}_i^{t} - \tilde{\mathbf{z}}_i^{t-1}), \tilde{\mathbf{z}}_i^{t+1} - \tilde{\mathbf{z}}_i^t \rangle. \tag{B.14}
\end{align*}

Next, we transform both sides of equation (\ref{a31}). First, the left-hand side of equation (\ref{a31}) can be expressed as:
\[
\langle \mathbf{A}_i^T (\boldsymbol{\lambda}_i^{t+1} - \boldsymbol{\lambda}_i^t), \tilde{\mathbf{z}}_i^{t+1} - \tilde{\mathbf{z}}_i^t \rangle = \langle c \mathbf{A}_i^T \mathbf{A}_i \mathbf{z}_i^{t+1}, \tilde{\mathbf{z}}_i^{t+1} - \tilde{\mathbf{z}}_i^t \rangle,
\]
This equality follows from the dual variable update step (\ref{gengxin1}c). Additionally:
\begin{align*}\label{a32}
\langle c \mathbf{A}_i^T \mathbf{A}_i \mathbf{z}_i^{t+1}, \tilde{\mathbf{z}}_i^{t+1} - \tilde{\mathbf{z}}_i^t \rangle 
&= \langle c \mathbf{A}_i^T \mathbf{A}_i (\mathbf{z}_i^{t+1} - \mathbf{z}_i^t), \tilde{\mathbf{z}}_i^{t+1} - \tilde{\mathbf{z}}_i^t \rangle 
- \langle c \mathbf{A}_i^T \mathbf{A}_i \mathbf{z}_i^t, \tilde{\mathbf{z}}_i^{t+1} - \tilde{\mathbf{z}}_i^t \rangle \\
&= -\frac{c}{2} \| \mathbf{z}_i^{t+1} - \mathbf{z}_i^t - (\tilde{\mathbf{z}}_i^{t+1} - \tilde{\mathbf{z}}_i^t) \|_{\mathbf{A}_i^T \mathbf{A}_i}^2 
+ \frac{c}{2} \| \tilde{\mathbf{z}}_i^{t+1} - \tilde{\mathbf{z}}_i^t \|_{\mathbf{A}_i^T \mathbf{A}_i}^2 
+ \frac{c}{2} \| \mathbf{z}_i^{t+1} - \mathbf{z}_i^t \|_{\mathbf{A}_i^T \mathbf{A}_i}^2 \\
&\quad + \frac{c}{2} \| \mathbf{A}_i (\mathbf{z}_i^t - \tilde{\mathbf{z}}_i^t) \|^2 
- \frac{c}{2} \| \mathbf{A}_i (\mathbf{z}_i^t -\tilde{\mathbf{z}}_i^{t+1}) \|^2 
+ \frac{c}{2} \| \mathbf{A}_i \tilde{\mathbf{z}}_i^{t+1} \|^2 
- \frac{c}{2} \| \mathbf{A}_i \tilde{\mathbf{z}}_i^t \|^2,\tag{B.15}
\end{align*}
The last equality holds due to the identity: for any $\mathbf{a}, \mathbf{b} \in \mathbb{R}^{(2N_i+1)(M+1)n}$,
\[
\langle \mathbf{a}, \mathbf{b}\rangle= -\frac{1}{2} \| \mathbf{a} - \mathbf{b} \|^2 + \frac{1}{2} \| \mathbf{a} \|^2 + \frac{1}{2} \| \mathbf{b} \|^2.
\]
Moreover, based on norm compatibility, we have
\begin{align*}\label{a33}
    \frac{c}{2} \| \mathbf{A}_i (\mathbf{z}_i^{t} - \tilde{\mathbf{z}}_i^{t+1}) \|^2 \leq \frac{c (2N_i + 1)}{2} \| \mathbf{z}_i^{t} - \tilde{\mathbf{z}}_i^{t+1} \|^2, \tag{B.16}
\end{align*}
where $\| \mathbf{A}_i^T \mathbf{A}_i \| = 2N_i + 1$ is derived from equation (\ref{24}).

The right-hand side of equation~(\ref{a31}) can be expressed as:
\begin{align*}\label{a36}
&\langle -\mathbf{U}_i (\tilde{\mathbf{z}}_i^{t+1} - \tilde{\mathbf{z}}_i^t) + \mathbf{H}_i^T \mathbf{D}_i (\mathbf{w}_i^t - \mathbf{w}_i^{t-1}), \tilde{\mathbf{z}}_i^{t+1} - \tilde{\mathbf{z}}_i^t \rangle \\
&+ \langle c \mathbf{A}_i^T \mathbf{A}_i (\mathbf{z}_i^{t+1} - \mathbf{z}_i^t) + c \mathbf{B}_i^T \mathbf{B}_i (\mathbf{z}_i^{t+1} - \tilde{\mathbf{z}}_i^t), \mathbf{z}_i^{t+1} - \tilde{\mathbf{z}}_i^t \rangle\\
=& -\frac{1}{2} \| \mathbf{H}_i (\tilde{\mathbf{z}}_i^{t+1} - \tilde{\mathbf{z}}_i^t)- \mathbf{D}_i (\mathbf{w}_i^t - \mathbf{w}_i^{t-1}) \|^2- \frac{c}{2} \| \tilde{\mathbf{z}}_i^{t+1} - \tilde{\mathbf{z}}_i^t - ({z}_i^{t+1} - \mathbf{z}_i^t) \|_{\mathbf{A}_i^T \mathbf{A}_i}^2\\
&-\frac{1}{2}\|\tilde{\mathbf{z}}_i^{t+1}-\tilde{\mathbf{z}}_i^t\|_{\mathbf{U}_i}^2- \frac{c}{2} \| \tilde{\mathbf{z}}_i^{t+1} -\mathbf{z}_i^t - (\tilde{\mathbf{z}}_i^{t} - \mathbf{z}_i^{t-1}) \|_{\mathbf{B}_i^T \mathbf{B}_i}^2 + \frac{c}{2} \|\mathbf{A}_i( \mathbf{z}_i^{t+1} - \mathbf{z}_i^t) \|^2 \\
&+ \frac{c}{2} \| \mathbf{z}_i^{t} -\mathbf{z}_i^{t+1} \|_{\mathbf{B}_i^T \mathbf{B}_i}^2+\frac{1}{2}\|\mathbf{D}_i(\mathbf{w}_i^t-\mathbf{w}_i^{t+1})\|^2
. \tag{B.17}
\end{align*}
Moreover, we obtain the following inequality:
\begin{align*}\label{a38}
    \sum_{i \in \mathcal{N}} \| \mathbf{z}_i^t - \tilde{\mathbf{z}}_i^{t+1} \|_{\mathbf{U}_i}^2 \leq \sum_{i \in \mathcal{N}} \frac{1}{2} \| \tilde{\mathbf{z}}_i^{t+1} - \mathbf{z}_i^t \|_{\mathbf{U}_i}^2 - \frac{1}{2} \|\tilde{\mathbf{z}}_i^{t+1}-\tilde{\mathbf{z}}_i^t-(\mathbf{z}_i^{t+1} - \mathbf{z}_i^t )\|_{\mathbf{U}_i}^2. \tag{B.18}
\end{align*}
Substituting equations (\ref{a33}) and (\ref{a38}) into equations (\ref{a32}) and (\ref{a36}) respectively, and again using $\mathbf{U}_i = \mathbf{H}_i^T \mathbf{H}_i + c \mathbf{A}_i^T \mathbf{A}_i + c \mathbf{B}_i^T \mathbf{B}_i$ together with equation (\ref{a31}), we can derive the final result.
\end{proof} 

However, note that in Lemma \ref{lemc1}, it is still necessary to restrict the increment term $\| \mathbf{z}_i^t - \tilde{\mathbf{z}}_i^{t+1} \|^2$. In the following, we establish a lemma to ensure its descent.

\begin{Lemma}\label{lemc2}
    Let $\{(\mathbf{z}_i^t, \mathbf{w}_i^t, \boldsymbol{\lambda}_i^t)\}_{t \geq 1}$ be the sequence generated by Algorithm \ref{al1}, and suppose $c\mathbf{B}_i^T \mathbf{B}_i$ is defined as in equation (\ref{23}). Then the following holds:
\begin{align*}
&\sum_{i \in \mathcal{N}} \left[ \frac{c}{2} \| \mathbf{A}_i \mathbf{z}_i^{t+1} \|^2 + \frac{c}{2} \| \mathbf{z}_i^{t+1} - \mathbf{z}_i^t \|_{\mathbf{B}_i^T \mathbf{B}_i}^2 \right]\\
\leq &\sum_{i \in \mathcal{N}} \Bigg[ \frac{c}{2} \| \mathbf{A}_i \mathbf{z}_i^t \|^2 + \frac{c}{2} \| \mathbf{z}_i^t - \mathbf{z}_i^{t-1} \|_{\mathbf{B}_i^T \mathbf{B}_i}^2 + \frac{d_{\max}^2}{2} \| \mathbf{w}_i^t - \mathbf{w}_i^{t-1} \|^2 \\
&- \frac{c\cdot\tilde{\tau}_{\min} }{2(2M+1)N_{\text{sum}}n(c+1)^2} \| \mathbf{z}_i^{t} - \tilde{\mathbf{z}}_i^{t+1} \|^2\Bigg],
\end{align*}
where
\[
\tilde{\tau}_{min} :=\min\quad\{(c+1)^2[(6M+1)N_i^2+N_i], i \in \mathcal{N}\}, \quad N_{\text{sum}} := \sum_{i \in \mathcal{N}} N_i.
\]
\end{Lemma}
\begin{proof}
Using equation (\ref{a27}) again, we obtain
\begin{align*}\label{a39}
&\langle \mathbf{A}_i^T (\boldsymbol{\lambda}_i^{t+1} - \boldsymbol{\lambda}_i^t), \mathbf{z}_i^{t+1} - \tilde{\mathbf{z}}_i^t \rangle\\
= &\langle -\mathbf{U}_i ( \tilde{\mathbf{z}}_i^{t+1} - \tilde{\mathbf{z}}_i^t), \mathbf{z}_i^{t+1} -\mathbf{z}_i^t \rangle+  c\|\mathbf{A}_i (\mathbf{z}_i^{t+1} - \mathbf{z}_i^t)\|^2 \\
&+\langle \mathbf{H}_i^T\mathbf{D}_i(\mathbf{w}_i^t-\mathbf{w}_i^{t-1})+c \mathbf{B}_i^T \mathbf{B}_i (\mathbf{z}_i^{t} -\mathbf{z}_i^{t-1}), \mathbf{z}_i^{t+1} -\mathbf{z}_i^t \rangle. \tag{B.19}
\end{align*}

First, the left-hand side of equation (\ref{a39}) can be written as:
\begin{align*}
&\langle \mathbf{A}_i^T (\boldsymbol{\lambda}_i^{t+1} - \boldsymbol{\lambda}_i^t), \mathbf{z}_i^{t+1} -\tilde{\mathbf{z}}_i^t \rangle\\
=& \langle c \mathbf{A}_i \mathbf{z}_i^{t+1}, \mathbf{A}_i\mathbf{z}_i^{t+1}-\mathbf{A}_i\mathbf{z}_i^t\rangle \\
=&\frac{c}{2}\|\mathbf{A}_i\mathbf{z}_i^{t+1}\|^2-\frac{c}{2}\|\mathbf{A}_i\mathbf{z}_i^t\|^2+\frac{c}{2}\|\mathbf{A}_i(\mathbf{z}_i^{t+1}-\mathbf{z}_i^t)\|^2.\tag{B.20}
\end{align*}
The first equality holds due to equation (\ref{gengxin1}), and the second equality follows from the fact that for any $\mathbf{a}, \mathbf{b} \in \mathbb{R}^{(2N_i+1)(M+1)n}$, we have
\[
\langle \mathbf{a}, \mathbf{a} - \mathbf{b} \rangle = \frac{1}{2} \|\mathbf{a} - \mathbf{b}\|^2 + \frac{1}{2} \|\mathbf{a}\|^2 - \frac{1}{2} \|\mathbf{b}\|^2.
\]

For the right-hand side of equation (\ref{a39}), using equation (\ref{a27}), we obtain
\[
\sum_{i \in \mathcal{N}} \langle -\mathbf{U}_i (\tilde{\mathbf{z}}_i^{t+1} - \tilde{\mathbf{z}}_i^t), \mathbf{z}_i^{t+1} - \mathbf{z}_i^t \rangle 
\leq -\sum_{i \in \mathcal{N}} \|\mathbf{z}_i^{t+1} - \mathbf{z}_i^t\|_{\mathbf{U}_i}^2. \tag{B.21}
\]
Applying the Cauchy–Schwarz inequality yields
\begin{align*} \label{a312}
&\langle \mathbf{H}_i^T \mathbf{D}_i (\mathbf{w}_i^t - \mathbf{w}_i^{t-1}) + c \mathbf{B}_i^T \mathbf{B}_i (\mathbf{z}_i^t - \mathbf{z}_i^{t-1}), \mathbf{z}_i^{t+1} - \mathbf{z}_i^t \rangle\\
\leq& -\frac{1}{2} \|\mathbf{D}_i (\mathbf{w}_i^t - \mathbf{w}_i^{t-1})\|^2 
+ \frac{1}{2} \|\mathbf{H}_i (\mathbf{z}_i^{t+1} - \mathbf{z}_i^t)\|^2 
+ \frac{c}{2} \|\mathbf{z}_i^t - \mathbf{z}_i^{t-1}\|_{\mathbf{B}_i^T \mathbf{B}_i}^2 
+ \frac{c}{2} \|\mathbf{z}_i^{t+1} - \mathbf{z}_i^t\|_{\mathbf{B}_i^T \mathbf{B}_i}^2. \tag{B.22}
\end{align*}

Combining equations (\ref{a39})–(\ref{a312}) and substituting $\mathbf{U}_i = c\mathbf{B}_i^T \mathbf{B}_i + c \mathbf{A}_i^T \mathbf{A}_i + \mathbf{H}_i^T \mathbf{H}_i$, we derive:
\begin{align*}\label{a313}
&\sum_{i \in \mathcal{N}} \left[ \frac{c}{2} \| \mathbf{A}_i \mathbf{z}_i^{t+1} \|^2 + \frac{c}{2} \| \mathbf{z}_i^{t+1} - \mathbf{z}_i^t \|_{\mathbf{B}_i^T \mathbf{B}_i}^2 \right]\\
\leq& \sum_{i \in \mathcal{N}} \left[ \frac{c}{2} \| \mathbf{A}_i \mathbf{z}_i^t \|^2 + \frac{c}{2} \| \mathbf{z}_i^t - \mathbf{z}_i^{t-1} \|_{\mathbf{B}_i^T \mathbf{B}_i}^2 + \frac{1}{2} \| \mathbf{D}_i (\mathbf{w}_i^{t} - \mathbf{w}_i^{t-1}) \|^2 - \frac{c}{2} \| \mathbf{A}_i(\mathbf{z}_i^{t+1} - \mathbf{z}_i^t )\|^2 \right].\tag{B.23}
\end{align*}
Next, we use the dual residual term $\sum_{i \in \mathcal{N}} \|\mathbf{A}_i (\mathbf{z}_i^{t+1} - \mathbf{z}_i^t)\|^2$ to bound $\sum_{i \in \mathcal{N}} \|\mathbf{z}_i^t - \tilde{\mathbf{z}}_i^{t+1}\|_0$. Recalling the definitions of $\mathbf{z}_i$ in equation (\ref{fz}) and $\mathbf{A}_i$ in equation (\ref{ys1}), we have
\begin{align*}
\mathbf{A}_i (\mathbf{z}_i^{t+1} - \mathbf{z}_i^t) = \text{vec} \Big( &\mathbf{x}_i^{t+1} - \mathbf{x}_i^t - \big((\mathbf{p}_{i,j}^-)^{t+1} - (\mathbf{p}_{i,j}^-)^t\big)\\
+&\sum_{k\in\mathcal{M}}\big(\mathbf{y}_{i, \textcolor{black}{k}}^{t+1} - \mathbf{y}_{i, \textcolor{black}{k}}^t - \big((\mathbf{q}_{i, j, \textcolor{black}{k}}^-)^{t+1} - (\mathbf{q}_{i, j, \textcolor{black}{k}}^-)^t\big)\\
+&\sum_{k\in\mathcal{M}}\big(\mathbf{y}_{i, \textcolor{black}{k}}^{t+1} - \mathbf{y}_{i, \textcolor{black}{k}}^t - \big((\mathbf{q}_{i, j, \textcolor{black}{k}}^+)^{t+1} - (\mathbf{q}_{i, j, \textcolor{black}{k}}^+)^t\big)\big)\big), \, j \in \mathcal{N}_i \Big).
\end{align*}
Let $\mathbf{1}_{(1+2M)N_in} \in \mathbb{R}^{N_i n}$ denote a $(1+2M)N_in$-dimensional column vector of all ones. Then it follows:
\begin{align*}\label{a314}
\mathbf{1}_{(1+2M)N_in}^T \mathbf{A}_i (\mathbf{z}_i^{t+1} - \mathbf{z}_i^t) = \sum_{j \in \mathcal{N}_i} \Big[& \mathbf{x}_i^{t+1} - \mathbf{x}_i^t - ((\mathbf{p}_{i,j}^-)^{t+1} - (\mathbf{p}_{i,j}^-)^t) \\
+&\sum_{k\in\mathcal{M}}\Big(\mathbf{y}_{i, \textcolor{black}{k}}^{t+1} - \mathbf{y}_{i, \textcolor{black}{k}}^t - \big((\mathbf{q}_{i, j, \textcolor{black}{k}}^-)^{t+1} - (\mathbf{q}_{i, j, \textcolor{black}{k}}^-)^t\\
+&\mathbf{y}_{i, \textcolor{black}{k}}^{t+1} - \mathbf{y}_{i, \textcolor{black}{k}}^t - \big((\mathbf{q}_{i, j, \textcolor{black}{k}}^+)^{t+1} - (\mathbf{q}_{i, j, \textcolor{black}{k}}^+)^t\big)\big)\Big)\Big]. \tag{B.24}
\end{align*}
Substituting equations (\ref{28}a)–(\ref{28}f) into equation (\ref{a314}) and summing over nodes $i \in \mathcal{N}$ yields:
\begin{align*}\label{a315}
    \sum_{i \in \mathcal{N}} \mathbf{1}_{(1+2M)N_in}^T \mathbf{A}_i (\mathbf{z}_i^{t+1} - \mathbf{z}_i^t) = \sum_{i \in \mathcal{N}} \sum_{j \in \mathcal{N}_i}& \Big[ \frac{c}{c+1} \big( \tilde{\mathbf{x}}_i^{t+1} - \mathbf{x}_i^t - \big( (\tilde{\mathbf{p}}_{i,j}^-)^{t+1} - (\mathbf{p}_{i,j}^-)^t\big)\\
    &+ \sum_{k\in\mathcal{M}}\Big(\tilde{\mathbf{y}}_{i, \textcolor{black}{k}}^{t+1} - \mathbf{y}_{i, \textcolor{black}{k}}^t - \big( (\tilde{\mathbf{q}}_{i, j, \textcolor{black}{k}}^-)^{t+1} - (\mathbf{q}_{i, j, \textcolor{black}{k}}^-)^t\big)\\
    &+ \tilde{\mathbf{y}}_{i, \textcolor{black}{k}}^{t+1} - \mathbf{y}_{i, \textcolor{black}{k}}^t - \big( (\tilde{\mathbf{q}}_{i, j, \textcolor{black}{k}}^+)^{t+1} - (\mathbf{q}_{i, j, \textcolor{black}{k}}^+)^t\big)\Big) \\
    &+\frac{1}{c+1}\big( \tilde{\mathbf{x}}_i^{t+1} - \mathbf{x}_i^t - \big( (\tilde{\mathbf{p}}_{j,i}^-)^{t+1} - (\mathbf{p}_{i,j}^-)^t\big)\big)\Big]. \tag{B.25}
\end{align*}
For the first three lines on the right-hand side of equation (\ref{a315}), using equation (\ref{ys1}) again, we derive:
\begin{align*}\label{a316}
    \sum_{i \in \mathcal{N}} \sum_{j \in \mathcal{N}_i}& \Big[ \frac{c}{c+1} \big( \tilde{\mathbf{x}}_i^{t+1} - \mathbf{x}_i^t - \big( (\tilde{\mathbf{p}}_{i,j}^-)^{t+1} - (\mathbf{p}_{i,j}^-)^t\big)\\
    &+ \sum_{k\in\mathcal{M}}\Big(\tilde{\mathbf{y}}_{i, \textcolor{black}{k}}^{t+1} - \mathbf{y}_{i, \textcolor{black}{k}}^t - \big( (\tilde{\mathbf{q}}_{i, j, \textcolor{black}{k}}^-)^{t+1} - (\mathbf{q}_{i, j, \textcolor{black}{k}}^-)^t\big)\\
    &+ \tilde{\mathbf{y}}_{i, \textcolor{black}{k}}^{t+1} - \mathbf{y}_{i, \textcolor{black}{k}}^t - \big( (\tilde{\mathbf{q}}_{i, j, \textcolor{black}{k}}^+)^{t+1} - (\mathbf{q}_{i, j, \textcolor{black}{k}}^+)^t\big)\Big) \Big]\\
    =\sum_{i \in \mathcal{N}}&\frac{1}{c+1}\mathbf{1}_{N_i}\mathbf{C}_i\mathbf{A}_i(\mathbf{z}_i^{t+1} - \mathbf{z}_i^t)\tag{B.26}
\end{align*}
where $\mathbf{C}_i=[c\cdot\mathbf{I}_{N_in}, (c+1)\cdot\mathbf{I}_{N_in},...,(c+1)\cdot\mathbf{I}_{N_in}]\in\mathbb{R}^{n\cdot N_i \times n \cdot (1+2M)N_i}$.

For the fourth line on the right-hand side of equation (\ref{a315}), since $\mathbf{z}_i^t$ is the optimal solution to problem (\ref{22}), we have $\mathbf{z}_i^t \in \mathcal{Y}$ for all $t \geq 1$, and thus:
\begin{align*}\label{a317}
    (\mathbf{p}_{i,j}^-)^{t} = (\mathbf{p}_{j,i}^+)^t, \quad \forall i \in \mathcal{N}, \, j \in \mathcal{N}_i. \tag{B.27}
\end{align*}
Combining equation (\ref{a317}), the fourth line on the right-hand side of equation (\ref{a315}) can be written as:
\begin{align*}\label{a318}
    &\sum_{i \in \mathcal{N}} \sum_{j \in \mathcal{N}_i} \frac{1}{c+1} \left( \tilde{\mathbf{x}}_i^{t+1} - \mathbf{x}_i^t - ((\tilde{\mathbf{p}}_{j,i}^+)^{t+1} - (\mathbf{p}_{i,j}^-)^t) \right)\\
    =& \sum_{i \in \mathcal{N}} \sum_{j \in \mathcal{N}_i} \frac{1}{c+1} \left( \tilde{\mathbf{x}}_i^{t+1} - \mathbf{x}_i^t - ((\tilde{\mathbf{p}}_{j,i}^+)^{t+1} - (\mathbf{p}_{j,i}^+)^t) \right)\\
    =& \sum_{i \in \mathcal{N}}\frac{1}{c+1}\mathbf{1}_{N_in}^T\mathbf{M}_i\mathbf{H}_i(\tilde{\mathbf{z}}_i^{t+1}-\mathbf{z}_i^t). \tag{B.28}
\end{align*}
where $\mathbf{M}_i=[\mathbf{I}_{N_in},\mathbf{0}_{N_in},...,\mathbf{0}_{N_in}]\in\mathbb{R}^{n \cdot N_i\times n \cdot (N_i+M)}$. The last equality holds based on the expressions for matrix $\mathbf{H}_i$ in equation (\ref{ai}) and the definitions of vectors $\mathbf{z}_i^{t+1}, \mathbf{z}_i^t$.

Substituting equations (\ref{a316}) and (\ref{a318}) into equation (\ref{a315}), we obtain:
\begin{align*}\label{a319}
\left\| \sum_{i \in \mathcal{N}} \mathbf{1}_{(1+2M)N_in}^T \mathbf{A}_i (\mathbf{z}_i^{t+1} - \mathbf{z}_i^t) \right\|^2 = \frac{1}{(c+1)^2} \left\| \sum_{i \in \mathcal{N}} \mathbf{1}_{N_in}^T\left[ \mathbf{C}_i \mathbf{A}_i+ \mathbf{M}_i\mathbf{H}_i \right] (\tilde{\mathbf{z}}_i^{t+1} - \mathbf{z}_i^t) \right\|^2. \tag{B.29}
\end{align*}

Let $\mathbf{A} = \text{diag}(\mathbf{A}_i, i \in \mathcal{N})$ denote the block-diagonal matrix with blocks $\mathbf{A}_i$, and let  $\mathbf{CA} = \text{diag}(\mathbf{C}_i\mathbf{A}_i, i \in \mathcal{N})$, $\mathbf{MH}= \text{diag}(\mathbf{M}_i\mathbf{H}_i, i \in \mathcal{N})$, $\mathbf{z}^{t+1} = \text{vec}(\mathbf{z}_i^{t+1}, i \in \mathcal{N})$. Then equation (\ref{a319}) can be compactly written as:
\begin{align*}\label{a320}
\| \mathbf{1}_{(1+2M)N_{\text{sum}} n}^T \mathbf{A}(\mathbf{z}^{t+1} - \mathbf{z}^t) \|^2 = \frac{1}{(c+1)^2}\|\mathbf{1}_{N_{sum}n}^T[\mathbf{CA}+ \mathbf{MH}] (\mathbf{z}^{t+1} - \mathbf{z}^t)  \|^2. \tag{B.30}
\end{align*}
where $N_{\text{sum}} = \sum_{i \in \mathcal{N}} N_i$. According to the Cauchy–Schwarz inequality, the left-hand side of equation (\ref{a320}) can be lower bounded by:
\begin{align*}\label{a321}
    3N_{\text{sum}}n \| \mathbf{A} (\mathbf{z}^{t+1} - \mathbf{z}^t) \|^2 \geq \| \mathbf{1}_{(1+2M)N_{\text{sum}}n}^T A (\mathbf{z}^{t+1} - \mathbf{z}^t) \|^2. \tag{B.31}
\end{align*}
For the right-hand side of equation (\ref{a320}), we have:
\begin{align*}\label{a322}
    \|\mathbf{1}_{N_{sum}n}^T[\mathbf{CA}+ \mathbf{MH}] (\mathbf{z}^{t+1} - \mathbf{z}^t)  \|^2 \geq \tilde{\tau}_{\min} \| \tilde{\mathbf{z}}^{t+1} - \mathbf{z}^t \|^2, \tag{B.32}
\end{align*}
where $\tilde{\tau}_{\min}$ is defined as the smallest nonzero eigenvalue of $[\mathbf{1}_{N_{sum}n}^T(\mathbf{CA}+ \mathbf{MH})]^T [\mathbf{1}_{N_{sum}n}^T(\mathbf{CA}+ \mathbf{MH})]$. Based on the definitions of $\mathbf{A}_i, \mathbf{H}_i$ in equations (\ref{ai}) and (\ref{qi}), it can be derived that:
\[
\tilde{\tau}_{\min} = \min\quad\{(c+1)^2[(6M+1)N_i^2+N_i] , i \in \mathcal{N}\}.
\]
Substituting equations (\ref{a321}) and (\ref{a322}) into equation (\ref{a320}), we obtain
\[
\frac{c}{2}\| \mathbf{A} (\mathbf{z}^{t+1} - \mathbf{z}^t) \|^2\geq\frac{c\cdot\tilde{\tau}_{\min} }{2(2M+1)N_{sum}n(c+1)^2}\| \tilde{\mathbf{z}}^{t+1} - \mathbf{z}^t \|^2.\tag{B.33}\label{a323}
\]
Rewriting the above into a per-node expression, we get:
\begin{align*}\label{a324}
    \sum_{i \in \mathcal{N}}\frac{c}{2}\| \mathbf{A}_i (\mathbf{z}^{t+1}_i - \mathbf{z}_i^t) \|^2\geq\sum_{i \in \mathcal{N}}\frac{c\cdot\tilde{\tau}_{\min} }{2(2M+1)N_{\text{sum}}n(c+1)^2}\| \tilde{\mathbf{z}}^{t+1}_i - \mathbf{z}^t_i \|^2.\tag{B.34}
\end{align*}
Finally, substituting equation (\ref{a324}) into equation (\ref{a313}) completes the proof. 
    
\end{proof}
The result of Lemma~\ref{lem12} follows from Lemmas~\ref{lem1}, \ref{lemc1} and \ref{lemc2}, along with the potential function $\mathcal{P}^t$ defined in Equation~(\ref{43}).

\end{proof}

\section{(Proof of Theorem \ref{the3})}\label{sa5}
\begin{proof}
According to~\cite{zhang2023distributed}, we have that the function value \( \mathcal{P}^t \) is lower-bounded, i.e.,
\begin{equation}\label{xiajie}\tag{C.1}
	\exists \, \underline{ \mathcal{P}} > -\infty \quad \text{such that} \quad  \mathcal{P}^t \geq \underline{ \mathcal{P}}, \, \forall t > 0.
\end{equation}

Combining equation (\ref{44}), equation (\ref{xiajie}), and the fact that $\mathbf{U}_i$ is a positive definite diagonal matrix, we obtain:
\[
\lim_{t \to \infty} \mathbf{z}_{i}^{t+1} - \mathbf{z}_{i}^t \to 0, \quad \lim_{t \to \infty} \mathbf{w}_{i}^{t+1} - \mathbf{w}_{i}^t \to 0.
\]
The update step of $\boldsymbol{\lambda}_i^{t+1}$ in equation (\ref{gengxin1}c), together with equations (\ref{a25}), (\ref{a28}), and (\ref{xiajie}), leads to
\[
\lim_{t \to \infty} \mathbf{A}_i\mathbf{z}_i^{t+1} \to 0, \quad \lim_{t \to \infty} \boldsymbol{\lambda}_i^{t+1}-\boldsymbol{\lambda}_i^{t}\to 0. \tag{C.2}
\]
\end{proof}

\section{(Proof of Lemma \ref{lem15})}\label{sa4}
\begin{proof}

Suppose $\mathcal{G}(\mathbf{z}^t, \mathbf{w}^t, \mathbf{\lambda}^t) = 0$, then for any node $i \in \mathcal{N}$,
\begin{align}
    \mathbf{z}_i^t - \mathrm{proj}_{\mathcal{X}, \mathcal{Y}} \left(\mathbf{z}_i^t - \left(\nabla_{\mathbf{z}_i} G_i(\mathbf{z}_i^t, \mathbf{w}_i^t) + \mathbf{A}_i^T \boldsymbol{\lambda}_i^t \right) \right) &= 0, \tag{D.1} \label{a41}\\
    \mathbf{A}_i \mathbf{z}_i^t &= 0, \tag{D.2} \label{a42}\\
    \mathbf{w}_i^t - \mathbf{w}_i^{t-1} &= 0. \tag{D.3}\label{a43}
\end{align}

According to the definition of the projection operator, we have
\begin{align*}\label{a44}
    &\mathrm{proj}_{\mathcal{X}, \mathcal{Y}} \left(\mathbf{z}_i^t - \left(\nabla_{\mathbf{z}_i} G_i(\mathbf{z}_i^t, \mathbf{w}_i^t) + \mathbf{A}_i^T \boldsymbol{\lambda}_i^t \right) \right)\\
    =& \arg\min_{\substack{\mathbf{z} \in \mathcal{Y} \\ \mathbf{z} \in \mathcal{X}}} \left\| \mathbf{z}_i - \mathbf{z}_i^t + \left( \nabla_{\mathbf{z}_i} G_i(\mathbf{z}_i^t, \mathbf{w}_i^t) + \mathbf{A}_i^T \boldsymbol{\lambda}_i^t \right) \right\|^2. \tag{D.4}
\end{align*}
Combining the above with equation (\ref{a41}), we get that $\mathbf{z}_i^t$ is the optimal solution of problem (\ref{a44}). Therefore, using the first-order optimality condition of problem (\ref{a44}), we obtain
\begin{align}\label{a45}
    \langle \nabla_{\mathbf{z}_i} G_i(\mathbf{z}_i^t, \mathbf{w}_i^t) + \mathbf{A}_i^T \boldsymbol{\lambda}_i^t, \mathbf{x} - \mathbf{z}_i^t \rangle &\geq 0, \quad \forall \mathbf{x} \in \mathcal{X}, \mathbf{x} \in \mathcal{Y}. \tag{D.5}
\end{align}
Since with fixed $\mathbf{w}_i$, the function $G_i(\mathbf{z}_i, \mathbf{w}_i) + \langle\boldsymbol{\lambda}_i, \mathbf{A}_i \mathbf{z}_i\rangle$ is convex in $\mathbf{z}_i$ over sets $\mathcal{X}$ and $\mathcal{Y}$, equation (\ref{a45}) implies that $\mathbf{z}_i^t$ is the optimal solution of minimizing $G_i(\mathbf{z}_i, \mathbf{w}_i) + \langle\boldsymbol{\lambda}_i, \mathbf{A}_i \mathbf{z}_i\rangle$ over $\mathcal{X}$ and $\mathcal{Y}$, i.e.,
\begin{align*}
    \mathbf{z}_i^t \in \arg \min_{\substack{\mathbf{z} \in \mathcal{Y} \\ \mathbf{z} \in \mathcal{X}}} G_i(\mathbf{z}_i, \mathbf{w}_i^t) + \langle\boldsymbol{\lambda}_i^t, \mathbf{A}_i \mathbf{z}_i\rangle.
\end{align*}Since the objective function is separable across nodes $i$, $\mathbf{z}^t$ must satisfy the Karush-Kuhn-Tucker conditions:
\begin{align}\label{a46}
    \mathbf{z}^t \in \arg \min_{\substack{\mathbf{z} \in \mathcal{Y} \\ \mathbf{z} \in \mathcal{X}}} \sum_{i \in \mathcal{N}} G_i(\mathbf{z}_i, \mathbf{w}_i^t) + \langle\boldsymbol{\lambda}_i^t, \mathbf{A}_i \mathbf{z}_i\rangle. \tag{D.6}
\end{align}
Similarly, since the update form of variable $\mathbf{w}_i$ is
\begin{align*}
    \mathbf{w}_i^t = \mathrm{proj}_{\mathcal{B}^{\textcolor{black}{N_i+M}}} \left(\mathbf{w}_i^{t-1} - \nabla_{\mathbf{w}_i} G_i(\mathbf{z}_i^t, \mathbf{w}_i^{t}) \right),
\end{align*}
Combining it with equation (\ref{a43}), we obtain
\begin{align*}
    0 = \| \mathbf{w}_i^{t-1} - \mathbf{w}_i^t \|^2 = \| \mathbf{w}_i^{t-1} - \mathrm{proj}_{\mathcal{B}^{\textcolor{black}{N_i+M}}} \left(\mathbf{w}_i^{t-1} - \nabla_{\mathbf{w}_i} G_i(\mathbf{z}_i^t, \mathbf{w}_i^{t}) \right) \|^2.
\end{align*}
Therefore, $\mathbf{w}_i^{t-1}$ is the optimal solution of the following problem:
\begin{align}\label{a47}
    \min_{\mathbf{w}_i \in \mathcal{B}^{\textcolor{black}{N_i+M}}} \| \mathbf{w}_i - \mathbf{w}_i^{t-1} + \nabla_{\mathbf{w}_i} G_i(\mathbf{z}_i^t, \mathbf{w}_i^{t}) \|^2. \tag{D.7}
\end{align}
By the optimality condition of problem (\ref{a47}) and equation (\ref{a43}), we have
\begin{align}
    0 \leq \langle \nabla_{\mathbf{w}_i} G_i(\mathbf{z}_i^t, \mathbf{w}_i^t), \mathbf{w}_i - \mathbf{w}_i^t+ \mathbf{w}_i^t - \mathbf{w}_i^{t-1} \rangle = \langle \nabla_{\mathbf{w}_i} G_i(\mathbf{z}_i^t, \mathbf{w}_i^t), \mathbf{w}_i - \mathbf{w}_i^t \rangle, \quad \forall \mathbf{w}_i \in \mathcal{B}^{\textcolor{black}{N_i+M}}. \tag{D.8}
\end{align}
Since $G_i(\mathbf{z}_i, \mathbf{w}_i)$ is convex in $\mathbf{w}_i$ when $\mathbf{z}_i$ is fixed, this implies
\begin{align}
    \mathbf{w}_i^t \in \arg\min_{\mathbf{w}_i \in \mathcal{B}^{\textcolor{black}{N_i+M}}} G_i(\mathbf{z}_i^t, \mathbf{w}_i^t), \tag{D.9}
\end{align}
which can also be written as
\begin{align}\label{a410}
    0 \in \nabla_{\mathbf{w}_i} G_i(\mathbf{z}_i^t, \mathbf{w}_i^t) + \delta_{\mathcal{B}^{\textcolor{black}{N_i+M}}}(\mathbf{w}_i^t). \tag{D.10}
\end{align}

\textbf{Combining equations (\ref{a46}), (\ref{a410}), and (\ref{a42}), we conclude that $(\mathbf{z}^t, \mathbf{w}^t, \mathbf{\lambda}^t)$ satisfies the Karush-Kuhn-Tucker condition (\ref{49}).}

Next, we prove that if $(\mathbf{z}^t, \mathbf{w}^t, \mathbf{\lambda}^t)$ is a KKT point of problem (\ref{jin}), i.e., satisfies condition (\ref{49}), then $(\mathbf{x}^t, \mathbf{y}_1^t)$ is a critical point of the non-convex problem (\ref{ss}), namely,
\begin{align}
    \mathbf{0} &= \sum_{j \in \mathcal{N}_i} \left(\mathbf{x}_i^t - \mathbf{x}_j^t - d_{i,j} \mathbf{v}_{i,j}^t + \mathbf{x}_i^t - \mathbf{x}_j^t + d_{j,i} \mathbf{v}_{j,i}^t\right) +\sum_{\textcolor{black}{k \in \mathcal{M}}} (\mathbf{x}_i^t - \mathbf{y}_{1,k}^t- r_{i, \textcolor{black}{k}} \mathbf{u}_{i, \textcolor{black}{k}}^t), \quad \forall i \in (\mathcal{N} / \mathcal{A}), \tag{D.11} \label{a411}\\
    \mathbf{0} &= \sum_{i \in \mathcal{N}/\mathcal{A}}\sum_{\textcolor{black}{k \in \mathcal{M}}}(\mathbf{y}_{i, \textcolor{black}{k}}^t - \mathbf{x}_i^t + r_{i, \textcolor{black}{k}} \mathbf{u}_{i, \textcolor{black}{k}}^t)
+ \sum_{i \in \mathcal{A}}\sum_{\textcolor{black}{k \in \mathcal{M}}}(\mathbf{y}_{i, \textcolor{black}{k}}^t - \mathbf{a}_i^t + r_{i, \textcolor{black}{k}} \mathbf{u}_{ i,k}^t), \tag{D.12} \label{a412}\\
    \mathbf{x}_l^t &= \mathbf{a}_l, \quad \forall l \in \mathcal{A}. \tag{D.13}\label{a413}
\end{align}

First, since $\mathbf{z}^t$ is the optimal solution of (\ref{49}a), we have $\mathbf{z}^t \in \mathcal{X}$, and thus equation (\ref{a413}) holds.

Next, we prove that $(\mathbf{x}^t, \mathbf{y}_1^t)$ satisfies equations (\ref{a411}) and (\ref{a412}). Since $\mathbf{z}^t$ is the optimal solution of (\ref{49}a), we have $\mathbf{z}^t \in \mathcal{Y}$. Combining this with the KKT condition (\ref{49}c), equation (\ref{a411}) holds, and we can further derive:
\begin{align}\label{a414}
    (\mathbf{p}_{i,j}^+)^t = \mathbf{x}_j^t = (\mathbf{p}_{j,i}^-)^t, \quad \forall i \in \mathcal{N}, \, \forall j \in \mathcal{N}_i. \tag{D.14}
\end{align}
Using the definition of $G_i(\mathbf{z}_i^t, \mathbf{w}_i^t)$ and equation (\ref{a414}), we get:
\begin{align}\label{a415}
\nabla_{\mathbf{w}_i} G_i(\mathbf{z}_i^t, \mathbf{w}_i^t) &= -\mathbf{D}_i \mathbf{H}_i \mathbf{z}_i^t\nonumber \\
&= -\text{vec} \big(-d_{i,j}(\mathbf{x}_i^t - (\mathbf{z}_{i,j}^+)^t),\quad j \in \mathcal{N}_i\quad,-r_{i, \textcolor{black}{k}}(\mathbf{x}_i^t -\mathbf{y}_{i, \textcolor{black}{k}}^t),\quad k\in\mathcal{M} \big)\nonumber \\
&= -\text{vec} \big(-d_{i,j}(\mathbf{x}_i^t - \mathbf{x}_j^t),\quad j \in \mathcal{N}_i \quad,-r_{i, \textcolor{black}{k}}(\mathbf{x}_i^t -\mathbf{y}_{i, \textcolor{black}{k}}^t),\quad k\in\mathcal{M}\big), \tag{D.15}
\end{align}

Recall the following definition:
\begin{align}\label{a417}
G_i(\mathbf{z}_i, \mathbf{w}_i) &= \frac{1}{2}\|\mathbf{H}_i \mathbf{z}_i \|^2 - \mathbf{w}_i^T \mathbf{D}_i \mathbf{H}_i \mathbf{z}_i, \quad \mathcal{X} := \{ \mathbf{z} \mid \mathbf{E}_i \mathbf{z}_i =\mathbf{a}_i, \, \forall \, i \in \mathcal{A} \},\tag{D.17}\end{align}
where $\mathbf{E}_i := \begin{bmatrix} 1, & \mathbf{0}_{(2N_i+1)M+2N_i}^T \end{bmatrix} \otimes \mathbf{I}_n \in \mathbb{R}^{n \times (2N_i+1)(1+M)n}.$

Substituting equations (\ref{a417}) into (\ref{49}a), the optimization problem (\ref{49}a) is rewritten as
\begin{align*}\label{a418}
&\arg \min_{\substack{\mathbf{z}_i \in \mathcal{Y} \\ \mathbf{z}_i \in \mathcal{X}}} 
\sum_{i \in \mathcal{N}} G_i(\mathbf{z}_i, \mathbf{w}_i) + \langle \boldsymbol{\lambda}_i, \mathbf{A}_i \mathbf{z}_i \rangle \\
=& \arg \min_{\substack{\mathbf{p}_{i,j}^+ = \mathbf{p}_{i,j}^- \\\mathbf{q}_{i, j, \textcolor{black}{k}}^+ = \mathbf{q}_{i, j, \textcolor{black}{k}}^-, \forall k\in\mathcal{M}\\ \mathbf{x}_i = \mathbf{a}_i, \forall i \in \mathcal{A}}} 
\sum_{i \in \mathcal{N}} \bigg( 
\sum_{j \in \mathcal{N}_i} \bigg( 
\frac{1}{2} \|\mathbf{x}_i - \mathbf{p}_{i,j}^+ \|^2 
- d_{i,j} \mathbf{v}_{i,j}^T (\mathbf{x}_i - \mathbf{p}_{i,j}^+) \\
& \quad + {\mathbf{\lambda}_{i,j}^0}^T (\mathbf{x}_i - \mathbf{p}_{i,j}^-) 
+\sum_{k\in\mathcal{M}}\big( {\mathbf{\lambda}_{i,j}^{2k-1}}^T (\mathbf{y}_{i, \textcolor{black}{k}} - \mathbf{q}_{i, j, \textcolor{black}{k}}^-) 
+ {\mathbf{\lambda}_{i,j}^{2k}}^T (\mathbf{y}_{i, \textcolor{black}{k}} - \mathbf{q}_{i, j, \textcolor{black}{k}}^+) \big)\bigg) \\
& \quad + \sum_{k\in\mathcal{M}_i}\big(\frac{1}{2} \|\mathbf{x}_i - \mathbf{y}_{i, \textcolor{black}{k}} \|^2 
- r_{i, \textcolor{black}{k}} \mathbf{u}_{i, \textcolor{black}{k}}^T (\mathbf{x}_i - \mathbf{y}_{i, \textcolor{black}{k}}) \big)\bigg) \\
=& \arg \min_{\substack{\mathbf{p}_{i,j}^+ = \mathbf{p}_{i,j}^-\\\mathbf{q}_{i, j, \textcolor{black}{k}}^+ = \mathbf{q}_{i, j, \textcolor{black}{k}}^-\\ \mathbf{x}_i = \mathbf{a}_i, i \in \mathcal{A}}} 
\sum_{i \in (\mathcal{N}/\mathcal{A})} \bigg( 
\sum_{j \in \mathcal{N}_i} \bigg( 
\frac{1}{2} \|\mathbf{x}_i - \mathbf{p}_{j,i}^- \|^2 
- d_{i,j} \mathbf{v}_{i,j}^T (\mathbf{x}_i - \mathbf{p}_{j,i}^-) \\
& \quad + {\mathbf{\lambda}_{i,j}^0}^T (\mathbf{x}_i - \mathbf{p}_{i,j}^-) 
+ \sum_{k\in\mathcal{M}}\big( {\mathbf{\lambda}_{i,j}^{2k-1}}^T (\mathbf{y}_{i, \textcolor{black}{k}} - \mathbf{q}_{i, j, \textcolor{black}{k}}^-) 
+ {\mathbf{\lambda}_{i,j}^{2k}}^T (\mathbf{y}_{i, \textcolor{black}{k}} - \mathbf{q}_{j, i, \textcolor{black}{k}}^-) \big)\bigg) \\
& \quad +\sum_{k\in\mathcal{M}_i}\big( \frac{1}{2} \|\mathbf{x}_i - \mathbf{y}_{i, \textcolor{black}{k}} \|^2 
- r_{i, \textcolor{black}{k}} \mathbf{u}_{i, \textcolor{black}{k}}^T (\mathbf{x}_i - \mathbf{y}_{i, \textcolor{black}{k}})\big) \bigg)\\
& \quad +\sum_{i \in \mathcal{A}} \bigg( 
\sum_{j \in \mathcal{N}_i} \bigg( 
\frac{1}{2} \|\mathbf{a}_i - \mathbf{p}_{j,i}^- \|^2 
- d_{i,j} \mathbf{v}_{i,j}^T (\mathbf{a}_i - \mathbf{p}_{j,i}^-) \\
& \quad + {\mathbf{\lambda}_{i,j}^0}^T (\mathbf{a}_i - \mathbf{p}_{i,j}^-) 
+ \sum_{k\in\mathcal{M}}\big( {\mathbf{\lambda}_{i,j}^{2k-1}}^T (\mathbf{y}_{i, \textcolor{black}{k}} - \mathbf{q}_{i, j, \textcolor{black}{k}}^-) 
+ {\mathbf{\lambda}_{i,j}^{2k}}^T (\mathbf{y}_{i, \textcolor{black}{k}} - \mathbf{q}_{j, i, \textcolor{black}{k}}^-) \big)\bigg) \\
& \quad + \sum_{k\in\mathcal{M}_i}\big(\frac{1}{2} \|\mathbf{a}_i - \mathbf{y}_{i, \textcolor{black}{k}} \|^2 
- r_{i, \textcolor{black}{k}} \mathbf{u}_{i, \textcolor{black}{k}}^T (\mathbf{a}_i - \mathbf{y}_{i, \textcolor{black}{k}})\big) \bigg). \tag{D.18}
\end{align*}
The final equality above holds by substituting the constraints $\mathbf{p}_{i,j}^+=\mathbf{p}_{j,i}^-, \mathbf{q}_{i, j, \textcolor{black}{k}}^+=\mathbf{q}_{j, i, \textcolor{black}{k}}^-$ for all $i \in \mathcal{N}, j \in \mathcal{N}_i$ and $\mathbf{x}_i=\mathbf{a}_i, \forall i\in\mathcal{A}$ into the objective function. Using the KKT condition (\ref{49}a), we conclude that $\mathbf{z}_i^t$ is the optimal solution of problem (\ref{a418}). From the optimality conditions of problem (\ref{a418}) with respect to $\mathbf{x}_i$ and $\mathbf{p}_{i,j}^-$, we obtain:
\begin{align*}
\mathbf{0} &= \sum_{j \in \mathcal{N}_i} \big(\mathbf{x}_i^t - (\mathbf{p}_{j,i}^-)^t - d_{i,j} \mathbf{v}_{i,j}^t + {\mathbf{\lambda}_{i,j}^0}^t\big)+\sum_{k\in\mathcal{M}_i}\big( \mathbf{x}_i^t-\mathbf{y}_{i, \textcolor{black}{k}}^t-r_{i, \textcolor{black}{k}}\mathbf{u}_{i, \textcolor{black}{k}}^t\big), \quad \forall i \in (\mathcal{N}/\mathcal{A}),\\
\mathbf{0} &=\sum_{j\in \mathcal{N}_i}\big((\mathbf{p}_{i,j}^-)^t - \mathbf{x}_j^t + d_{j,i}\mathbf{v}_{j,i}^t- {\mathbf{\lambda}_{i,j}^0}^t\big), \quad \forall i \in (\mathcal{N}/\mathcal{A})
\end{align*}
Summing the above equations yields:
\begin{equation}\label{a419}
\mathbf{0} = \sum_{j \in \mathcal{N}_i} \left( \mathbf{x}_i^t - (\mathbf{p}_{j,i}^-)^t - d_{i,j} \mathbf{v}_{i,j}^t + (\mathbf{p}_{i,j}^-)^t - \mathbf{x}_j^t + d_{j,i} \mathbf{v}_{j,i}^t \right)+\sum_{k\in\mathcal{M}_i}\big(\mathbf{x}_i^t-\mathbf{y}_{i, \textcolor{black}{k}}^t-r_{i, \textcolor{black}{k}}\mathbf{u}_{i, \textcolor{black}{k}}^t\big), \quad \forall i \in (\mathcal{N}/\mathcal{A}). \tag{D.19}
\end{equation}
Using the KKT condition (\ref{49}c), we have $(\mathbf{z}_{i,j}^-)^t = \mathbf{x}_i^t, \, (\mathbf{z}_{j,i}^-)^t = \mathbf{x}_j^t$. Therefore, equation (\ref{a419}) can be rewritten as:
\begin{equation}\label{a420}
\mathbf{0} = \sum_{j \in \mathcal{N}_i} \left( \mathbf{x}_i^t - \mathbf{x}_j^t - d_{i,j} \mathbf{v}_{i,j}^t + \mathbf{x}_i^t - \mathbf{x}_j^t + d_{j,i} \mathbf{v}_{j,i}^t \right)+\sum_{k\in\mathcal{M}_i}\big(\mathbf{x}_i^t-\mathbf{y}_{i, \textcolor{black}{k}}^t-r_{i, \textcolor{black}{k}}\mathbf{u}_{i, \textcolor{black}{k}}^t\big), \quad \forall i \in (\mathcal{N}/\mathcal{A}). \tag{D.20}
\end{equation}
This indicates that $(\mathbf{x}^t, \mathbf{y}_i^t)$ satisfies equation (\ref{a411}).

Similarly, using the optimality conditions of problem (\ref{a418}) with respect to variables $\mathbf{y}_{i, \textcolor{black}{k}}$ and $\mathbf{q}_{i, j, \textcolor{black}{k}}^-$, we obtain:
\begin{align*}
&\begin{cases}
\mathbf{0} &= \sum_{j \in \mathcal{N}_i} \sum_{k\in\mathcal{M}} \big( {\mathbf{\lambda}_{i,j}^{2k-1}}^t + {\mathbf{\lambda}_{i,j}^{2k}}^t \big) +\sum_{k\in\mathcal{M}_i} \big(\mathbf{y}_{i, \textcolor{black}{k}}^t - \mathbf{x}_i^t + r_{i, \textcolor{black}{k}} \mathbf{u}_{i, \textcolor{black}{k}}^t\big), \quad \forall i \in (\mathcal{N} / \mathcal{A}), \\
\mathbf{0} &= \sum_{j \in \mathcal{N}_i} \big( -{\mathbf{\lambda}_{i,j}^{2k-1}}^t - {\mathbf{\lambda}_{j,i}^{2k}}^t \big), \quad \forall i \in (\mathcal{N} / \mathcal{A}),
\end{cases}\\
&\begin{cases}
\mathbf{0} &= \sum_{j \in \mathcal{N}_i} \sum_{k\in\mathcal{M}} \big( {\mathbf{\lambda}_{i,j}^{2k-1}}^t + {\mathbf{\lambda}_{i,j}^{2k}}^t \big) + \sum_{k\in\mathcal{M}_i} \big(\mathbf{y}_{i, \textcolor{black}{k}}^t - \mathbf{a}_i^t + r_{i, \textcolor{black}{k}} \mathbf{u}_{i, \textcolor{black}{k}}^t\big), \quad \forall i \in \mathcal{A}, \\
\mathbf{0} &= \sum_{j \in \mathcal{N}_i} \big( {-\mathbf{\lambda}_{i,j}^{2k-1}}^t - {\mathbf{\lambda}_{j,i}^{2k}}^t \big), \quad \forall i \in \mathcal{A}.
\end{cases}
\end{align*}
Summing the above equations, we have:
\begin{align*}
\mathbf{0} &= \sum_{j \in \mathcal{N}_i} \sum_{k\in\mathcal{M}} \textcolor{black}{\left( {\mathbf{\lambda}_{i,j}^{2k}}^T- {\mathbf{\lambda}_{j,i}^{2k}}^T \right) }
+ \sum_{k\in\mathcal{M}_i} \big(\mathbf{y}_{i, \textcolor{black}{k}}^t - \mathbf{x}_i^t + r_{i, \textcolor{black}{k}} \mathbf{u}_{i }^t\big) , \quad \forall i \in (\mathcal{N} / \mathcal{A}), \\
\mathbf{0} &= \sum_{j \in \mathcal{N}_i} \sum_{k\in\mathcal{M}} \big(\textcolor{black}{\left( {\mathbf{\lambda}_{i,j}^{2k}}^T- {\mathbf{\lambda}_{j,i}^{2k}}^T \right) }
+\sum_{k\in\mathcal{M}_i} \big( \mathbf{y}_{i, \textcolor{black}{k}}^t - \mathbf{a}_i^t + r_{i, \textcolor{black}{k}} \mathbf{u}_{i, \textcolor{black}{k}}^t\big), \quad \forall i \in \mathcal{A} , 
\tag{D.21}
\end{align*}
Therefore,
\begin{align*}
    \mathbf{0} = \sum_{i \in \mathcal{N}} \sum_{j \in \mathcal{N}_i}\sum_{k\in\mathcal{M}} \left( \mathbf{\lambda}_{i,j}^{2k,t} - \mathbf{\lambda}_{j,i}^{2k,t} \right) 
+ \sum_{i \in \mathcal{N}/\mathcal{A}}\sum_{k\in\mathcal{M}_i} (\mathbf{y}_{i, \textcolor{black}{k}}^t - \mathbf{x}_i^t + r_{i, \textcolor{black}{k}} \mathbf{u}_{ i,k}^t)
+ \sum_{i \in \mathcal{A}}\sum_{k\in\mathcal{M}_i} (\mathbf{y}_{i, \textcolor{black}{k}}^t - \mathbf{a}_i^t + r_{i, \textcolor{black}{k}} \mathbf{u}_{ i,k}^t),
\end{align*}
By symmetry, we easily obtain:
\begin{align*}
     \mathbf{0} = \sum_{i \in \mathcal{N}} \sum_{j \in \mathcal{N}_i}\sum_{k\in\mathcal{M}}\left( \mathbf{\lambda}_{i,j}^{2k,t} - \mathbf{\lambda}_{j,i}^{2k,t} \right) 
\end{align*}
Using the KKT condition (\ref{49}a), the constraint $\mathbf{y}_{i, \textcolor{black}{k}}=\mathbf{y}_{j,k},\forall(i,j)\in\mathcal{E}$ holds. Given the sensor network is connected, it follows that all nodes reach consensus on the target position estimate, i.e., $\mathbf{y}_{i, \textcolor{black}{k}}^t=\mathbf{y}_{1,k}^t$ for all $i\in\mathcal{N}$. Thus:
\begin{align*}
    \mathbf{0} =  \sum_{i \in \mathcal{N}/\mathcal{A}}\sum_{k\in\mathcal{M}_i}(\mathbf{y}_{1,k}^t - \mathbf{x}_i^t + r_{i, \textcolor{black}{k}} \mathbf{u}_{ i,k}^t)
+ \sum_{i \in \mathcal{A}}\sum_{k\in\mathcal{M}_i}(\mathbf{y}_{1,k}^t - \mathbf{a}_i^t + r_{i, \textcolor{black}{k}} \mathbf{u}_{ i,k}^t),
\end{align*}
This indicates that $(\mathbf{x}^t, \mathbf{y}_1^t)$ satisfies equation (\ref{a412}), where $\mathbf{y}_1^t=\text{vec}(\mathbf{y}_{1,k}, k\in\mathcal{M})$.
Hence, we obtain that \( (\mathbf{x}^t, \mathbf{y}_1^t, \mathbf{w}^t) \) is a critical point of problem~(\ref{ss}). Consequently, it follows naturally that \( (\mathbf{x}^t, \mathbf{y}_1^t) \) is a critical point of the original problem~(\ref{yuanwenti}).

\end{proof}


\section{ (Proof of Lemma \ref{lem17})}\label{sa6}
\begin{proof}

According to Lemma \ref{lem12}, the function value $\mathcal{P}^t$ is monotonically decreasing. This result, combined with Lemma \ref{lem15}, implies that for any $t > 1$, there exists $s \in \{1, 2, \ldots, t-1\}$ such that $\mathcal{P}^s - \mathcal{P}^{s+1} \leq \frac{\mathcal{P}^0 - \underline{\mathcal{P}}}{t-1}$. Let
\begin{align*}
    K := (\mathcal{P}^0 - \underline{\mathcal{P}}) \cdot \max \big\{ &\frac{1}{2}, \frac{\kappa_1 - 1}{2}, \frac{c (\kappa_1 - 1)}{2}, \frac{\rho}{4}, \frac{c (\kappa_1 - 6 (2N_{\max} + 1))}{2}, \frac{c \kappa_1 - 6 (N_{\max} + M_{\max}+1)}{2c},\\ & \frac{c \kappa_2\cdot \tilde{\tau}_{\min}}{2 (2M+1)N_{\text{sum}} n (c+1)^2} - \frac{(2N_{\max}+1)c\kappa_1}{2} , \frac{\rho}{4} - d^2_{\max} (\kappa_1 + \kappa_2) \big\}.
\end{align*}
Combining this with inequality (\ref{44}), we obtain
\begin{align*}\label{a61}
    &\| \mathbf{z}_i^{s+1} - \mathbf{z}_i^s \|^2 < \frac{K}{t-1}, \quad \| \mathbf{w}_i^{s+1} - \mathbf{w}_i^s \|^2 < \frac{K}{t-1},\quad
\| \tilde{\mathbf{z}}_i^{s+1} - \mathbf{z}_i^{s+1} - (\tilde{\mathbf{z}}_i^s - \mathbf{z}_i^s) \|^2_{\mathbf{A}_i^T \mathbf{A}_i} < \frac{K}{t-1}, \\ 
&\|\mathbf{H}_i(\tilde{\mathbf{z}}_i^{s+1}-\tilde{\mathbf{z}}_i^{s})-\mathbf{D}_i(\mathbf{w}_i^s-\mathbf{w}_i^{s-1})\|<\frac{K}{t-1},\quad\| \tilde{\mathbf{z}}_i^{s+1} - \mathbf{z}_i^s - (\tilde{\mathbf{z}}_i^s - \mathbf{z}_i^{s-1}) \|^2_{\mathbf{B}_i^T \mathbf{B}_i} < \frac{K}{t-1}.\tag{E.1}
\end{align*}
From problem (\ref{gengxin1}a), we know that
\[
\mathbf{z}_i^{s+1} = \arg \min_\mathbf{z} G_i(\mathbf{z}_i, \mathbf{w}_i^s) + \langle \boldsymbol{\lambda}^s, \mathbf{A}_i \mathbf{z}_i \rangle + \frac{c}{2} \| \mathbf{z}_i - \mathbf{z}_i^s \|^2_{\mathbf{B}_i^T \mathbf{B}_i} + \delta_{\mathcal{X}\times \mathcal{Z}}(\mathbf{z}_i).
\]
Applying the optimality condition, we obtain
\[
\mathbf{0} \in \nabla_{\mathbf{z}_i} G_i(\mathbf{z}_i^{s+1}, \mathbf{w}_i^s) + \mathbf{A}_i^T \boldsymbol{\lambda}_i^s + c \mathbf{B}_i^T \mathbf{B}_i (\mathbf{z}_i^{s+1} - \mathbf{z}_i^s) + \partial\delta_{\mathcal{X}\times \mathcal{Z}}(\mathbf{z}_i^{s+1}).
\]
Define the vector $\boldsymbol{\eta}_1$ as
\[
\boldsymbol{\eta}_1 := \nabla_{\mathbf{z}_i} G_i(\mathbf{z}_i^{s+1}, \mathbf{w}_i^s) + \mathbf{A}_i^T \boldsymbol{\lambda}_i^{s+1} - \nabla_{\mathbf{z}_i} G_i(\mathbf{z}_i^{s+1}, \mathbf{w}_i^s) - \mathbf{A}_i^T \boldsymbol{\lambda}_i^s - c \mathbf{B}_i^T \mathbf{B}_i (\mathbf{z}_i^{s+1} - \mathbf{z}_i^s).
\]
Then, using the triangle inequality of norms, we can derive
\begin{align*}\label{a62}
    \| \boldsymbol{\eta}_1 \| =& \| -\mathbf{H}_i^T \mathbf{D}_i(\mathbf{w}_i^{s+1} - \mathbf{w}_i^s) +  \mathbf{A}_i^T (\boldsymbol{\lambda}_i^{s+1} - \boldsymbol{\lambda}_i^s)  + c \mathbf{B}_i^T \mathbf{B}_i (\mathbf{z}_i^{s+1} - \mathbf{z}_i^s) \|\\
    \leq &d_{\max} \sqrt{N_i + M_i+1} \| \mathbf{w}_i^{s+1} - \mathbf{w}_i^s \| + \sqrt{2N_i + 1} \| \boldsymbol{\lambda}_i^{s+1} - \boldsymbol{\lambda}_i^s \|
\\&+ (1 + c)(2 + N_{\max}) \| \mathbf{z}_i^{s+1} - \mathbf{z}_i^s \|,\tag{E.2}
\end{align*}
The last inequality holds because $\|\mathbf{H}_i^T \mathbf{H}_i\| = N_i + M_i+1$ and $\|\mathbf{A}_i^T \mathbf{A}_i\| = 2N_i + 1$, and based on (\ref{a29}).

Substituting equations (\ref{a61}) and (\ref{a28}) into equation (\ref{a62}) yields
\begin{align*}
    \|\boldsymbol{\eta}_1\| &\leq \sqrt{\frac{K}{t-1}} \left( d_{\max} \sqrt{N_i + M_i+1} + \sqrt{(2N_i + 1)[3(N_i + M_i+1) + 3c^2(4N_i + 2) + 3c(1+c)(1+M_{\max} + N_{\max})]} \right)  \\
    &\le \sqrt{\frac{S_1}{t-1}},
\end{align*}
where 
\[
S_1 := K\left( d_{\max} \sqrt{N_i + M_i+1} + \sqrt{(2N_i + 1)[3(N_i + M_i+1) + 3c^2(4N_i + 2) + 3c(1+c)(1+M_{\max}  + N_{\max})]} \right) > 0.
\]
Recalling the update step of \(\mathbf{w}\), we have
\[
\mathbf{w}_i^{s+1} = \arg\min G_i(\mathbf{z}_i^{s+1}, \mathbf{w}_i) + \frac{\rho}{2} \|\mathbf{w}_i - \mathbf{w}_i^{s}\|^2 + \delta_{\mathbf{B}^{N_i+M}}(\mathbf{w}_i).
\]
By the optimality condition, it follows that
\[
\mathbf{0} \in \nabla_{\mathbf{w}_i} G_i(\mathbf{z}_i^{s+1}, \mathbf{w}_i^{s+1}) + \rho(\mathbf{w}_i^{s+1} - \mathbf{w}_i^{s}) + \partial \delta_{\mathbf{B}^{N_i+M}} (\mathbf{w}_i^{s+1}).
\]
Define the vector \(\boldsymbol{\eta}_2\) as
\[
\boldsymbol{\eta}_2 := \nabla_{\mathbf{w}_i} G_i(\mathbf{z}_i^{s+1}, \mathbf{w}_i^{s+1}) - \nabla_{\mathbf{w}_i} G_i(\mathbf{z}_i^{s+1}, \mathbf{w}_i^{s}) - \rho (\mathbf{w}_i^{s+1} - \mathbf{w}_i^{s})=-\rho(\mathbf{w}_i^{s+1}-\mathbf{w}_i^s),
\]
then we obtain
\[
\boldsymbol{\eta}_2 \in \nabla_{\mathbf{w}_i} G_i(\mathbf{z}_i^{s+1}, \mathbf{w}_i^{s+1}) + \partial \delta_{\mathbf{B}^{N_i+M}} (\mathbf{w}_i^{s+1}).
\]
Combining this with inequality (\ref{a61}), we get
\[
\|\boldsymbol{\eta}_2\| \leq \sqrt{\frac{K}{t-1}} \cdot \rho = \frac{S_2}{\sqrt{t}},
\]
where 
\[
S_2 := K \rho^2 > 0.
\]
Therefore, Lemma \ref{lem17} holds for \(\boldsymbol{\eta}_1, \boldsymbol{\eta}_2\), and \(S := \max\{S_1, S_2\}\), which implies that \((\mathbf{z}^t, \mathbf{w}^t, \boldsymbol{\lambda}^t)\) is an \(S/\sqrt{t-1}\)-solution to problem (\ref{jin}). This shows that the proposed algorithm converges to a Karush-Kuhn-Tucker point of problem (\ref{jin}) at a sublinear rate.
    
\end{proof}

\section{(Proof of Theorem \ref{the1})}\label{sa7}
\begin{proof}

We first prove the first part. Using the definition of the function $G_i$ given in equations (\ref{zwan}) and (\ref{jin}), we have
\begin{align*}\tag{F.1}
\nabla_{\mathbf{z}_i} G_i (\mathbf{z}_i^t, \mathbf{w}_i^t) + \mathbf{A}_i^T \boldsymbol{\lambda}_i^t = \mathbf{H}_i^T \mathbf{H}_i \mathbf{z}_i^t - \mathbf{H}_i^T \mathbf{D}_i \mathbf{w}_i^t + \mathbf{A}_i^T \boldsymbol{\lambda}_i^t 
= -\mathbf{U}_i (\tilde{\mathbf{z}}_i^{t+1} - \mathbf{z}_i^t) - c \mathbf{A}_i^T \mathbf{A}_i \mathbf{z}_i^t.
\end{align*}
Therefore, combining the above equations yields:
\begin{align*}
&\| \mathbf{z}_i^t - \mathrm{proj}_{\mathcal{X}, \mathcal{Y}} \big( \mathbf{z}_i^t - (\nabla_{\mathbf{z}_i} G_i (\mathbf{z}_i^t, \mathbf{w}_i^t) + \mathbf{A}_i^T \boldsymbol{\lambda}_i^t) \big) \| \nonumber\\
=& \| \mathbf{z}_i^t - \mathrm{proj}_{\mathcal{X}, \mathcal{Y}} \big( \tilde{\mathbf{z}}_i^{t+1}) - \mathrm{proj}_{\mathcal{X}, \mathcal{Y}} (\mathbf{z}_i^t + (\mathbf{U}_i (\tilde{\mathbf{z}}_i^{t+1} - \mathbf{z}_i^t) - c \mathbf{A}_i^T \mathbf{A}_i \mathbf{z}_i^t) \big) \|\nonumber\\
\leq& \| \mathbf{z}_i^t - \mathrm{proj}_{\mathcal{X}, \mathcal{Y}} (\tilde{\mathbf{z}}_i^{t+1}) \| + \| (\mathbf{I}_i - \mathbf{U}_i) (\tilde{\mathbf{z}}_i^{t+1} - \mathbf{z}_i^t) - c \mathbf{A}_i^T \mathbf{A}_i \mathbf{z}_i^t \|. \tag{F.2}\label{a72}
\end{align*}
The above inequality follows from the triangle inequality and the non-expansiveness property of the projection operator. Squaring both sides of inequality (\ref{a72}) and applying the Cauchy–Schwarz inequality further gives:
\begin{align*}
    &\| \mathbf{z}_i^t - \mathrm{proj}_{\mathcal{X}, \mathcal{Y}} \big( \mathbf{z}_i^t - (\nabla_{\mathbf{z}_i} G_i (\mathbf{z}_i^t, \mathbf{w}_i^t) + \mathbf{A}_i^T \boldsymbol{\lambda}_i^t) \big) \|^2\\
\leq &3 \| \mathbf{z}_i^t - \mathrm{proj}_{\mathcal{X}, \mathcal{Y}} (\tilde{\mathbf{z}}_i^{t+1}) \|^2 + 3 c^2 \| \mathbf{A}_i^T \mathbf{A}_i \mathbf{z}_i^t \|^2 + 3 \| (\mathbf{I}_i - \mathbf{U}_i) (\tilde{\mathbf{z}}_i^{t+1} - \mathbf{z}_i^t) \|^2. \tag{F.3}\label{a73}
\end{align*}
Now consider computing $\mathrm{proj}_{\mathcal{X}, \mathcal{Y}} (\tilde{\mathbf{z}}_i^{t+1})$. Since
\[
\mathrm{proj}_{\mathcal{X}, \mathcal{Y}} (\tilde{\mathbf{z}}_i^{t+1}) := 
\begin{cases}
\arg \min_\mathbf{z} \frac{1}{2} \sum_{i \in \mathcal{N}} \| \mathbf{z}_i - \tilde{\mathbf{z}}_i^{t+1} \|^2 \\
\text{subject to }\quad z \in \mathcal{X}, z \in \mathcal{Y}.
\end{cases}\tag{F.4}
\]
Following a derivation similar to Remark \ref{remark1}, we can also obtain
\[
\mathrm{proj}_{\mathcal{X}, \mathcal{Y}} (\tilde{\mathbf{z}}_i^{t+1}) = \mathbf{\widetilde{U}}_i (\tilde{\mathbf{z}}_i^{t+1} - \mathbf{z}_i^t) + \tilde{\mathbf{z}}_i^{t+1}, \quad \text{and } \tilde{\mathbf{x}}_i = \mathbf{a}_i, \ \text{if } i \in \mathcal{A}.\tag{F.5}\label{a75}
\]
where
\[
\mathbf{\widetilde{U}}_i :=- \frac{1}{2} \cdot \mathrm{Diag} \left( \left[ 0, (c+1) \cdot \mathbf{1}_{N_i}, \frac{c+1}{c} \cdot \mathbf{1}_{N_i},0,\mathbf{1}_{N_i}, \mathbf{1}_{N_i},...,0,\mathbf{1}_{N_i}, \mathbf{1}_{N_i}  \right] \right) \otimes \mathbf{I}_n. \tag{F.6}\label{a76}
\]
Equation (\ref{a75}), in combination with equation (\ref{28}) and the Cauchy–Schwarz inequality, leads to
\begin{align*}
    &\| \mathbf{z}_i^t - \mathrm{proj}_{\mathcal{X},\mathcal{Y}} (\tilde{\mathbf{z}}_i^{t+1}) \|^2 \nonumber\\
=& \| \mathbf{z}_i^t - \tilde{\mathbf{z}}_i^{t+1} - \mathbf{\widetilde{U}}_i (\tilde{\mathbf{z}}_i^{t+1} - \mathbf{z}_i^t +\mathbf{z}_i^{t}- \mathbf{z}_i^{t+1}) \|^2 \nonumber\\
\leq& 2 \| (\mathbf{I}_i + \mathbf{\widetilde{U}}_i) (\mathbf{z}_i^t - \tilde{\mathbf{z}}_i^{t+1}) \|^2 + 2 \| \mathbf{\widetilde{U}}_i (\mathbf{z}_i^t - \mathbf{z}_i^{t+1} )\|^2. \tag{F.7}\label{a77}
\end{align*}

For the third term on the left-hand side of inequality (\ref{a73}), applying the Cauchy–Schwarz inequality again gives
\[
\| (\mathbf{I}_i - \mathbf{U}_i) (\tilde{\mathbf{z}}_i^{t+1} - \mathbf{z}_i^t) \|^2 \leq \| \mathbf{I}_i - \mathbf{U}_i \|^2 \| \tilde{\mathbf{z}}_i^{t+1} - \mathbf{z}_i^t \|^2. \tag{F.8}\label{a78}
\]
Substituting equations (\ref{a77})–(\ref{a78}) into equation (\ref{a73}) yields
\begin{align*}
&\| \mathbf{z}_i^t - \mathrm{proj}_{\mathcal{X},\mathcal{Y}} (\mathbf{z}_i^t - \nabla_{\mathbf{z}_i} G_i (\mathbf{z}_i^t, \mathbf{w}_i^t) + \mathbf{A}_i^T \boldsymbol{\lambda}_i^t) \|^2\\
\leq& (6 \| \mathbf{I}_i + \mathbf{\widetilde{U}}_i \|^2 + 3 \| \mathbf{I}_i - \mathbf{U}_i \|^2) \| \mathbf{z}_i^t - \tilde{\mathbf{z}}_i^{t+1} \|^2 + 6 \| \mathbf{\widetilde{U}}_i \|^2 \| \mathbf{z}_i^t - \mathbf{z}_i^{t+1} \|^2 + 3c^2 \| \mathbf{A}_i^T \mathbf{A}_i \mathbf{z}_i^t \|^2\\
=& \sigma_1 \| \mathbf{z}_i^t - \tilde{\mathbf{z}}_i^{t+1} \|^2 + \sigma_2 \| \mathbf{z}_i^t - \mathbf{z}_i^{t+1} \|^2 + 3c^2 (2N_{\max} + 1)^2 \| \mathbf{A}_i \mathbf{z}_i^t \|^2, \tag{F.9}\label{a79}
\end{align*}
This equality holds based on the definitions of matrices \( \mathbf{A}_i^T \mathbf{A}_i \), \( \mathbf{U}_i \), and \( \mathbf{\widetilde{U}}_i \) in equations (\ref{24}), (\ref{27}), and (\ref{a76}), respectively, and
\begin{align*}\label{a710}
\sigma_1 &:= 3 \max\left\{4\big( (c+1) N_{\max} +M_{\max} \big)^2-1,(4cN_\text{max})^2\right\} + 6 \max \left\{ (\frac{1-c}{2})^2, (\frac{1-c}{2c})^2,(\frac{1}{2})^2\right\},\\
\sigma_2 &:= \frac{3}{2} \max \left\{ (1+c)^2, \left( 1 + \frac{1}{c} \right)^2 \right\}. \tag{F.10}
\end{align*}
The definition of the function \( \mathcal{G}(\mathbf{z}^t, \mathbf{w}^t, \boldsymbol{\lambda}^t) \) in equation (\ref{48}), combined with equation (\ref{a79}) and the first part of Theorem \ref{the1}, implies that \( \mathcal{G}(\mathbf{z}^t, \mathbf{w}^t, \boldsymbol{\lambda}^t) \to 0 \). Together with Lemma \ref{lem15}, we conclude that when the parameters satisfy equations (\ref{45})–(\ref{47}), the sequence \( \{ (\mathbf{z}_i, \mathbf{w}_i, \boldsymbol{\lambda}_i) \}_{t \geq 1} \) generated by the SP-ADMM-JCNl converges to a Karush-Kuhn-Tucker point of problem (\ref{jin}).

Finally, we prove the second of Theorem \ref{the1} by further bounding the upper limit of \( \mathcal{G}(\mathbf{z},\mathbf{v}, \lambda) \). Applying the Cauchy–Schwarz inequality gives:
\[
\| \mathbf{A}_i \mathbf{z}_i^t \|^2 = \| \mathbf{A}_i (\mathbf{z}_i^t - \mathbf{z}_i^{t+1} + \mathbf{z}_i^{t+1}) \|^2 \leq 2 \| \mathbf{A}_i (\mathbf{z}_i^t - \mathbf{z}_i^{t+1}) \|^2 + 2 \| \mathbf{A}_i \mathbf{z}_i^{t+1} \|^2. \tag{F.11}\label{a711}
\]
Combining equation (\ref{48}) with \textcolor{black}{equation (\ref{a79})} and equation (\ref{a711}), we obtain:
\begin{align*}
\mathcal{G}(\mathbf{z}^t, \mathbf{w}^t, \boldsymbol{\lambda}^t) \leq &\sum_{i \in \mathcal{N}} \left( \sigma_1 \| \tilde{\mathbf{z}}_i^{t+1} - \mathbf{z}_i^t \|^2 + \sigma_2 \| \mathbf{z}_i^t - \mathbf{z}_i^{t+1} \|^2 + \| \mathbf{w}_i^t - \mathbf{w}_i^{t+1} \|^2 \right.\\
&\left. + 2\sigma_3 \| \mathbf{A}_i (\mathbf{z}_i^t - \mathbf{z}_i^{t+1}) \|^2 + 2\sigma_3 \| \mathbf{A}_i \mathbf{z}_i^{t+1} \|^2 \right)\\
\leq &\sum_{i \in \mathcal{N}} \big( \| \tilde{\mathbf{z}}_i^{t+1} - \mathbf{z}_i^t \|^2 + (\sigma_2 + 2\sigma_3 (2N_{\max} + 1)) \| \mathbf{z}_i^t - \mathbf{z}_i^{t+1} \|^2 \\
&+ \| \mathbf{w}_i^t - \mathbf{w}_i^{t+1} \|^2 + \frac{2\sigma_3}{c^2} \| \boldsymbol{\lambda}_i^{t+1} - \boldsymbol{\lambda}_i^t \|^2 \big), \tag{F.12}
\end{align*}
Here, \(\sigma_1\) and \(\sigma_2\) are defined in equation (\ref{a710}), and \(\sigma_3 = 3c^2 (2N_{\max} + 1)^2 + 1\). The last inequality follows from the update step of \(\boldsymbol{\lambda}_i^{t+1}\) and the fact that \(\| \mathbf{A}_i^T \mathbf{A}_i \| = 2N_i + 1\). Matching the upper bound of \(\mathcal{P}^t\) differences from equation (\ref{44}) with the upper bound of \(\mathcal{G}(\mathbf{z}^t, \mathbf{w}^t, \boldsymbol{\lambda}^t)\) in equation (\ref{a812}), we get:
\[
\mathcal{G}(\mathbf{z}^t, \mathbf{w}^t, \boldsymbol{\lambda}^t) \leq \epsilon\left( \mathcal{P}^t - \mathcal{P}^{t+1} \right), \tag{F.13}\label{a713}
\]
where $\epsilon = \frac{\min\{ c, 1, \frac{\rho}{4}, C_0, C_1, \frac{\rho}{4}-d_{\max}^2 (\kappa_1 + \kappa_2)\}}{\max \{ \sigma_1, \sigma_2 + 2\sigma_3 (1+2N_{\max}),1, \frac{2\sigma_3}{c^2}  \}}$,  $C_0 := \frac{\min\{ c\kappa_1 - 6 (N_{\max} + M_{\max}+1), c (c\kappa_1 - 6 (1+c) (N_{\max} + M_{\max}+1)) \}}{\max \{ 3 (N_i + M_i+1), 3c (1+c) (2+N_{\max}) \}}$, and $C_1 := \frac{c \kappa_2\cdot\tilde{\tau}_{\min\quad}}{6 N_{\mathrm{sum}} n(c+1)^2}-\frac{(2N_{max}+1)c\kappa_1}{2}$.

Suppose \(\mathcal{G}(\mathbf{z}^t, \mathbf{w}^t, \boldsymbol{\lambda}^t)\) reaches the lower bound \(\epsilon_1\) for the first time at iteration \(T\). Summing inequality (\ref{a713}) over the first \(T\) iterations and combining it with equation (\ref{xiajie}), we have:
\begin{align*}
\epsilon_1 &\leq \frac{1}{T-1} \sum_{t=1}^T \mathcal{G}(\mathbf{z}^t, \mathbf{w}^t, \boldsymbol{\lambda}^t)\\
&\leq \frac{1}{T-1} \epsilon_0 \left( \mathcal{P}^1 - \mathcal{P}^{T+1} \right)\\
&\leq \frac{1}{T-1} \epsilon_0 \left( \mathcal{P}^1 - \underline{\mathcal{P}} \right) = \frac{\epsilon_2}{T-1}. \tag{F.14}
\end{align*}
Since \(\epsilon_2 := \epsilon_0 ( \mathcal{P}^1 - \underline{\mathcal{P}}) > 0\) is a constant, it follows that \(\mathcal{G}(\mathbf{z}^t, \mathbf{w}^t, \boldsymbol{\lambda}^t)\) converges at a rate of \(O(1/T)\). This completes the proof of Theorem \ref{the1}.     
\end{proof}

\section{(Proof of Theorem \ref{the2})}\label{sa8}
The following derivation is primarily based on Section 6 of~\cite{zhang2023distributed}. To ensure completeness, we include the full proof here. We begin with a concise overview of the proof strategy for Theorem \ref{the2}:

1. First, we show that the potential function $\mathcal{P}$ is lower bounded and possesses a sufficient decrease property;

2. We prove that the sequence $\{\mathbf{s}^t\}_{t \geq 0}$ is bounded. By combining this result with Step 1, we infer that the set of limit points of $\{\mathbf{s}^t\}_{t \geq 0}$ is nonempty and compact;

3. To approximate the subgradient set of $\mathcal{P}^t$, we provide a lower bound on the subgradient in terms of the inter-iteration gap. We also state some key properties of the limit set;

4. Since the potential function $\mathcal{P}^t$ is semi-algebraic, it satisfies the Kurdyka-Łojasiewicz property. Based on the previous steps and this property, we conclude that the sequence $\{\mathbf{s}^t\}_{t \geq 0}$ is a Cauchy sequence and converges to a Karush-Kuhn-Tucker point of problem (\ref{jin}).

The following lemmas correspond to the steps outlined above.

\textbf{Lemma A.8.1 (Sufficient Decrease Condition)}  
Assume that the sequence $\{\mathbf{z}^t, \mathbf{w}^t, \boldsymbol{\lambda}^t\}$ is generated by Algorithm \ref{al1}, the form of $c\mathbf{B}_i^T \mathbf{B}_i$ satisfies (\ref{23}), and the parameters satisfy conditions (\ref{45})–(\ref{47}). Then:
\begin{align*}
\mathcal{P}^{t+1} -\mathcal{P}^t \leq \sum_{i \in \mathcal{N}}\big[ &-\min\{c, 1\} \|\mathbf{z}_i^{t+1} - \mathbf{z}_i^t\|^2 - \frac{\rho}{4} \|\mathbf{w}_i^t - \mathbf{w}_i^{t-1}\|^2 - C_0 \|\boldsymbol{\lambda}_i^{t+1} - \boldsymbol{\lambda}_i^t\|^2 \\
&-C_1\|\tilde{\mathbf{z}}_i^{t+1}-\mathbf{z}_i^{t+1}\|^2 \big].
\end{align*}
Here, $C_0$ and $C_1$ are defined in equation (\ref{a713}), and $d_{\max} := \max\{r_{i, \textcolor{black}{k}},d_{i,j}, i \in \mathcal{N}, j \in \mathcal{N}_i\}$, $N_{\max} := \max\{N_i, i \in \mathcal{N}\}$, $N_{\text{sum}} := \sum_{i \in \mathcal{N}} N_i$, $\tilde{\mathbf{\tau}}_{\min} = \min\{(c+1)^2 N_i^2 + c^2 N_i + N_i, i \in \mathcal{N}\}$.

\textbf{Proof:}  
Substituting equation (\ref{a29}) into equation (\ref{a28}) yields:
\begin{align*}
\|\boldsymbol{\lambda}_i^{t+1} - \boldsymbol{\lambda}_i^t\|^2 
\leq \max \big\{ &3(N_i + M_i+1), \,3c^2(2N_\text{max}+1)\, 3c(1+c)(1+M_{\text{max}} + N_{\text{max}}) \big\} \cdot \big( 
\|\mathbf{Q}_i (\tilde{\mathbf{z}}_i^{t+1} - \tilde{\mathbf{z}}_i^t) - \mathbf{D}_i (\mathbf{w}_i^t - \mathbf{w}_i^{t-1}) \|^2 \\
&+ \|\tilde{\mathbf{z}}_i^{t+1} - \mathbf{z}_i^{t+1} - (\tilde{\mathbf{z}}_i^t - \mathbf{z}_i^t)\|_{\mathbf{A}_i^T \mathbf{A}_i}^2 
+ \|\tilde{\mathbf{z}}_i^{t+1} - \mathbf{z}_i^{t+1} - (\tilde{\mathbf{z}}_i^t - \mathbf{z}_i^t)\|_{\mathbf{B}_i^T \mathbf{B}_i}^2 \big), \quad \forall \, i \in \mathcal{N}. \tag{G.1}\label{a81}
\end{align*}
Combining equation (\ref{42}), equation (\ref{a81}), and equation (\ref{27}) completes the proof. $\hfill \Box$

\textbf{Lemma A.8.2 (Boundedness of the Sequence)}  
Suppose the parameters satisfy conditions (\ref{45})–(\ref{47}), the sensor network is connected, and there exists at least one anchor node. Then the sequence $\{\mathbf{z}^t, \mathbf{w}^t, \boldsymbol{\lambda}^t, \tilde{\mathbf{z}}^t\}$ generated by Algorithm \ref{al1} is bounded.

\textbf{Proof:}  
We establish boundedness of the sequences $\{\mathbf{w}^t\}_{t \geq 1}$, $\{\mathbf{z}^t\}_{t \geq 1}$, $\{\boldsymbol{\lambda}^t\}_{t \geq 1}$, and $\{\tilde{\mathbf{z}}^t\}_{t \geq 1}$ as follows:

1. From Algorithm \ref{al1}, we have $\{\mathbf{w}^t\}_{t \geq 1} \subset \mathcal{B}^{\sum_{i\in\mathcal{N}}(N_i+M)}$, where $\mathcal{B}$ is the unit ball constraint set. Hence, $\{\mathbf{w}^t\}$ is bounded.

2. Using the update rule for $\boldsymbol{\lambda}_i^{t+1}$, we obtain:
\begin{align*}
c \sum_{t=1}^\infty \sum_{i \in \mathcal{N}} \|\mathbf{A}_i \mathbf{z}_i^{t+1}\|^2 &= \sum_{t=1}^\infty \sum_{i \in \mathcal{N}} \|\boldsymbol{\lambda}_i^{t+1} - \boldsymbol{\lambda}_i^t\|^2\\
&\leq \frac{1}{C_0} \sum_{t=1}^\infty (\mathcal{P}^t - \mathcal{P}^{t+1}) \leq \frac{1}{C_0} (\mathcal{P}^1 - \underline{\mathcal{P}}) < \infty. \tag{G.2}
\end{align*}
The first inequality follows from Lemma A.8.1, and the second from equation (\ref{xiajie}). Hence, for all $i \in \mathcal{N}$, the sequence $\{\|\mathbf{A}_i \mathbf{z}_i^{t+1}\|\}_{t\geq1}$ is bounded. Similarly, applying Lemma A.8.1 and equation (\ref{xiajie}) again, we find that $\{\|\mathbf{z}_i^{t+1} - \mathbf{z}_i^t\|\}_{t\geq1}$ is also bounded. 

 \textcolor{black}{We construct a matrix $\mathbf{B}_i$ of the form }
\textcolor{black}{\begin{align*}\tag{G.3}
    \mathbf{B}_i&:=\begin{bmatrix}
        \mathbf{0}_{N_i} & \mathbf{O}_{N_i} & -\mathbf{I}_{N_i} & \mathbf{0}_{N_i} & \mathbf{O}_{N_i} & \mathbf{O}_{N_i} &...&\mathbf{0}_{N_i} & \mathbf{O}_{N_i} & \mathbf{O}_{N_i} \\
        0 & \mathbf{0}_{N_i}^T & \mathbf{0}_{N_i}^T & -1 & \mathbf{1}_{N_i}^T & \mathbf{1}_{N_i}^T &...& -1 & \mathbf{1}_{N_i}^T & \mathbf{1}_{N_i}^T \\
        \vdots & \vdots  & \vdots  & \vdots   & \vdots  & \vdots & &\vdots   & \vdots  & \vdots  \\
        0 & \mathbf{0}_{N_i}^T & \mathbf{0}_{N_i}^T & -1 & \mathbf{1}_{N_i}^T & \mathbf{1}_{N_i}^T& ...&-1 & \mathbf{1}_{N_i}^T & \mathbf{1}_{N_i}^T 
    \end{bmatrix} \otimes \mathbf{I}_n\in\mathbb{R}^{(N_i+M)n\times(1+M)(1+2N_i)n}
\end{align*}}
 \textcolor{black}{This matrix is designed to satisfy the condition, for all $i \in \mathcal{N}$:}
 \textcolor{black}{\begin{align*}
   \mathbf{B}_i \cdot \mathbf{H}_i^T\mathbf{H}_i &= \mathbf{H}_i, \quad\\
    \mathbf{B}_i\cdot \mathbf{A}_i^T\mathbf{A}_i &= \mathbf{O}_{(N_i+M)n\times(1+2N_i)n}. \tag{G.4}\label{a84}
\end{align*}}
Moreover, substituting equation (\ref{27}) into the update rule of $\tilde{\mathbf{z}}_i^{t+1}$ yields:
\[
\mathbf{H}_i^T \mathbf{H}_i \mathbf{z}_i = -\mathbf{U}_i (\tilde{\mathbf{z}}_i^{t+1} - \mathbf{z}_i^t) + \mathbf{H}_i^T \mathbf{D}_i \mathbf{w}_i^t - c \mathbf{A}_i^T \mathbf{A}_i \mathbf{z}_i^t - c \mathbf{A}_i^T \boldsymbol{\lambda}_i, \quad \forall i \in \mathcal{N}. \tag{G.5}\label{a85}
\]
Multiplying both sides of equation (\ref{a85}) by the matrix$\mathbf{B}_i$
and then taking the norm, we obtain:
\begin{align*}
\|\mathbf{H}_i \mathbf{z}_i^{t+1}\| &= \| \mathbf{B}_i\cdot \mathbf{U}_i (\tilde{\mathbf{z}}_i^{t+1} - \mathbf{z}_i^t) + \mathbf{D}_i \mathbf{w}_i^t\| \\
&\leq 2 \|\tilde{\mathbf{z}}_i^{t+1} - \mathbf{z}_i^t\| + d_{\max} \|\mathbf{w}_i^t\|, 
\end{align*}
where $\textcolor{black}{d_{\max} = \max\left(\{d_{i,j} \mid i \in \mathcal{N}, j \in \mathcal{N}_i \} \cup \{r_{i, \textcolor{black}{k}} \mid i \in \mathcal{N},k\in\mathcal{M}_i\}\right)}$. The first inequality follows from equation (\ref{a84}), and the second is based on the triangle inequality. Since the right-hand side of the inequality is bounded, it follows that $\{\|\mathbf{H}_i \mathbf{z}_i^t\|\}_{t\geq1}$ is bounded for all $i \in \mathcal{N}$. Recalling the definitions $\mathbf{z}_i $, along with $\mathbf{H}_i$ and $\mathbf{A}_j$, it follows that:
\begin{align*}
\mathbf{H}_i \mathbf{z}_i^t &= \text{vec}\left(\mathbf{x}_i^t - (\mathbf{p}_{i,j}^+)^t \mid j \in \mathcal{N}_j \quad,\quad \mathbf{x}_i^t - \mathbf{y}_{i, \textcolor{black}{k}}^t \mid \textcolor{black}{k\in\mathcal{M}_i}\right),\\
\mathbf{A}_{j} \mathbf{z}_{j}^t &= \text{vec} \left( \mathbf{x}_{j}^{t} - (\mathbf{p}_{j,i}^{-})^{t}\mid i\in \mathcal{N}_{j} \quad,\quad\mathbf{y}_{j}^{t} - (\mathbf{q}_{j, i, \textcolor{black}{k}}^{-})^{t}\mid i \in \mathcal{N}_{j}, k\in\mathcal{M} \quad,\quad\mathbf{x}_{j}^{t} - (\mathbf{q}_{j, i, \textcolor{black}{k}}^{+})^{t}\mid i \in \mathcal{N}_{j}, k\in\mathcal{M} \right). \tag{G.6}\label{a86}
\end{align*}
Using the triangle inequality, for all $t \geq 0$, we obtain:
\[
\| (\mathbf{p}_{i,j}^{+})^{t} \| \leq \| \mathbf{x}_{i}^{t} \| + \| \mathbf{x}_{i}^{t} - (\mathbf{p}_{i,j}^{+})^{t} \|, \quad \| \mathbf{x}_{j}^{t} \| \leq \| (\mathbf{p}_{j,i}^{-})^{t} \| + \| \mathbf{x}_{j}^{t} - (\mathbf{p}_{j,i}^{-})^{t} \|.\tag{G.7}\label{a87}
\]
From the update rule of $\mathbf{z}_i^t$, we also have:
\[
(\mathbf{p}_{i,j}^{+})^{t} = (\mathbf{p}_{j,i}^{-})^{t}, \, \forall j \in \mathcal{N}_{i}, \, i \in \mathcal{N}_{j}, \quad \text{and} \quad \mathbf{p}_{i}^{t} = \mathbf{a}_i, \, \forall i \in \mathcal{A} \subseteq \mathcal{N}. \tag{G.8}\label{a88}
\]
Since at least one anchor exists, combining equation (\ref{a87}) and equation (\ref{a88}) yields:
\[
\| \mathbf{x}_{j}^{t} \| \leq \| \mathbf{a}_i \| + \| \mathbf{x}_{i}^{t} - (\mathbf{p}_{i,j}^{+})^{t} \| + \| \mathbf{x}_{j}^{t} - (\mathbf{p}_{j,i}^{-})^{t} \|, \quad \forall j \in \mathcal{N}_{i}, \, i \in \mathcal{A}. \tag{G.9}\label{a89}
\]
Because both $\| \mathbf{H}_i^t \mathbf{z}_i \|$ and $\| \mathbf{A}_j^t \mathbf{z}_j \|$ are bounded for all \( i \in \mathcal{N} \), it follows from equation (\ref{a86}) that the right-hand side of inequality (\ref{a89}) is bounded. Hence, the sequence $\{ \| \mathbf{p}_j^t \| \}_{t \geq 1}$ is bounded for all $j \in \mathcal{N}_i$, $i \in \mathcal{A}$.

By a similar argument as in equations (\ref{a87})–(\ref{a89}), we can also prove that for the neighbor sensor nodes $j \in \mathcal{N}_i \cap \mathcal{A}$, the sequence $\{ \mathbf{x}_i^t \}_{t \geq 1}$ is bounded. Since the graph is connected, we conclude that the sequence $\{ \mathbf{x}_i^t \}_{t \geq 1}$ is bounded for all $i \in \mathcal{N}$.

Similarly, using the triangle inequality:
\begin{align*}
   \| \mathbf{y}_{i, \textcolor{black}{k}}^{t} \| \leq \| \mathbf{x}_{i}^{t} \| + \| \mathbf{x}_{i}^{t} - \mathbf{y}_{i, \textcolor{black}{k}}^{t}  \|, \forall k\in\mathcal{M},  i\in\mathcal{N},
\end{align*}
and the fact that $\| \mathbf{H}_i^t \mathbf{z}_i \|$ is bounded for all \( i \in \mathcal{N} \), we conclude that the sequence $\{ \mathbf{y}_{i, \textcolor{black}{k}}^t \}_{t \geq 1}$ is also bounded.

Furthermore, based on the definition of $\mathbf{z}_j^t$ and the triangle inequality:
\begin{align*}
    \| \mathbf{z}_{i}^{t}\|
    =&\|\mathbf{x}_i^t\|+\|\mathbf{y}_{i, \textcolor{black}{k}}^t\|+\sum_{j\in\mathcal{N}_i}(\|(\mathbf{p}_{i,j}^-)^t\|+\|(\mathbf{p}_{i,j}^+)^t\|+\|(\mathbf{q}_{i, j, \textcolor{black}{k}}^-)^t\|+\|(\mathbf{q}_{i, j, \textcolor{black}{k}}^+)^t\|) \\
    \leq& 2\|\mathbf{x}_i^t\|+\sum_{k\in\mathcal{M}_i}\|\mathbf{x}_i^t-\mathbf{y}_{i, \textcolor{black}{k}}^t\|+\sum_{i \in \mathcal{N}_{j}} \Big(2\| \mathbf{x}_{i}^{t} \| + \| \mathbf{x}_{i}^{t} - (\mathbf{p}_{i,j}^{-})^{t} \| + \| \mathbf{x}_{i}^{t} - (\mathbf{p}_{i,j}^{+})^{t} \|\\
    &+\sum_{k\in\mathcal{M}}(2\| \mathbf{y}_{i, \textcolor{black}{k}}^{t} \| + \| \mathbf{y}_{i, \textcolor{black}{k}}^{t} - (\mathbf{q}_{i, j, \textcolor{black}{k}}^{-})^{t} \| + \| \mathbf{y}_{i, \textcolor{black}{k}}^{t} - (\mathbf{q}_{i, j, \textcolor{black}{k}}^{+})^{t} \|)\Big), \quad \forall i \in \mathcal{N}. \tag{G.10}\label{a810}
\end{align*}
Since the right-hand side of equation (\ref{a810}) is bounded, we conclude that the sequence $\{ \| \mathbf{z}_i^t \| \}_{t \geq 1}$ is bounded for all \( i \in \mathcal{N} \).

3. Rearranging the terms in equation (\ref{a85}) and taking the norm yields:
\[
\|\boldsymbol{\lambda}_i^t\| \leq \|\mathbf{A}_i^T \boldsymbol{\lambda}_i^t\| \leq \|\mathbf{H}_i^T \mathbf{H}_i \mathbf{z}_i^t\| + \|U_i (\tilde{\mathbf{z}}_i^{t+1} - \mathbf{z}_i^t)\| + \|\mathbf{H}_i^T \mathbf{D}_i \mathbf{w}_i^t\| + \|c \mathbf{A}_i^T \mathbf{A}_i \mathbf{z}_i^t\|. \tag{G.11}\label{a811}
\]
The first inequality above holds because the smallest eigenvalue of $\mathbf{A}_i \mathbf{A}_i^T$ is 1. The second inequality follows from the triangle inequality. Since the right-hand side of the second inequality in equation (\ref{a811}) is bounded, it follows that the sequence $\{\boldsymbol{\lambda}_i^t\}_{t \geq 1}$ is bounded for all $i \in \mathcal{N}$.

4. By the triangle inequality, we have
\[
\|\tilde{\mathbf{z}}_i^{t+1}\| \leq \|\tilde{\mathbf{z}}_i^t - \mathbf{z}_i^t\| + \|\mathbf{z}_i^t\|.
\]
Since the right-hand side is bounded, the sequence $\{\tilde{\mathbf{z}}_i^t\}_{t \geq 1}$ is also bounded for all $i \in \mathcal{N}$.  
$\hfill \Box$

\textbf{Lemma A.8.3 (Bounded Subgradient)}  
Assume the sequence $\{\mathbf{s}^t = (\mathbf{z}^t, \mathbf{w}^t, \boldsymbol{\lambda}^t)\}_{t \geq 0}$ is generated by Algorithm \ref{al1}, where $\kappa_1$, $\kappa_2$, $c$, and $\rho$ are the parameters in the potential function $\mathcal{P}$. Then the following holds:
\[
\boldsymbol{\xi}^{t+1} \in \partial \mathcal{P}^{t+1}, \quad \forall t \geq 0,
\]
where $
\boldsymbol{\xi}^{t+1} := \big(\boldsymbol{\xi}_{\mathbf{z}^{t+1}}, \boldsymbol{\xi}_{\mathbf{w}^{t+1}}, \boldsymbol{\xi}_{\boldsymbol{\lambda}^{t+1}}, \boldsymbol{\xi}_{\tilde{\mathbf{z}}^{t+1}}, \boldsymbol{\xi}_{\mathbf{w}^{t}}, \boldsymbol{\xi}_{\mathbf{z}^{t}}\big),
$
and
\begin{align*}\label{a812}
\boldsymbol{\xi}_{\mathbf{z}_i^{t+1}} &:= \mathbf{U}_i \left( \mathbf{z}_i^t - \tilde{\mathbf{z}}_i^{t+1} \right) 
+ \mathbf{Q}_i^T \mathbf{D}_i \left( \mathbf{w}_i^t - \mathbf{w}_i^{t+1} \right) 
+ (\kappa_2 + 1) \mathbf{A}_i^T \left( \boldsymbol{\lambda}_i^{t+1} - \boldsymbol{\lambda}_i^t \right) \\
&\quad + \left( c \mathbf{A}_i^T \mathbf{A}_i + \mathbf{Q}_i^T \mathbf{Q}_i + c (\kappa_1 + \kappa_2) \mathbf{B}_i^T \mathbf{B}_i \right) \left( \mathbf{z}_i^{t+1} - \mathbf{z}_i^t \right), \\
\boldsymbol{\xi}_{\mathbf{w}_i^{t+1}} &:= \frac{\rho}{2} \left( \mathbf{w}_i^t - \mathbf{w}_i^{t+1} \right), \quad 
\boldsymbol{\xi}_{\boldsymbol{\lambda}_i^{t+1}} := \frac{1}{c} \left( \boldsymbol{\lambda}_i^{t+1} - \boldsymbol{\lambda}_i^t \right), \\
\boldsymbol{\xi}_{\tilde{\mathbf{z}}_i^{t+1}} &:= \kappa_1 \mathbf{A}_i^T \left( c \mathbf{A}_i \left( \tilde{\mathbf{z}}_i^{t+1} - \mathbf{z}_i^t \right) + \boldsymbol{\lambda}_i^{t+1} - \boldsymbol{\lambda}_i^t \right), \\
\boldsymbol{\xi}_{\mathbf{w}_i^t} &:= \frac{\rho}{2} \left( \mathbf{w}_i^t - \mathbf{w}_i^{t+1} \right), \quad 
\boldsymbol{\xi}_{\mathbf{z}_i^t} := c \left( \kappa_1 + \kappa_2 \right) \mathbf{B}_i^T \mathbf{B}_i \left( \mathbf{z}_i^t - \mathbf{z}_i^{t+1} \right), \quad i \in \mathcal{N}.\tag{G.12}
\end{align*}
Furthermore, for any $t \geq 0$, the following inequality holds:
\begin{align*}
\|\boldsymbol{\xi}^{t+1}\| \leq &\sum_{i \in \mathcal{N}} \alpha_i \|\mathbf{z}_i^{t+1} - \mathbf{z}_i^t\| + \sum_{i \in \mathcal{N}} \left(\rho + d_{\max} \sqrt{N_i + M_i+1}\right) \|\mathbf{w}_i^{t+1} - \mathbf{w}_i^{t}\|\\
&+ \sum_{i \in \mathcal{N}} \beta_i \|\boldsymbol{\lambda}_i^{t+1} - \boldsymbol{\lambda}_i^t\|\\
\leq& C_2 \|\mathbf{s}^{t+1} - \mathbf{s}^t\|, \tag{G.13}\label{a813}
\end{align*}
where
\begin{align*}
\alpha_i :=& c(2N_i+1)+N_i + M_i+1 + 2 (\kappa_1 + \kappa_2) (1+c) (2 + N_{\max}) \\
&+ \frac{2[(c+1) N_i+1] (c + 1) \sqrt{3n N_{\text{sum}} (2N_i + 1)}}{\sqrt{\tilde{\tau}_{\min}}},\\
\beta_i :=& \sqrt{2N_i + 1} (\kappa_1 + 1 + \kappa_2) + \frac{1}{c},\\
C_2 :=& \sqrt{3N} \cdot \max \{\alpha_i, \rho + d_{\max} \sqrt{N_i + M_i+1}, \beta_i, i \in \mathcal{N}\}. \tag{G.14}\label{a814}
\end{align*}

\textbf{Proof:}  
Taking partial derivatives of $\mathcal{P}^{t+1}$ with respect to $\mathbf{z}_i^{t+1}, \mathbf{w}_i^{t+1}, \boldsymbol{\lambda}_i^{t+1},\tilde{\mathbf{z}}_i^{t+1}, \mathbf{w}_i^t, \mathbf{z}_i^t$, and combining with the update step of $\boldsymbol{\lambda}_i^{t+1}$, we obtain
\begin{align*}\label{a815}
\nabla_{\mathbf{z}_i^{t+1}} \mathcal{P}^{t+1} &= 
\nabla_{\mathbf{z}_i} \mathcal{L}_i \left( \mathbf{z}_i^{t+1}, \mathbf{w}_i^{t+1}, \boldsymbol{\lambda}_i^{t+1} \right)
+ \kappa_2 \mathbf{A}_i^T \left( \boldsymbol{\lambda}_i^{t+1} - \boldsymbol{\lambda}_i^t \right) \\
&\quad + c (\kappa_1 + \kappa_2) \mathbf{B}_i^T \mathbf{B}_i \left( \mathbf{z}_i^{t+1} - \mathbf{z}_i^t \right),\\
\partial_{\mathbf{w}_i^{t+1}} \mathcal{P}^{t+1} &= 
\partial_{\mathbf{w}_i} \mathcal{L}_i \left( \mathbf{z}_i^{t+1}, \mathbf{w}_i^{t+1}, \boldsymbol{\lambda}_i^{t+1} \right) 
+ \frac{\rho}{2} \left( \mathbf{w}_i^{t+1} - \mathbf{w}_i^t \right), \\
\nabla_{\boldsymbol{\lambda}_i^{t+1}} \mathcal{P}^{t+1} &= 
\nabla_{\boldsymbol{\lambda}_i} \mathcal{L}_i \left( \mathbf{z}_i^{t+1}, \mathbf{w}_i^{t+1}, \boldsymbol{\lambda}_i^{t+1} \right), \\
\nabla_{\tilde{\mathbf{z}}_i^{t+1}} \mathcal{P}^{t+1} &= \boldsymbol{\xi}_{\tilde{\mathbf{z}}_i^{t+1}}, \quad 
\nabla_{\mathbf{w}_i^t} \mathcal{P}^{t+1} = \boldsymbol{\xi}_{\mathbf{w}_i^t}, \quad 
\nabla_{\mathbf{z}_i^t} \mathcal{P}^{t+1} = \boldsymbol{\xi}_{\mathbf{z}_i^t}.
\tag{G.15}
\end{align*}
where the definitions of $\boldsymbol{\xi}_{\tilde{\mathbf{z}}_i^{t+1}}, \boldsymbol{\xi}_{\mathbf{w}_i^{t}}, \boldsymbol{\xi}_{\mathbf{z}_i^t}$ can be found in equation (\ref{a812}). Applying the update rule of $\tilde{\mathbf{z}}_i^{t+1}$ yields
\begin{align*}
\nabla_{\mathbf{z}_i} \mathcal{L}_i (\mathbf{z}_i^{t+1}, \mathbf{w}_i^{t+1}, \boldsymbol{\lambda}_i^{t+1}) =& \mathbf{U}_i (\mathbf{z}_i^t - \tilde{\mathbf{z}}_i^{t+1}) + \mathbf{H}_i^T \mathbf{D}_i (\mathbf{w}_i^t - \mathbf{w}_i^{t+1}) + \mathbf{A}_i^T (\boldsymbol{\lambda}_i^{t+1} - \boldsymbol{\lambda}_i^t)\\
&+ \left(c \mathbf{A}_i^T \mathbf{A}_i + \mathbf{H}_i^T \mathbf{H}_i\right) (\mathbf{z}_i^{t+1} - \mathbf{z}_i^t). \tag{G.16}\label{a816}
\end{align*}

Combining equations (\ref{a816}) and (\ref{a815}) gives $\boldsymbol{\xi}_{\mathbf{z}_i^{t+1}} = \nabla_{\mathbf{z}_i^{t+1}} \mathcal{P}^{t+1}$. According to the optimality condition of problem (\ref{gengxin1}b) with respect to variable $\mathbf{v}$, we have
\begin{align*}
\mathbf{0} \in& \partial_{\mathbf{w}_i} \mathcal{L}_i (\mathbf{z}_i^{t+1}, \mathbf{w}_i^{t+1}, \boldsymbol{\lambda}_i^{t}) + \rho (\mathbf{w}_i^{t+1} - \mathbf{w}_i^t)\\
=&  \partial_{\mathbf{w}_i}\mathcal{L}_i (\mathbf{z}_i^{t+1}, \mathbf{w}_i^{t+1}, \boldsymbol{\lambda}_i^{t+1}) + \rho (\mathbf{w}_i^{t+1} - \mathbf{w}_i^t)\tag{G.17}\label{a817}
\end{align*}
Combining equations (\ref{a817}) and (\ref{a815}) shows that $\boldsymbol{\xi}_{\mathbf{z}_i^{t+1}} \in \partial_{\mathbf{z}_i^{t+1}} \mathcal{P}^{t+1}$. Using the update rule of $\boldsymbol{\lambda}_i^{t+1}$, we have
\[
\mathbf{w}_{\boldsymbol{\lambda}_i^{t+1}} = \frac{1}{c} (\boldsymbol{\lambda}_i^{t+1} - \boldsymbol{\lambda}_i^t) = \mathbf{A}_i \mathbf{z}_i^{t+1} = \nabla_{\boldsymbol{\lambda}_i} \mathcal{L}_i (\mathbf{z}_i^{t+1}, \mathbf{w}_i^{t+1}, \boldsymbol{\lambda}_i^{t+1}) = \nabla_{\boldsymbol{\lambda}_i^{t+1}} \mathcal{P}^{t+1}.
\]
Applying the expression of $\boldsymbol{\xi}^{t+1}$ and using the triangle inequality, we derive
\begin{align*}
\|\boldsymbol{\xi}^{t+1}\| \leq& \|\mathbf{U}_i \|\|(\tilde{\mathbf{z}}_i^{t+1} - \mathbf{z}_i^t)\| + (2c(\kappa_1 + \kappa_2) \|\mathbf{B}_i^T \mathbf{B}_i\| + c \|\mathbf{A}_i^T \mathbf{A}_i\| + \|\mathbf{H}_i^T \mathbf{H}_i\|) \|\mathbf{z}_i^{t+1} - \mathbf{z}_i^t\|\\
&+ (\sigma_{\max}^{1/2} (\mathbf{H}_i \mathbf{H}_i^T) (\|\mathbf{D}_i\| + \rho)) \|\mathbf{w}_i^{t+1} - \mathbf{w}_i^t\|\\
&+( (\sigma_{\max}^{1/2} (\mathbf{A}_i \mathbf{A}_i^T) (\kappa_2+ 1 +\kappa_1) + \frac{1}{c}) \|\boldsymbol{\lambda}_i^{t+1} - \boldsymbol{\lambda}_i^t\|. \tag{G.18}\label{a818}
\end{align*}
where $\sigma_{\max} (\mathbf{H}_i \mathbf{H}_i^T)$ and $\sigma_{\max} (\mathbf{A}_i \mathbf{A}_i^T)$ denote the largest eigenvalues of $\mathbf{H}_i \mathbf{H}_i^T$ and $\mathbf{A}_i \mathbf{A}_i^T$, respectively. Taking square roots of inequality (\ref{a323}) and using the fact that $\|\mathbf{x}\|_2 \leq \|\mathbf{x}\|_1 \leq \sqrt{N} \|\mathbf{x}\|_2$ for any $\mathbf{x} \in \mathbb{R}^N$, we get
\[
\sum_{i \in \mathcal{N}} \|\tilde{\mathbf{z}}_i^{t+1} - \mathbf{z}_i^t\| \leq \frac{(c+1) \sqrt{3n N_{\text{sum}}}}{\sqrt{\tilde{\tau}_{\min}}} \sum_{i \in \mathcal{N}} \|\mathbf{A}_i (\mathbf{z}_i^{t+1} - \mathbf{z}_i^t)\|. \tag{G.19}\label{a819}
\]

Moreover, from equations (\ref{27})  we obtain
\begin{align*}
&\|\mathbf{U}_i\| = \max\big\{2 [(c + 1) N_i+M_i],2(2cN_i+1)\big\}\leq2 [(2c + 1) N_i+M_i+1], \\ &\sigma_{\max} (\mathbf{A}_i \mathbf{A}_i^T) = 2N_i+1,\quad\sigma_{\max} (\mathbf{H}_i \mathbf{H}_i^T) = N_i + M_i+1. \tag{G.20}\label{a820}
\end{align*}
Substituting equations (\ref{a89}), (\ref{a819}), and (\ref{a820}) into equation (\ref{a818}), we obtain
\begin{align*}
\|\boldsymbol{\xi}^{t+1}\| \leq& \sum_{i \in \mathcal{N}} \alpha_i \|\mathbf{z}_i^{t+1} - \mathbf{z}_i^t\| + \sum_{i \in \mathcal{N}} \left( \rho + d_{\max} \sqrt{N_i + M_i+1} \right) \|\mathbf{w}_i^{t+1} - \mathbf{w}_i^t\| \\
&+ \sum_{i \in \mathcal{N}} \beta_i \|\boldsymbol{\lambda}_i^{t+1} - \boldsymbol{\lambda}_i^t\|, \tag{G.21}
\end{align*}
where the constants $\alpha_i$ and $\beta_i$ are defined in equation (\ref{a814}). Using the fact that $\|x\|_1 \leq \sqrt{3N} \|x\|_2$ for all $x \in \mathbb{R}^{3N}$, we conclude
\[
\|\boldsymbol{\xi}^{t+1}\| \leq C_2 \|\mathbf{s}^{t+1} - \mathbf{s}^t\|,
\]
where the constant $C_2$ is defined in equation (\ref{a814}). This completes the proof.  
$\hfill \Box$

In the following, let $\sigma(\{x^t\}_{t \geq 1})$ denote the set of limit points of the sequence $\{x^t\}_{t \geq 1}$, and define$\mathrm{crit}\, \mathcal{L} := \{\mathbf{s} : 0 \in \partial \mathcal{L}(\mathbf{s})\}$ as the set of critical points of $\mathcal{L}$.\\

\textbf{Lemma A.8.4 (Properties of the Limit Set).}  
Let $\{\mathbf{z}^t, \mathbf{w}^t, \boldsymbol{\lambda}^t\}_{t \geq 1}$ be the sequence generated by Algorithm \ref{al1}. Suppose the parameters satisfy conditions (\ref{45})–(\ref{47}), and define $\Omega := \sigma(\{(\mathbf{z}^t, \mathbf{w}^t, \boldsymbol{\lambda}^t, \tilde{\mathbf{z}}^t, \mathbf{w}^{t-1}, \mathbf{z}^{t-1})\}_{t \geq 1})$. Then the following statements hold:

1. The set $\Omega$ is a nonempty compact set;\\
2. $\lim_{t \to \infty} \mathrm{dist}\left[(\mathbf{z}^t, \mathbf{w}^t, \boldsymbol{\lambda}^t, \tilde{\mathbf{z}}^t, \mathbf{w}^{t-1}, \mathbf{z}^{t-1}), \Omega\right] = 0$;\\
3. $\Omega \subseteq \{(\mathbf{z}, \mathbf{w}, \lambda, \mathbf{z}, \mathbf{w}, \mathbf{z}) : (\mathbf{z}, \mathbf{w}, \lambda) \in \mathrm{crit}\, \mathcal{L} \}$;\\
4. Any critical point of $\mathcal{L}$ is a KKT point of problem (\ref{jin});\\
5. The potential function $\mathcal{P}^t$ is constant and bounded over the set $\Omega$.
\begin{proof}
We prove the above statements one by one:

1. From Lemma A.8.2, the sequence $\{(\mathbf{z}^t, \mathbf{w}^t, \boldsymbol{\lambda}^t, \tilde{\mathbf{z}}^t)\}_{t \geq 1}$ is bounded. Therefore, $\Omega$ is nonempty and bounded. Since $\Omega$ is defined as a limit set, it is closed by definition. Hence, $\Omega$ is compact.

2. This follows directly from the definition of the limit set.

3. Since the sequence $\{(\mathbf{z}^t, \mathbf{w}^t, \boldsymbol{\lambda}^t)\}_{t \geq 1}$ is bounded, there exists a convergent subsequence $\{(\mathbf{z}^{t_k}, \mathbf{w}^{t_k}, \boldsymbol{\lambda}^{t_k})\}_{k\geq 1}$ that converges to a limit point $(\mathbf{z}^*, \mathbf{w}^*, \boldsymbol{\lambda}^*)$. As $t \to \infty$, combining Lemma A.8.1 and equation (\ref{a28}), we obtain
\[
\sum_{i \in \mathcal{N}} \|\tilde{\mathbf{z}}_i^{t+1} - \mathbf{z}_i^t\|^2 \to 0, \quad \sum_{i \in \mathcal{N}} \|\mathbf{w}_i^{t+1} - \mathbf{w}_i^t\|^2 \to 0, \quad \sum_{i \in \mathcal{N}} \|\mathbf{z}_i^{t+1} - \mathbf{z}_i^t\|^2 \to 0.
\]
which implies that $\{(\tilde{\mathbf{z}}^{t_k}, \mathbf{w}^{t_k-1}, \mathbf{z}^{t_k-1})\}_{k \geq 1}$ converges to $(\mathbf{z}^*, \mathbf{w}^*, \mathbf{z}^*) \in \Omega$. Therefore,
\begin{align*}
    (\mathbf{z}^*, \mathbf{w}^*, \boldsymbol{\lambda}^*,\mathbf{z}^*, \mathbf{w}^*,\mathbf{z}^*) \in \Omega
\end{align*}
By continuity of the function $\mathcal{P}^t$, we obtain
\[
\lim_{k \to \infty} \mathcal{P}^{t_k} = \mathcal{P}^*.
\]
On the other hand, by Lemma A.8.1, equation (\ref{a28}), and Lemma A.8.3, we know that $\mathbf{w}^{t_k} \in \partial \mathcal{P}^{t_k}$ and $\mathbf{w}^{t_k} \to 0$ as $k \to \infty$. Since the subdifferential $\partial \mathcal{P}$ is closed see ([\cite{bolte2014proximal}, Remark 1(ii)), we have
\[
\mathbf{0} \in \partial \mathcal{P}^*. \tag{G.22}\label{a822}
\]
Furthermore, from equation (\ref{a815}), it follows that
\[
\partial \mathcal{P}^* = \partial \mathcal{L}(\mathbf{z}^*, \mathbf{w}^*, \boldsymbol{\lambda}^*). \tag{G.23}\label{a823}
\]
Combining equations (\ref{a822}) and (\ref{a823}) shows that $(\mathbf{z}^*, \mathbf{w}^*, \boldsymbol{\lambda}^*)$ is a critical point of $\mathcal{L}$.

4. Suppose $(\mathbf{z}^*, \mathbf{w}^*, \boldsymbol{\lambda}^*) \in \mathrm{crit} \mathcal{L}$. Since $\mathcal{L}$ is separable across nodes $i$, we have
\begin{align*}
\mathbf{0} &= \nabla_{\mathbf{z}_i} \mathcal{L}_i(\mathbf{z}_i^*, \mathbf{w}_i^*, \boldsymbol{\lambda}_i^*) = \nabla_{\mathbf{z}_i} G_i(\mathbf{z}_i^*,  \mathbf{w}_i^*) + \mathbf{A}_i^T \boldsymbol{\lambda}_i^* + c \mathbf{A}_i^T\mathbf{A}_i \mathbf{z}_i^*,\tag{G.24}\label{a824}\\
\mathbf{0} & \in \partial_{\mathbf{w}_i} \mathcal{L}_i(\mathbf{z}_i^*, \mathbf{w}_i^*, \boldsymbol{\lambda}_i^*) = \nabla_{\mathbf{w}_i} G_i(\mathbf{z}_i^*, \mathbf{w}_i^*) + \delta_{\mathcal{B}^{\textcolor{black}{N_i+M}}}(\mathbf{w}_i^*),\tag{G.25}\label{a825}\\
\mathbf{0} &= \nabla_{\boldsymbol{\lambda}_i} \mathcal{L}_i(\mathbf{z}_i^*, \mathbf{w}_i^*, \boldsymbol{\lambda}_i^*) = \mathbf{A}_i \mathbf{z}_i^*. \tag{G.26}\label{a826}
\end{align*}
Combining equations (\ref{a824}) and (\ref{a826}), we obtain
\[
\nabla_{\mathbf{z}_i} G_i(\mathbf{z}_i^*, \mathbf{w}_i^*) + \mathbf{A}_i^T \boldsymbol{\lambda}_i^* = \mathbf{0}. \tag{G.27}\label{a827}
\]
Since for any fixed $\mathbf{w}_i$, the function $G_i(\mathbf{z}_i, \mathbf{w}_i) + \langle \boldsymbol{\lambda}_i, \mathbf{A}_i \mathbf{z}_i \rangle$ is convex in $\mathbf{z}_i$ over the set $\mathcal{X} \cap \mathcal{Y}$, equation (\ref{a827}) implies that $\mathbf{z}_i^*$ also satisfies
\[
\mathbf{z}_i^* = \arg\min_{\mathbf{z}_i \in \mathcal{X} \cap \mathcal{Y}} \left[ G_i(\mathbf{z}_i, \mathbf{w}_i^*) + \langle \boldsymbol{\lambda}_i^*, \mathbf{A}_i \mathbf{z}_i \rangle \right]. \tag{G.28}\label{a828}
\]
Combining (\ref{a825}), (\ref{a826}), and (\ref{a828}), we conclude that $(\mathbf{z}^*, \mathbf{w}^*, \boldsymbol{\lambda}^*)$ is a KKT point of problem (\ref{jin}).

5. From Lemma A.8.1, the sequence $\{\mathcal{P}^t\}_{t \geq 1}$ is monotonically decreasing. In addition, equation (\ref{xiajie}) ensures that $\mathcal{P}^t$ has a lower bound. Therefore, the sequence $\{\mathcal{P}^t\}_{t \geq 1}$ converges to a finite limit $\mathcal{P}^*$. It follows that $\mathcal{P}^t$ is constant on the set $\Omega$.

To establish the global convergence of the entire sequence, we require the Kurdyka–Łojasiewicz (KŁ) property, which has been extensively used in the literature \cite{bolte2014proximal, guo2017convergence, yashtini2022convergence, bolte2018first}. The proof of Theorem \ref{the2} is given below.

As in Lemma A.8.4, let $\Omega := \sigma(\{(\mathbf{z}^t, \mathbf{w}^t, \boldsymbol{\lambda}^t, \tilde{\mathbf{z}}^t, \mathbf{w}^{t-1}, \mathbf{z}^{t-1})\}_{t \geq 1})$. The set $\Omega$ is nonempty and compact. From Lemma A.8.4 and the continuity of $\mathcal{P}^t$, we have
\[
\lim_{t \to \infty} \mathcal{P}^t = \mathcal{P}^*, \quad \forall(\mathbf{z}^*, \mathbf{w}^*, \boldsymbol{\lambda}^*, \mathbf{z}^*, \mathbf{w}^*, \mathbf{z}^*) \in \Omega. \tag{G.29}\label{a829}
\]
Moreover, by Lemma A.8.1, we obtain
\[
C_3 \|\mathbf{s}^{t+1} - \mathbf{s}^t\|^2 \leq \mathcal{P}^t - \mathcal{P}^{t+1}, \quad \forall t \geq 1, \tag{G.30}\label{a830}
\]
where $C_3 := \min\{c, 1, \frac{\rho}{4}, C_0\}$. Since $\mathcal{P}^t$ is monotonically non-increasing, it holds that $\mathcal{P}^t \geq \mathcal{P}^*$ for all $t$. We consider two possible cases:

1. There exists an integer $\overline{t} \geq 0$ such that $\mathcal{P}^{\overline{t}} = \mathcal{P}^*$. Then from Lemma A.8.1 and equation (\ref{a830}), we have
\[
C_3 \|\mathbf{s}^{t+1} - \mathbf{s}^t\|^2 \leq \mathcal{P}^t - \mathcal{P}^{t+1} = \mathcal{P}^{\overline{t}} - \mathcal{P}^* = 0, \quad \forall t \geq \overline{t}.
\]
Therefore, for all $t \geq \overline{t}$, we have $\mathbf{s}^{t+1} = \mathbf{s}^t$, and consequently $\sum_{t=0}^\infty \|\mathbf{s}^{t+1} - \mathbf{s}^t\|^2 < \infty$.

2. For all $t$, we have $\mathcal{P}^t > \mathcal{P}^*$. From equation (\ref{a829}), for any $\eta > 0$, there exists $t_0 > 0$ such that
\[
\mathcal{P}^t < \mathcal{P}^* + \eta, \quad \forall t \geq t_0.
\]
From Lemma A.8.4, we know that $\lim_{t \to \infty} \mathrm{dist}((\mathbf{z}^t, \mathbf{w}^t, \boldsymbol{\lambda}^t, \tilde{\mathbf{z}}^{t}, \mathbf{w}^{t-1}, \mathbf{z}^{t-1}), \Omega) = 0$. This implies that for any $\gamma > 0$, there exists $t_1 \geq 1$ such that for all $t \geq t_1$,
\[
\mathrm{dist}[(\mathbf{z}^t, \mathbf{w}^t, \boldsymbol{\lambda}^t, \tilde{\mathbf{z}}^{t}, \mathbf{w}^{t-1}, \mathbf{z}^{t-1}), \Omega] < \gamma.
\]
Thus, for all $t \geq t_2 := \max\{t_0, t_1\}$ and for any $\eta, \gamma > 0$, we have
\[
\mathrm{dist}[(\mathbf{z}^t, \mathbf{w}^t, \boldsymbol{\lambda}^t, \tilde{\mathbf{z}}^{t}, \mathbf{w}^{t-1}, \mathbf{z}^{t-1}), \Omega] < \gamma, \quad \text{and} \quad \mathcal{P}^* < \mathcal{P}^t < \mathcal{P}^* + \eta.
\]

Since $\mathcal{P}$ is constant over $\Omega$ and semi-algebraic (i.e., a polynomial), it satisfies the Kurdyka–Lojasiewicz (KL) property (see \cite{bolte2018first}, Theorem 6.1). Applying Lemma \ref{lem12} to $\Omega$, we obtain:
\[
\varphi'(\mathcal{P}^t - \mathcal{P}^*) \cdot \mathrm{dist}(\mathbf{0}, \partial \mathcal{P}^t) \geq 1, \quad \forall t \geq t_2. \tag{G.31}\label{a831}
\]
Due to the concavity of $\varphi$, it follows that
\[
\varphi(\mathcal{P}^t - \mathcal{P}^*) - \varphi(\mathcal{P}^{t+1} - \mathcal{P}^*) \geq \varphi'(\mathcal{P}^t - \mathcal{P}^*) (\mathcal{P}^t - \mathcal{P}^{t+1}). \tag{G.32}\label{a832}
\]
From Lemma A.8.3 (equation (\ref{a813})), we have
\[
\mathrm{dist}(\mathbf{0}, \partial \mathcal{P}^t) \leq C_2 \|\mathbf{s}^t - \mathbf{s}^{t+1}\|, \quad C_2 > 0.
\]
Combining this with equations (\ref{a831}) and (\ref{a832}), and the fact that $\varphi'(\mathcal{P}^t - \mathcal{P}^*) > 0$, yields:
\begin{align*}\label{a833}
\mathcal{P}^t - \mathcal{P}^{t+1} &\leq \frac{\varphi(\mathcal{P}^t - \mathcal{P}^*) - \varphi(\mathcal{P}^{t+1} - \mathcal{P}^*)}{\varphi'(\mathcal{P}^t - \mathcal{P}^*)} \\
&\leq C_2 \|\mathbf{s}^t - \mathbf{s}^{t+1}\| \cdot [\varphi(\mathcal{P}^t - \mathcal{P}^*) - \varphi(\mathcal{P}^{t+1} - \mathcal{P}^*)]. \tag{G.33}
\end{align*}

Define for any non-negative integers $a, b$:
\[
\Delta_{a,b} := \varphi(\mathcal{P}^a - \mathcal{P}^*) - \varphi(\mathcal{P}^b - \mathcal{P}^*).
\]
Substituting equation (\ref{a833}) into (\ref{a830}) gives
\[
C_3 \|\mathbf{s}^{t+1} - \mathbf{s}^t\|^2 \leq \mathcal{P}^t - \mathcal{P}^{t+1} \leq C_2 \|\mathbf{s}^t - \mathbf{s}^{t+1}\| \Delta_{t,t+1}, \quad \forall t \geq t_2.
\]
Using the inequality $2\sqrt{ab} \leq a + b$ for all $a, b \geq 0$, we derive
\[
2 \|\mathbf{s}^{t+1} - \mathbf{s}^t\| \leq \|\mathbf{s}^t - \mathbf{s}^{t+1}\| + \frac{C_2}{C_3} \Delta_{t, t+1}. \tag{G.34}\label{a834}
\]

We now show that for any $t > t_2$,
\[
\sum_{k=t_2+1}^t \|\mathbf{s}^{k+1} - \mathbf{s}^k\| \leq \|\mathbf{s}^{t_2+1} - \mathbf{s}^{t_2}\| + \frac{C_2}{C_3} \Delta_{t_2, t+1}.
\]
Summing equation (\ref{a834}) with $k = t_2+1,t_2+2,\cdots,t$, we obtain
\begin{align*}
2 \sum_{k=t_2+1}^t \|\mathbf{s}^{k+1} - \mathbf{s}^k\| &\leq \sum_{k=t_2+1}^{t} \|\mathbf{s}^{k} - \mathbf{s}^{k-1}\| + \frac{C_2}{C_3} \sum_{k=t_2+1}^t \Delta_{k,k+1} \\
&\leq \sum_{k=t_2+1}^t \|\mathbf{s}^{k+1} - \mathbf{s}^{k}\| + \|\mathbf{s}^{t_2+1} - \mathbf{s}^{t_2}\| + \frac{C_2}{C_3} \Delta_{t_2+1,t+1}.
\end{align*}
The last inequality uses the fact that for any non-negative integers $a, b, c$, we have $\Delta_{a,b} + \Delta_{b,c} = \Delta_{a,c}$. Since $\varphi \geq 0$, for any $t > t_2$,
\[
\sum_{k=t_2+1}^t \|\mathbf{s}^{k+1} - \mathbf{s}^k\| \leq \|\mathbf{s}^{t_2+1} - \mathbf{s}^{t_2}\| + \frac{C_2}{C_3} \varphi(\mathcal{P}^{t_2+1} - \mathcal{P}^*).
\]
Since the right-hand side does not depend on $t$, the sequence $\{\mathbf{s}^t\}_{t \geq 1}$ has finite length:
\[
\sum_{t=1}^\infty \|\mathbf{s}^{t+1} - \mathbf{s}^t\| < \infty.
\]
Thus, the sequence $\{\mathbf{s}^t\}_{t \geq 1}$ is a Cauchy sequence in Euclidean space and hence convergent.

Combining the two cases above completes the proof of the first part of Theorem \ref{the2}. For the second part, Lemma A.8.4 guarantees the existence of a point $\mathbf{s}^* \in \mathrm{crit}\, \mathcal{L}$ such that $\lim_{t \to \infty} \mathbf{s}^t = \mathbf{s}^*$. Moreover, $\mathbf{s}^*$ is a Karush-Kuhn-Tucker point of problem (\ref{jin}).  
    
\end{proof}                      

\end{document}